\documentclass[a4paper,11pt,amstex]{amsart}
\usepackage{t1enc}

\usepackage{epsf}
\usepackage{amssymb}
\usepackage{latexsym}
\usepackage{amsmath}
\usepackage{eucal}
\usepackage{xypic,eepic,eepicemu}
\usepackage{palatino}
\xyoption{curve}
\input xy

\newtheorem{corollaire}{\underline{Corollaire}}

\newtheorem{Def}{\fbox{D\'efinition}}
\newtheorem{prop}{\fbox{Proposition}}
\newtheorem{lem}{\fbox{Lemme}}
\newtheorem{theo}{\fbox{Th\'eor\`eme}}

\newenvironment{ex}{\underline{Exemple}:}{$\Box$}

\newenvironment{pv}{\fbox{Preuve} \hspace{0.1cm}}{$\rule{2.5mm}{2.5mm}$}

\def\I{{\mathrm{I}}}
\def\S{{\mathrm{S}}}
\def\H{{\mathrm{H}}}

\def\C{{\mathbb C}}
\def\T{{\mathrm{T}}}
\def\Z{{\mathbb Z}}
\def\P{{\mathbf P}}

\def\Grass{{\mathrm{Grass}}}
\def\O{{\mathcal{O}}}
\def\M{{\mathcal{M}}}
\def\dim{{\mathrm{dim}\hspace{0.2cm}}}

\def\Hom{{\mathcal{H}\mathrm{om}}}
\def\Ext{{\mathcal{E}\mathrm{xt}}}

\def \ta{\tau}
\def \ta1{\tau_1}

\def \fc{\frac}

\def \ri{\longrightarrow}

\def \O{{\mathcal O}}

\title{Le degr\'e de la vari\'et\'e des courbes de Poncelet}
\author{Yann S\'epulcre}

\begin{document}

\maketitle

\begin{abstract}
 Nous calculons pour tout entier $n \geq 4$ le degr\'e de la vari\'et\'e des courbes
de Poncelet, qui est une sous-vari\'et\'e de dimension $2n +5$ de
l'espace projectif des courbes planes de degr\'e $n$. Pour cela nous
interpr\'etons tout d'abord ce degr\'e comme nombre d'intersection
sur un espace de modules de faisceaux semi-stables sur $\P_2$. Nous
faisons alors appel \`a des espaces de modules de syst\`emes
coh\'erents semi-stables pour diff\'erentes valeurs d'un param\`etre
rationnel positif, qui sont birationnels au premier, et qui sont
li\'es entre eux par une suite d'\'eclatements et de contractions.
\end{abstract}

\section{Introduction}

 La d\'efinition traditionnelle d'une courbe de Poncelet g\'en\'erale est
 la suivante:

 \begin{Def}
  \label{courbeponc}Soit $C$ une courbe plane lisse de degr\'e $n$. Alors $C$ est de
Poncelet si et seulement si il existe une conique lisse $E$ et un
polyg\^one non d\'eg\'en\'er\'e \`a $n+1$ c\^ot\'es tangents \`a $E$
dont les $\fc{n(n+1)}{2}$ sommets sont sur $C$.
  \end{Def}

 On a le th\'eor\`eme suivant d\^u \`a Poncelet et Darboux (cf [D]).

 \begin{theo}
 Soit $C$ une courbe de Poncelet au sens d\'efini ci-dessus, et soit $E$ une conique
fix\'ee v\'erifiant la propri\'et\'e ci-dessus pour un polyg\^one
donn\'e $\mathcal{P}$. Alors $\mathcal{P}$ appartient \`a un pinceau
de polyg\^ones dont l'\'el\'ement g\'en\'eral v\'erifie la m\^eme
propri\'et\'e pour $E$ et $C$.
 \end{theo}

  L'ensemble des courbes de Poncelet g\'en\'erales forme un sous-ensemble localement
ferm\'e $P_n$ dans l'espace projectif des courbes planes de degr\'e
$n$. Dans [T], l'auteur montre que si $E$ est lisse, et $\Theta$ est
le fibr\'e inversible de degr\'e $n+1$ sur $E$ il existe un
morphisme $\Grass(2, \H^0 (\Theta)) --> P_n$ birationnel sur son
image, constitu\'ee des courbes de Poncelet g\'en\'erales
v\'erifiant avec $E$ la propri\'et\'e de la d\'efinition
\ref{courbeponc}. Ceci permet de voir que $P_n$ est irr\'eductible
de dimension $2n +5$.

  \begin{Def}
  \label{courbesPoncelet}La vari\'et\'e des courbes de Poncelet est l'adh\'erence
$\overline{P_n}$ de la sous-vari\'et\'e localement ferm\'ee $P_n$
dans l'espace projectif des courbes planes de degr\'e $n$.
  \end{Def}

   On cherche \`a conna\^itre le degr\'e de $\overline{P_n}$,
   c'est-\`a-dire \`a
\'evaluer le nombre d'intersection $h^{2 n + 5} \cap
[\overline{P_n}]$, o\`u $h$ est la restriction \`a $\overline{P_n}$
de la classe du fibr\'e $\O(1)$ sur l'espace projectif des courbes
planes de degr\'e $n$. Ceci va \^etre fait en ramenant le calcul de
ce nombre \`a celui d'un nombre d'intersection sur l'espace de
modules $M_n$ des faisceaux coh\'erents sur $\P_2$, semi-stables de
rang 2 et de classes de Chern $(0,n)$.

    En effet, les courbes de Poncelet peuvent \^etre interpr\'et\'ees comme courbes
des droites de saut de tels faisceaux semi-stables $F$, v\'erifiant
la condition $h^0 (F(1)) \geq 2$. On appelle $\mathcal{P}_n
\subseteq M_n$ le lieu, ferm\'e et irr\'eductible, des classes de
tels faisceaux coh\'erents. On a un morphisme naturel
    $$\beta: \mathcal{P}_n \ri \overline{P_n}$$
    qui \`a la classe d'un faisceau associe sa courbe des droites de saut. Dans [T]
il est prouv\'e que pour $n \geq 5$ ce morphisme est birationnel.
Pour $n=4$ c'est encore vrai gr\^ace au travail fait dans [LT]. Si
l'on note $\mathcal{D}=\beta^* (h)$, ceci implique que
    $$ \mathrm{deg} \hspace{0.1cm} \overline{P_n} = c_1 (\mathcal{D})^{2 n + 5} \cap
[\mathcal{P}_n] \hspace{0.4cm}(*)$$
  Nous introduisons ensuite, pour chaque nombre rationnel $\alpha >0$, un espace de
modules de syst\`emes coh\'erents $\alpha-$semi-stables, que l'on
note $S_{\alpha,n}$ (l'entier $n$ \'etant fix\'e, on notera plus
simplement $S_{\alpha}$). Chaque $S_{\alpha}$ est birationnel \`a la
vari\'et\'e $\mathcal{P}_n$, et pour $\alpha >>0$ on a un morphisme
birationnel de $S_{\alpha}$ sur $\mathcal{P}_n$. On montre sinon
qu'il existe un ensemble fini de valeurs critiques de $\alpha$
telles qu'entre deux valeurs critiques $\alpha_i$ et $\alpha_{i+1}$
cons\'ecutives, la classe d'isomorphisme de $S_{\alpha}$ reste
constante. \\

   Pour chaque $\alpha$ non critique, on peut d\'efinir un fibr\'e inversible
$\mathcal{D}_{\alpha}$ sur $S_{\alpha}$ tel que si $\alpha >>0$, le
fibr\'e $\mathcal{D}_{\alpha}$ soit l'image r\'eciproque de
$\mathcal{D}$ par le morphisme birationnel $S_{\alpha,n} \ri
\mathcal{P}_n$. Le calcul de (*) se m\`ene donc en deux \'etapes:
   \begin{enumerate}
   \item \'evaluer pour chaque param\`etre $\alpha$ {\bf{critique}} la diff\'erence
$c_{1}(\mathcal{D}_{\alpha_{-}})^{2n+5} \cap [S_{\alpha_{-}}] -
c_{1}(\mathcal{D}_{\alpha_{+}})^{2n+5} \cap [S_{\alpha_{+}}]$, o\`u
$\alpha_{-}$ et $\alpha_{+}$ sont deux valeurs tr\`es proches de
$\alpha$ telles que $\alpha_{-} < \alpha < \alpha_{+}$;

   \item \'evaluer $c_{1}(\mathcal{D}_{\epsilon})^{2n+5} \cap [S_{\epsilon}]$ pour
$\epsilon >0$ et strictement inf\'erieur \`a la plus petite valeur
critique.
   \end{enumerate}

    Le premier point occupe une large partie de cet article, le dernier occupant la
derni\`ere section. On est amen\'e \`a mettre en relation les
espaces de modules $S_{\alpha_{-}}$ et $S_{\alpha_{+}}$: on montre
que si l'on \'eclate ces deux vari\'et\'es, \`a savoir
$S_{\alpha_{-}}$ le long du ferm\'e des classes de syst\`emes
coh\'erents $\alpha_{-}$-semi-stables mais non
$\alpha_{+}-$semi-stables, et $S_{\alpha_{+}}$ le long du ferm\'e
des classes de syst\`emes coh\'erents $\alpha_{+}-$semi-stables mais
non $\alpha_{-}-$semi-stables, on obtient une m\^eme vari\'et\'e
$\widetilde{S_{\alpha}}$ \`a isomorphisme pr\`es. On \'etudie
ensuite plus pr\'ecis\'ement le diagramme de morphismes projectifs
birationnels
    $$\xymatrix{
     &  \widetilde{S_{\alpha}} \ar[dl] \ar[dr] & \\
     S_{\alpha_{-}} \ar[dr] & & S_{\alpha_{+}} \ar[dl] \\
     & S_{\alpha} & \\
     } $$
    qui est commutatif mais non cart\'esien, ce qui ne rend pas si ais\'e l'usage de
la formule de projection pour le point (1). Ajoutons \`a cela que
les espaces de modules $S_{\alpha}$ sont en g\'en\'eral singuliers,
et ce m\^eme si $\alpha$ est non critique. Nous montrerons ainsi
qu'en g\'en\'eral le diviseur exceptionnel de
$\widetilde{S_{\alpha}}$ contient deux composantes
irr\'eductibles, chacunes avec la multiplicit\'e 1. \\

     Nous obtenons finalement le th\'eor\`eme:
    \begin{theo}
    \label{theoprincipal}Le degr\'e de la vari\'et\'e
    $\overline{P_n }$ des courbes
    de Poncelet de degr\'e $n$ est \'egal
\`a: $54$ lorsque $n=4$, $6867$ lorsque $n=5$, $618$ lorsque
$n=6$.
    \end{theo}
    Nous obtenons en fait des formules g\'en\'erales pour
    calculer deg $\overline{P_n }$, mais dont les expressions
    sont lourdes. Nous nous sommes donc restreints \`a $n \leq 6$
    pour les applications num\'eriques.

\section{Faisceaux coh\'erents semi-stables de Poncelet}

 On fait d'abord quelques rappels sur la notion de semi-stabilit\'e et les espaces
de modules de faisceaux coh\'erents semi-stables.

 \subsection{Faisceaux semi-stables}
 \label{prelim}
  Soit $X$ une surface projective lisse, munie d'un fibr\'e tr\`es ample
$\O_{X}(1)$. Soit $K(X)$ le groupe de Grothendieck des classes de
$\O_{X}-$modules coh\'erents; il est \'equip\'e de la forme
quadratique d'intersection de cycles (cf [F]), et aussi de la forme
quadratique enti\`ere $q$ d\'efinie par $\chi(u^2)$, dont la forme
polaire est not\'ee $<,>$. Si $F$ est un faisceau alg\'ebrique
coh\'erent sur $X$, on d\'esigne par $\chi_F$ son polyn\^ome de
Hilbert relatif \`a la polarisation donn\'ee. Si $Y$ est une section
hyperplane de $X$, on note $h$ la classe de $O_Y$, de sorte que le
polyn\^ome de Hilbert est donn\'e par la formule
  $$ \chi_{F} (m) = \fc{1}{2}<c,h^2> m (m+1) + <c,h> m + \chi_{F}(0) $$
  o\`u $c$ est la classe de $F$ dans $K(X)$, et $\chi_{F}(0)$ est la
caract\'eristique d'Euler-Poincar\'e de $F$. La formule ci-dessus
permet de d\'efinir pour {\bf{toute}} classe de Grothendieck $c$ un
polyn\^ome $\chi_c$ \`a coefficients rationnels. On peut associer un
tel polyn\^ome \`a tout \'el\'ement de $K(X)_{\mathbb{Q}} = K(X)
\otimes_{\mathbb{Z}} \mathbb{Q}$. Dans toute la suite l'ordre que
nous consid\'erons sur $\mathbb{Q}[X]$ est l'ordre lexicographique.
Pour deux polyn\^omes $P$ et $Q$ on écrit $P \geq Q$ si et seulement
si $P-Q$ est de terme dominant positif. Le cône des classes de
$K(X)$ de polynôme de Hilbert strictement positif est noté
$K(X)_{\mathbb{Q}}^{+}$.

  \begin{Def}
  Soit $F$ un faisceau alg\'ebrique coh\'erent sur $X$, dont le support
sch\'ematique est de dimension $d$ (c'est encore le plus grand
entier $d$ tel que $<c,h^d> \neq 0$); on d\'efinit la
{\bf{multiplicit\'e}} de $F$ comme le nombre entier d\'efini par
$r=<F , h^d>$. Nous appelons le polyn\^ome $\fc{\chi_{F}}{r}$ le
polyn\^ome de Hilbert {\bf{r\'eduit}} de $F$.
  \end{Def}
  Lorsque $d=2$, il s'agit du rang de $F$, mais nous aurons besoin de la
multiplicit\'e de faisceaux de support une courbe.

  \begin{Def}
  \label{faiscsemistab}Supposons $F$ de classe de Grothendieck $c$, de
multiplicit\'e $r$, et de support de dimension $d$; alors $F$ est
{\bf{(semi-)stable}} relativement \`a la polarisation $\O_{X}(1)$
si:
  \begin{enumerate}
  \item tout sous-faisceau coh\'erent non nul $F^{'} \subset F$ a un support de
dimension $d$,

  \item si $0 < r^{'} < r$ est la multiplicit\'e d'un tel sous-faisceau, on a
  $$ \fc{\chi_{F^{'}}}{r^{'}} < \fc{\chi_{F}}{r} \hspace{0.5cm}(\mathrm{resp.}
\hspace{0.2cm} \leq)$$

  \end{enumerate}

  \end{Def}

  Chaque faisceau semi-stable $F$ admet une filtration de Jordan-H\"older, i.e une
filtration $\{0\} =F_0 \subset F_1 \subset ... \subset F_k = F$ par
des sous-faiseaux coh\'erents tels que les gradu\'es soient stables
de m\^eme polyn\^ome de Hilbert r\'eduit que $F$. Le gradu\'e total
d'une telle filtration ne d\'epend que de $F$ et d\'etermine la
classe dite de S-\'equivalence de $F$. Lorsque $F$ est stable,
classe de S-\'equivalence et classe d'isomorphisme co\"incident.

  \begin{theo}
  Soit $c \in K(X)$, il existe un espace de modules grossier $M_X (c)$ pour les
faisceaux semi-stables de classe de Grothendieck $c$, qui est un
sch\'ema projectif, et dont les points sont les classes de
S-\'equivalence de faisceaux semi-stables.
  \end{theo}

  \begin{pv}
  Voir [G].
  \end{pv}
  \hspace{0.3cm} \\

  Venons en au cas o\`u $X=\P_2$, $\O_{X}(1) = \O_{\P_2}(1)$. Alors on a $K(X)
\simeq \Z^3$, l'isomorphisme faisant correspondre \`a la classe $u$
d'un faisceau coh\'erent sur $\P_2$ le triplet donn\'e par son rang
$r$, sa premi\`ere classe de Chern $c_1$, et sa caract\'eristique
d'Euler-Poincar\'e $\chi$. La forme quadratique enti\`ere $q$ a pour
expression: $q(u)=2 r \chi + c_{1}^2 - r^2$. On note $h$ la classe
d'une droite du plan. Soit $n \geq 4$ un entier.

      \begin{Def}
      On note $\M_n$ l'espace de modules des faisceaux semi-stables de rang 2 et de
classes de Chern $(0,n)$ sur $\P_2$, c'est-\`a-dire de classe $c = 2
- n h^2$.
      \end{Def}
      Dans [L], on trouve un expos\'e de propri\'et\'es g\'en\'erales de $\M_n$: ce
sont des vari\'et\'es projectives int\`egres de dimension $4 n - 3$,
localement factorielles, singuli\`eres si et seulement si $n$ est
pair. Pour tout faisceau coh\'erent semi-stable de classe $c=2 - n
h^2$ on a la notion de {\bf{droite de saut}}: c'est une droite $l$
du plan telle que $h^{1} (F(-1) |_{l}) \neq 0$, ou de mani\`ere
\'equivalente telle que $F|_{l}$ n'est pas isomorphe \`a $\O_{l}
\oplus \O_{l}$. L'ensemble des droites de saut de $F$ forme une
courbe projective de $\P_{2}^*$ de degr\'e $n$, dont l'\'equation
est le d\'eterminant du morphisme de fibr\'es vectoriels
      $$ \O_{\P_{2}^*}(-1) \otimes \H^1(F(-2)) \ri \O_{\P_{2}^*} \otimes
\H^{1}(F(-1)) $$
      Ceci motive la d\'efinition suivante.

      \begin{Def}
      Il existe un morphisme bien d\'efini de vari\'et\'es projectives
      $$\beta: \M_n \ri |\O_{\P_{2}^*}(n)|$$
      qui \`a la classe d'un faisceau semi-stable associe sa courbe des droites de
saut. On nomme $\beta$ le morphisme de Barth.
      \end{Def}

      Pour plus de d\'etails on renvoie \`a [B] et [L3] (dans cette premi\`ere
r\'ef\'erence l'auteur n'utilise cependant que l'ouvert
quasiprojectif des classes de faisceaux stables). Le th\'eor\`eme de
Le Potier-Tikhomirov \'etablit que pour $n \geq 4$ le morphisme
$\beta$ est birationnel sur son image (voir [LT]).

  \subsection{Faisceaux semi-stables de Poncelet}

   Ils font l'objet de la d\'efinition suivante.

   \begin{Def}
   Un faisceau semi-stable de classe de Grothendieck $c=2- n h^2$, est dit de
Poncelet si et seulement si $h^0 (F(1)) \geq 2$, o\`u de mani\`ere
\'equivalente $h^1 (F(1)) \geq n-4$.
   \end{Def}

  Faisons deux remarques: pour $n=4$, tout faisceau semi-stable est de Poncelet; de
plus pour tout $n \geq 4$ on a $h^0(F(1)) \leq 3$ (cf [LH]). La
condition figurant dans la d\'efinition pr\'ec\'edente est de nature
d\'eterminantielle, on en d\'eduit que l'ensemble des classes de
faisceaux de Poncelet forme un ferm\'e dans $\M_n$, qu'on note
$\mathcal{P}_n$.

  \begin{prop}
  Le sous-sch\'ema ferm\'e $\mathcal{P}_n$ de $\M_n$ est irr\'eductible et normal de
dimension $2 n + 5$.
  \end{prop}

   \begin{pv}
   Voir [LT] proposition 4.8. On y montre de plus que l'ouvert de $\mathcal{P}_n$
des classes de faisceaux stables est de Cohen-Macaulay.
   \end{pv}
   \hspace*{0.2cm} \\

   Le th\'eor\`eme suivant est fondamental pour la suite.
   \begin{theo}
   \label{toma}La restriction du morphisme de Barth $\beta$ \`a $\mathcal{P}_n$ est
birationnelle sur son image, qui est $\overline{P_n}$ (voir la
d\'efinition \ref{courbesPoncelet} dans l'Introduction).
   \end{theo}

   \begin{pv} Voir [T] pour $n \geq 5$. Pour $n=4$ voir le th\'eor\`eme de Le
Potier-Tikhomirov. \end{pv} \hspace{0.2cm} \\

   La preuve du th\'eor\`eme pr\'ec\'edent, ainsi que celle du th\'eor\`eme de Le
Potier-Tikhomirov, utilise la notion de syst\`eme coh\'erent qui
fait l'objet du paragraphe suivant.

    \subsection{Syst\`emes coh\'erents et faisceaux de Poncelet}
   \label{genesystcoh}
    Les syst\`emes coh\'erents permettent une description explicite de certains
espaces de modules de faisceaux semi-stables car ils correspondent
\`a un "d\'evissage" de ces faisceaux obtenus \`a partir de
sous-espaces de sections globales (moyennant torsion). Il existe
de nombreux travaux faisant usage de cette notion: citons [L1],
[L2], [He], [Th]. Cette derni\`ere r\'ef\'erence traite d'espaces
de modules de fibr\'es semi-stables sur les courbes, avec pour
horizon la formule de Verlinde; dans [He] on cherche \`a calculer
les nombres de Donaldson du plan projectif, ce qui est plus proche
de notre cadre. Cependant notre travail est bien distinct et
s'appuie sur l'utilisation d'autres syst\`emes coh\'erents. Nous
allons commencer par r\'esumer certains des r\'esultats
g\'en\'eraux \'etablis dans [He] sur les syst\`emes coh\'erents et
leurs espaces de modules. \\

    Supposons \`a nouveau que $(X,\O_{X}(1))$ soit une surface projective lisse et
polaris\'ee.

    \begin{Def}
    \label{defsystalgcoh}On appelle {\bf{syst\`eme alg\'ebrique}} un triplet
$(\Gamma, i, F)$ o\`u $F$ est un faisceau de $\O_{X}-$modules,
$\Gamma$ un espace vectoriel et $i$ une application lin\'eaire
$\Gamma \ri \H^0(F)$. Lorsque $i$ est injectif et $F$ est
coh\'erent, un tel syst\`eme alg\'ebrique s'appelle un
{\bf{syst\`eme coh\'erent}}, et on le note plus simplement $(\Gamma,
F)$; il est alors implicite que l'on s'est donn\'e un sous-espace
vectoriel $\Gamma \subseteq \H^0 (F)$. L'entier $k=
    \dim \Gamma$ est appel\'e le {\bf{rang}} du syst\`eme coh\'erent, l'entier $d$
\'egal \`a la dimension du support de $F$ est appel\'e sa
{\bf{dimension}}.
    \end{Def}

    On d\'efinit la notion de morphisme de syst\`emes alg\'ebriques: c'est la
donn\'ee d'un morphisme de faisceaux alg\'ebriques coh\'erents $f:
F^{'} \ri F^{''}$ et d'une application lin\'eaire $g: \Gamma^{'} \ri
\Gamma^{''}$ tels que le diagramme suivant soit commutatif
    $$\xymatrix{
    \Gamma^{'} \ar[r]^{g} \ar[d]^{i^{'}} & \Gamma^{''} \ar[d]^{i^{''}} \\
    \H^0(F^{'}) \ar[r]^{h^0 (f)} & \H^0 (F^{''}) \\
    } $$
    A partir de cette d\'efinition on voit assez facilement que les syst\`emes
alg\'ebriques forment une cat\'egorie ab\'elienne; on montre qu'elle
poss\`ede assez d'objets injectifs. La sous-cat\'egorie des
syst\`emes coh\'erents est
additive mais non ab\'elienne. \\

    On associe \`a un syst\`eme {\bf{coh\'erent}} $\Lambda= (\Gamma, F)$ et une
classe $\alpha \in K(X) \otimes_{\mathbb{Z}} \mathbb{Q}$ la classe
de Grothendieck \`a coefficients rationnels
    $$ c_{\alpha}(\Lambda) = \dim \Gamma \cdot \alpha + c(F)$$
  \begin{Def}
  Un syst\`eme coh\'erent $\Lambda=(\Gamma,F)$ de dimension $d$ tel que $F$ soit de
multiplicit\'e $r$ est dit {\bf{$\alpha$-(semi)-stable}} si $F$ n'a
pas de sous-faisceau coh\'erent non nul de support de dimension $<
d$ et si pour un tel sous-module $F^{'}$ de multiplicit\'e $0 <
r^{'} < r$ on a, en posant $\Gamma^{'}=\Gamma \cap \H^0 (F^{'})$ et
$\Lambda^{'}=(\Gamma^{'},F^{'})$, l'in\'egalit\'e
  $$ \fc{\chi_{c_{\alpha}(\Lambda^{'})}}{r^{'}} < \fc{\chi_{c_{\alpha}(\Lambda)}}{r}
\hspace{0.4cm}(\mathrm{resp.} \hspace{0.3cm}\leq) $$
  \end{Def}

  Etant donn\'e $\Lambda$ un syst\`eme coh\'erent $\alpha-$semi-stable, il lui est
associ\'e une filtration de Jordan-H\"older, dont les gradu\'es ont
le m\^eme polyn\^ome de Hilbert r\'eduit
$\fc{\chi_{c_{\alpha}(\Lambda)}}{r}$. Deux syst\`emes coh\'erents de
dimension $d$ et de multiplicit\'e $r$, et $\alpha-$semi-stables,
sont dits S-\'equivalents si leurs gradu\'es totaux de
Jordan-H\"older sont isomorphes.

  \begin{theo}
  \label{theoexistespacesmodules}Soient $k$ un entier positif et $c$ une classe dans $K(X)$. Il existe pour les
syst\`emes coh\'erents $\Lambda=(\Gamma, F)$, tels que $\dim
\Gamma=k$ et $c(F)=c$, qui sont $\alpha-$semi-stables, un espace de
modules {\bf{grossier}} $\mathcal{S}_{X}(c,k)$; c'est un sch\'ema
projectif dont les points sont les classes de S-\'equivalence de
syst\`emes coh\'erents $\alpha-$semi-stables.
  \end{theo}

  \begin{pv}
  Voir [He].
  \end{pv} \\

  Pour d\'efinir les espace de modules $\mathcal{S}_{X}(c,k)$,
  on a besoin de d\'efinir des familles
  {\bf{plates}} de syst\`emes coh\'erents. C'est ce que nous faisons plus loin
  au paragraphe \ref{famillesrappels}. Alors $\mathcal{S}_{X}(c,k)$
  repr\'esente en un sens faible le foncteur qui \`a un sch\'ema $S$ associe
  les familles plates sur $S$ de syst\`emes coh\'erents de type $(c,k)$. \\

  Revenons maintenant au cas $X=\P_2$, $\O_{X}(1)=\O_{\P_2}(1)$. Soit $n \geq 4$, et
soit $\alpha$ un nombre {\bf{rationnel}} strictement {\bf{positif}},
correspondant \`a une classe dans $K(\P_2)_{\mathbb{Q}}$ de rang 0,
premi\`ere classe de Chern 0, et de caract\'eristique d'Euler
$\alpha$ (qui sont les trois composantes de l'isomorphisme $K(\P_2)
\otimes_{\mathbb{Z}} \mathbb{Q} \simeq \mathbb{Q}^3$). On
consid\`ere des syst\`emes coh\'erents $\Lambda=(\Gamma, \Theta)$
tels que $c(\Theta)=2 h + n h^2$, $\dim \Gamma = 2$. Le faisceau
$\Theta$ \'etant de multiplicit\'e 2 et de caract\'eristique
d'Euler-Poincar\'e $n+2$, la condition de $\alpha-$semi-stabilit\'e
s\'ecrit encore
  $$ \begin{array}{llcc}
   \alpha \cdot \dim \Gamma^{'} + \chi(\Theta^{'}) & \leq \alpha +
\fc{\chi(\Theta)}{2} \\
    & \leq \alpha + \fc{n}{2} + 1 \\
     \end{array} $$
  pour tout sous-faisceau coh\'erent $\Theta^{'}$ de multiplicit\'e 1 (i.e de
support sch\'ematique une droite), avec $\Gamma^{'} = \Gamma \cap
\H^0 (\Theta^{'})$.

  \begin{lem}
    \label{prealableconstructespacemodules}Sous les hypoth\`eses pr\'ec\'edentes, si
$(\Gamma, \Theta)$ est $\alpha-$semi-stable, on a $\H^1 (\Theta)=0$,
et $\Theta$ est engendr\'e par ses sections.
    \end{lem}

    \begin{pv}
    On observe que si le faisceau coh\'erent $\Theta$ est semi-stable (au sens de la
d\'efinition \ref{faiscsemistab}), alors $\H^1 (\Theta(i))=0$ pour
$i \geq -1$. Le faisceau $\Theta$ est alors 0-r\'egulier (au sens de
Castelnuovo-Mumford) et engendr\'e par ses sections.
    Si $\Theta$ n'est pas semi-stable, il existe une suite exacte
    $$0 \ri \O_{l^{'}}(a) \ri \Theta \ri \O_{l^{''}}(b) \ri 0$$
    o\`u $a > b$ sont deux entiers et $l^{'},l^{''}$ sont deux droites. On a $a >
0$; supposons $b \leq 0$ alors $\chi(\O_{l^{''}}) \leq 1$, et
$\chi(\Theta) \leq a + 2$. De plus $\Gamma^{'} = \Gamma \cap
\H^0(\O_{l^{'}}(a))$ est de dimension $\geq 1$, et la condition de
$\alpha-$semi-stabilit\'e implique
    $$ \alpha + a + 1 \leq \alpha + 1 + \fc{a}{2}$$
    ce qui est absurde. Donc $b > 0$, et le r\'esultat est alors imm\'ediat.  \end{pv}

  \begin{lem}
   \label{prealable2}Avec les notations qui pr\'ec\`edent, pour $\alpha > \fc{n}{2}
- 1$, la condition de $\alpha-$semi-stabilit\'e est \'equivalente
\`a la condition suivante: pour tout sous-faisceau coh\'erent
$\Theta^{'} \subset \Theta$ de multiplicit\'e 1, on a $\dim
\Gamma^{'} \leq \fc{\dim \Gamma}{2} = 1$, et si on a \'egalit\'e
alors $\chi (\Theta^{'}) \leq \fc{n}{2} + 1$.
   \end{lem}

   \begin{pv}
   Supposons que $(\Gamma,\Theta)$ soit $\alpha-$semi-stable. On a pour tout
sous-faisceau $\Theta^{'}$ de multiplicit\'e 1 l'in\'egalit\'e
   $$ \alpha (\dim \Gamma^{'} - 1) + \chi(\Theta^{'}) \leq \fc{n}{2} + 1 $$
   Si on a $\alpha > \fc{n}{2} - 1$, et $\dim \Gamma^{'} = 2$, alors
$\chi(\Theta^{'}) \leq \fc{n}{2} + 1 - \alpha < 2$. On aurait donc
$\chi(\Theta^{'}) \leq 1$. Mais comme $\Theta^{'}$ est un faisceau
de la forme $\O_{l}(m)$, avec $m$ un entier, qui admet un pinceau de
sections globales, on a n\'ecessairement $\chi(\Theta^{'}) \geq 2$,
contradiction. On a donc n\'ecessairement $\dim \Gamma^{'} \leq 1$,
ce qui prouve une implication. R\'eciproquement, supposons que la
deuxi\`eme condition soit v\'erifi\'ee, alors on doit v\'erifier que
dans le cas o\`u $\Gamma^{'} = 0$ on a $\chi (\Theta^{'}) \leq
\fc{n}{2} + \alpha + 1$. Mais on a $\chi(\Theta^{'}) \leq n$, car
sinon $h^0 (\Theta^{'}) \geq n + 1$, et comme $h^0 (\Theta)=n + 2$
d'apr\`es le lemme pr\'ec\'edent, on aurait $\Gamma^{'} \neq 0$,
contradiction. D'o\`u l'\'enonc\'e.
\end{pv} \\

  Ci-dessous nous pr\'esentons les espaces de modules de syst\`emes coh\'erents
auxquels nous ferons appel.

  \begin{Def}
  Pour $\alpha > 0$ rationnel on pose $\mathcal{S}_{\alpha,n}$ l'espace de modules
des syst\`emes coh\'erents $\alpha-$semi-stables $\Lambda=(\Gamma,
\Theta)$ tels que $\dim \Gamma=2$ et $c(\Theta)=2 h + n h^2$.
  \end{Def}

  L'entier $n$ \'etant fix\'e, on note plus simplement $\mathcal{S}_{\alpha} =
\mathcal{S}_{\alpha,n}$ dans la suite. Le lien entre faisceaux de
Poncelet et syst\`emes coh\'erents $\alpha-$semi-stables est
expos\'e dans ce qui suit.

  Soit $F$ un faisceau coh\'erent de Poncelet, suppos\'e donc de rang 2 et de
classes de Chern $(0,n)$, et consid\'erons un sous-espace vectoriel
$W \subseteq \H^0(F(1))$ de dimension 2; il existe alors une suite
exacte courte dont la fl\`eche de gauche est le morphisme
d'\'evaluation:
   $$0 \ri W \otimes \O_{\P_2} \ri F(1) \ri \Omega \ri 0$$
   Le conoyau $\Omega$ est un faisceau coh\'erent de support une conique, de
caract\'eristique d'Euler $\chi = 4 - n$. Posons $\Theta =
\underline{\Ext}^{1}_{\O_{\P_2}} (\Omega, \O_{\P_2})$, alors $\Theta$
est un faisceau de classe de Grothendieck $2 h + n h^2$, n'admettant
pas de sous-faisceau coh\'erent de support fini, et posant $\Gamma$
l'espace vectoriel dual de $W$, on a une inclusion naturelle $\Gamma
\subset \H^0 (\Theta)$.

   \begin{prop}
   \label{construcPonc}Avec les notations qui pr\'ec\`edent, pour $\alpha >
\fc{n}{2} - 1$, il y a \'equivalence entre la propri\'et\'e que $F$
soit semi-stable et la propri\'et\'e que $(\Gamma,\Theta)$ soit
$\alpha-$semi-stable.
   \end{prop}

   \begin{pv}
   En utilisant le lemme \ref{prealable2}, on peut appliquer la preuve de [LT],
Th.4.3 (i).
   \end{pv} \\

   Comme nous l'avons vu, pour $\alpha > \fc{n}{2} - 1$ la propri\'et\'e de
$\alpha-$semi-stabilit\'e ne d\'epend pas d'un choix particulier de
$\alpha$, la classe d'isomorphisme de l'espace de modules
$\mathcal{S}_{\alpha}$ reste donc constante.

    \begin{theo}
    \label{morphismesystcohfaisc}Pour $\alpha > \fc{n-2}{2}$ il existe un morphisme
bien d\'efini
    $$ \pi_{\alpha}: \S_{\alpha} \ri \mathcal{P}_n $$
    qui \`a la classe de S-\'equivalence de $(\Gamma, \Theta)$ associe la classe du
faisceau semi-stable $F$ d\'efini par l'extension figurant dans la
proposition \ref{construcPonc}. De plus $\pi_{\alpha}$ est surjectif
et est un isomorphisme au dessus de l'ouvert de $\mathcal{P}_n$ des
classes de faisceaux $F$ tels que $h^0 (F(1)) =2$.
    \end{theo}

    \begin{pv}
    Voir [LT] Th.4.3 (ii).
    \end{pv} \\

  Les morphismes $\pi_{\alpha}$ sont compatibles avec les isomorphismes canoniques
entre les $\mathcal{S}_{\alpha}$. Notons les $\pi$. Il ressort du
th\'eor\`eme pr\'ec\'edent et du th\'eor\`eme \ref{toma}, que le
degr\'e de la vari\'et\'e des courbes de Poncelet est \'egal \`a
  $$ \deg \overline{P_n} = c_{1} (\gamma^* \O(1))^{2n+5} \cap [\mathcal{S}_{\alpha}]$$
  pour $\alpha > \fc{n-2}{2}$, le morphisme $\gamma$ \'etant la compos\'ee $\beta
\circ \pi$, et $\O(1)$ \'etant la restriction du fibr\'e tr\`es
ample canonique \`a $\overline{P_n}$.

  \section{Etude d'espaces de modules de syst\`emes coh\'erents}

  Soit $n \geq 4$ un entier fix\'e, et $\alpha > 0$ un nombre rationnel. On a
introduit \`a la section pr\'ec\'edente les espaces de modules
$\mathcal{S}_{\alpha} = \mathcal{S}_{\alpha, n}$ (on omettra
l'indice $n$) qui sont des sch\'emas projectifs. Dans cette section
nous rappelons quelques propri\'et\'es de base de ces sch\'emas.
Puis au paragraphe \ref{famillesrappels} nous d\'efinissons les
notions de familles, familles plates de syst\`emes alg\'ebriques et
coh\'erents, et groupes et faisceaux $\Ext$ qui nous permettront
ult\'erieurement de mener une \'etude plus approfondie. \\

Le premier fait g\'en\'eral et important est que quelque soit
$\alpha$, l'espace de modules $\mathcal{S}_{\alpha}$ est int\`egre
de dimension $2 n + 5$. Le second fait g\'en\'eral que nous allons
exposer est qu'il existe des valeurs critiques $\alpha_0 < \alpha_1
< ... < \alpha_r$ de $\alpha$ telles que sur chaque intervalle
$]\alpha_{i}, \alpha_{i+1} [$, la classe d'isomorphisme de
$\mathcal{S}_{\alpha}$ reste constante.
  En $\alpha = \alpha_i$, on prend $\alpha_{i-1} < \alpha_{-} < \alpha_{i}$ et
$\alpha_{i} < \alpha_{+} < \alpha_{i+1}$ deux valeurs proches de
$\alpha$, et on montre qu'on dispose de deux morphismes
  $$\xymatrix{
   \mathcal{S}_{\alpha_{-}} \ar[dr] & & \mathcal{S}_{\alpha_{+}} \ar[dl] \\
     & \mathcal{S}_{\alpha} & \\
  } $$
  qui sont birationnels, donc surjectifs. On montre ensuite que si l'on \'eclate
$\mathcal{S}_{\alpha_{-}}$ et $\mathcal{S}_{\alpha_{+}}$ le long de
ferm\'es que nous explicitons, on obtient une m\^eme vari\'et\'e
projective \`a isomorphisme pr\`es.

  \subsection{Rappels sur la construction des espaces de modules}
  \label{rappelconstruct}
   La construction des espaces de modules de syst\`emes coh\'erents est tr\`es
semblable \`a celle, utilisant la th\'eorie g\'eom\'etrique des
invariants, servant \`a construire les espaces de modules de
faisceaux semi-stables de classe donn\'ee (on renvoie \`a [G] et [L]
pour les d\'etails de cette construction).

    Nous r\'esumons ici une construction sp\'ecifique aux espaces de modules
$\mathcal{S}_{\alpha}$, un peu plus simple pour des raisons
num\'eriques. Soit la classe de Grothendieck $c= 2 h + n h^2$ et
soit $\Lambda = (\Gamma, \Theta)$ un syst\`eme coh\'erent
$\alpha-$semi-stable tel que $c(\Theta)=c$, $\dim \Gamma = 2$. On
prend un espace vectoriel $H$ de dimension $\dim \H^0(\Theta) =
n+2$. On pose $V = \Grass(2, H) \times \mathrm{Quot}^{c}(H \otimes
\O_{\P_2})$ qui est un sch\'ema projectif dont les points ferm\'es
param\`etrent les couples donn\'es par un faisceau alg\'ebrique
quotient $\Theta$ de $H \otimes \O_{\P_2}$ de classe de Grothendieck
$c$, et un sous-espace $\Gamma \subset H$ de dimension $2$. Il
existe une {\bf{action}} naturelle du groupe $GL(H)$ sur $V$.

    On consid\`ere l'ouvert $V^{ss} \subset V$ des points param\'etrant des
faisceaux quotients sans sous-faisceaux de support fini, tels que
$\H^1 (\Theta)=0$, tels que l'application $H \ri \H^0(\Theta)$ soit
un isomorphisme, et tels que le faisceau coh\'erent $(\Gamma,
\Theta)$ obtenu soit $\alpha-$semi-stable. De telles conditions
proviennent du lemme \ref{prealableconstructespacemodules}.

    \begin{theo}
    Sous les hypoth\`eses pr\'ec\'edentes, l'espace de modules
$\mathcal{S}_{\alpha}$ est un bon quotient de l'ouvert $V^{ss}$ par
l'action naturelle de $GL(H)$.
    \end{theo}

    \begin{pv}
    Voir [He]. Le principe est d\'etablir une \'equivalence entre la propri\'et\'e
de  $\alpha-$semi-stabilit\'e du syst\`eme coh\'erent et la
propri\'et\'e de semi-stabilit\'e pour l'action de $GL(H)$ relative
\`a une polarisation du sch\'ema projectif $V$ d\'ependante de
$\alpha$. On peut donc voir les sch\'emas $\mathcal{S}_{\alpha}$
comme une famille de bons quotients de $V$ variant avec une
polarisation.
    \end{pv}

   \subsection{Irr\'eductibilit\'e et normalit\'e des espaces de modules}

    Les espaces de modules sont int\`egres de la bonne dimension $2 n + 5$ d'apr\`es
le th\'eor\`eme suivant.

    \begin{theo}
    \label{normaliteSalpha}Les espaces de modules $\mathcal{S}_{\alpha}$ sont int\`egres et normaux de
dimension $2 n + 5$.
    \end{theo}

     \begin{pv} On peut reprendre mot pour mot la preuve de la proposition $4.5$ de
[LT], qui fonctionne gr\^ace aux observations faites dans la preuve
de la proposition \ref{prealableconstructespacemodules} (qui
interviennent au lemme 4.6 de la r\'ef\'erence).
     \end{pv}

  \subsection{Familles de syst\`emes coh\'erents et d\'eformations
  infinit\'esimales}

    \label{famillesrappels}Dans la section \ref{genesystcoh} il a \'et\'e dit que la cat\'egorie des
syst\`emes alg\'ebriques poss\`ede assez d'objets injectifs.

    \begin{Def}
    Etant donn\'es deux syst\`emes coh\'erents $\Lambda^{'}=(\Gamma^{'}, F^{'})$ et
$\Lambda^{''}=(\Gamma^{''}, F^{''})$, consid\'er\'es comme des
syst\`emes alg\'ebriques, on d\'efinit pour tout entier $q \geq 0$,
le groupe $\Ext^q (\Lambda^{'},\Lambda^{''})$ comme le $q$-i\`eme
foncteur d\'eriv\'e du foncteur exact \`a gauche $\Hom$ appliqu\'e \`a
$ \Lambda^{'}$ et $\Lambda^{''}$.
    \end{Def}

    Il est assez facile de voir, gr\^ace aux th\'eor\`emes de finitude cohomologique
de Serre, que ces espaces vectoriels sont de dimension finie.

   \begin{prop}
   On a une suite exacte
   $$ \begin{array}{llcc}
   0 \ri  \Hom(\Lambda^{'},\Lambda^{''}) \ri \Hom(F^{'}, F^{''}) \ri
\Hom(\Gamma^{'}, \H^0 (F^{''}) / \Gamma^{''}) \ri \\
   \hspace{1.2cm} \Ext^1 ( \Lambda^{'},\Lambda^{''}) \ri \Ext^1 (F^{'}, F^{''}) \ri
\Hom(\Gamma^{'}, \H^1 (F^{''})) \ri ... \\
   \end{array} $$
 \end{prop}

 \begin{pv}
 Voir [He].
 \end{pv} \\

  Etant donn\'es deux syst\`emes alg\'ebriques $\Lambda^{'}$ et $\Lambda^{''}$,
  l'espace vectoriel $\Ext^1 (\Lambda^{''} , \Lambda^{'})$ s'interpr\`ete, comme dans le cas
  des faisceaux coh\'erents, comme espace param\'etrant les classes
  d'extensions
$$ 0 \ri \Lambda^{'} \ri \Lambda \ri \Lambda^{''} \ri 0 $$
  Nous utiliserons ce point de vue \`a de nombreuses reprises. \\

 Venons en \`a la notion de famille de syst\`emes alg\'ebriques
param\'etr\'ee et plate sur une vari\'et\'e alg\'ebrique. Nous
supposons que $(X, \O_{X}(1))$ est une surface projective complexe
polaris\'ee, et que $S$ est une vari\'et\'e alg\'ebrique connexe.
Nous notons $p: S \times X \ri S$ la projection sur le premier
facteur pour toute vari\'et\'e $S$.

   \begin{Def}
   Une famille {\bf{param\'etr\'ee}} par (et {\bf{plate}} sur) $S$ de syst\`emes
alg\'ebriques est la donn\'ee d'un triplet $(\mathbf{\Gamma}, i,
\mathbf{F})$ tel que $\mathbf{F}$ soit un faisceau de $\O_{S \times
X}-$modules (qui soit $\O_{S}-$plat), tel que $ \mathbf{\Gamma}$
soit un faisceau de $\O_{S}-$modules (plat), et tel que $i:
\mathbf{\Gamma} \ri p_{*}(\mathbf{F})$ soit un morphisme de
faisceaux de $\O_{S}-$modules.

    Une famille {\bf{param\'etr\'ee}} par (et {\bf{plate}} sur) $S$ de syst\`emes
coh\'erents est une famille param\'etr\'ee (et plate) de syst\`emes
alg\'ebriques sur $S$ telle $\mathbf{\Gamma}$ et $\mathbf{F}$ soient
des faisceaux coh\'erents, et telle que pour tout point $s \in S$,
l'application naturelle $\mathbf{\Gamma}_s \ri \H^0 (\mathbf{F}_s )$
soit injective.
   \end{Def}

   On emploiera parfois dans la suite l'expression "{\bf{syst\`eme coh\'erent sur $S$}}"
   ou "{\bf{syst\`eme alg\'ebrique sur $S$}}" pour parler d'une
   famille de syst\`emes coh\'erents ou alg\'ebriques param\'etr\'ee
   par $S$ au sens pr\'ec\'edent.
   On remarque que dans la d\'efinition de platitude
   de syst\`eme coh\'erent sur $S$ le faisceau $\mathbf{\Gamma}$ est
   automatiquement localement
libre, et pour chaque point ferm\'e $s \in S$, le syst\`eme
alg\'ebrique induit $(\mathbf{\Gamma}_s, i_s, \mathbf{F}_s)$ est un
syst\`eme coh\'erent.

 \begin{Def}
 Soit $g: S^{'} \ri S$ un morphisme de sch\'emas, l'{\bf{image
 r\'eciproque}} $g^* \mathbf{\Lambda}$ d'une famille $
\mathbf{\Lambda}$ d'une famille de syst\`emes alg\'ebriques ou
coh\'erents param\'etr\'ee par $S$ d\'esigne la famille
param\'etr\'ee par $S^{'}$ d\'efinie par le triplet $(g^*
\mathbf{\Gamma}, g^* i, g^* \mathbf{F})$.
 \end{Def}

En particulier quand $g$ est un plongement on parle plus simplement
de famille {\bf{restreinte}} \`a $S^{'}$. Si $L$ est un fibr\'e
inversible sur $S$, la famille de syst\`emes alg\'ebriques not\'ee
$L \otimes \mathbf{\Lambda}$ est d\'efinie comme $(L \otimes
\mathbf{\Gamma}, L \otimes i, L \otimes \mathbf{F})$. \\

 On d\'egage facilement la notion de morphisme de familles de
syst\`emes alg\'ebriques, ainsi que la d\'efinition du faisceau
$\underline{\Hom}_{p}(,)$ de $\O_{S}-$modules qui \`a un ouvert $U
\subseteq S$ associe le $\O_{S}(U)-$module des morphismes de
familles restreintes \`a $U$. On montre \`a nouveau que la
cat\'egorie des familles de syst\`emes alg\'ebriques param\'etr\'ees
par $S$ est une cat\'egorie ab\'elienne avec assez d'injectifs, dans
laquelle se plonge la cat\'egorie additive des familles de
syst\`emes coh\'erents.

   \begin{Def}
   Etant donn\'ees deux familles $\mathbf{\Lambda}^{'}$ et $\mathbf{\Lambda}^{''}$
de syst\`emes coh\'erents param\'etr\'ees par $S$, pour tout entier
$i \geq 0$ le faisceau coh\'erent
$\underline{\Ext}^{i}_{p}(\mathbf{\Lambda}^{'},
\mathbf{\Lambda}^{''})$ sur $S$ est le $i-$\`eme foncteur d\'eriv\'e
du foncteur $\underline{\Hom}_{p}(,)$ appliqu\'e \`a
$\mathbf{\Lambda}^{'}$ et $\mathbf{\Lambda}^{''}$.
   \end{Def}

   Nous utiliserons ces faisceaux pour d\'ecrire des faisceaux conormaux et des
faisceaux tangents \`a des sous-sch\'emas d'espaces de modules de
syst\`emes coh\'erents. \\

 Comme rappel\'e juste apr\`es le th\'eor\`eme
\ref{theoexistespacesmodules}
    la notion de platitude d'une famille est une notion essentielle pour
    parler de la propri\'et\'e d'espace de modules grossier.
    Plus exactement, si $\alpha
\in K(X)_{\mathbb{Q}}^{+}$, et $\mathcal{S}_{\alpha}(c,k)$ est
l'espace de modules des syst\`emes coh\'erents $(\Gamma, \Theta)$
sur $X$, avec $c(\Theta)=c$ et $\Gamma$ de dimension $k$, alors la
propri\'et\'e de module grossier de $\mathcal{S}_{\alpha}(c,k)$
implique la propri\'et\'e suivante: pour tout sch\'ema $S$ de type
fini et toute famille $\mathbf{\Lambda}$ de syst\`emes coh\'erents
de type $(c,k)$, param\'etr\'ee par $S$ et plate sur $S$, il existe
un morphisme de sch\'emas
     $$ f_{S}: S \ri \mathcal{S}_{\alpha}(c,k)$$
     qui \`a un point $s$ de $S$ associe la classe du syst\`eme coh\'erent
param\'etr\'e en $s$. \\

 Le th\'eor\`eme suivant explicite l'espace tangent \`a l'espace de
 module en un point $\alpha-$stable, et donne une description
 locale de ce sch\'ema au voisinage d'un tel point.

 \begin{theo}
  \label{critlissite}Soit $\Lambda=(\Gamma, \Theta)$ un syst\`eme coh\'erent
  $\alpha-${\bf{stable}}
d\'efinissant un point $p$ de l'espace de modules
$\mathcal{S}_{\alpha}$. Alors l'espace tangent en $p$ est isomorphe
\`a l'espace vectoriel $\Ext^1 (\Lambda, \Lambda)$; de plus la
vari\'et\'e $\mathcal{S}_{\alpha}$ est lisse en ce point si et
seulement si $\Ext^2 (\Lambda, \Lambda)=0$.

 De plus au voisinage de $p$ le sch\'ema $(\mathcal{S}_{\alpha} , p)$
se plonge dans une vari\'et\'e lisse $(V,p)$ d'espace tangent de
Zariski $\Ext^1 (\Lambda , \Lambda)$, et l'id\'eal $J$ de ce
plongement est engendr\'e par dim $\Ext^2 (\Lambda , \Lambda)$
\'equations, de diff\'erentielles nulles en $p$.
  \end{theo}

  \begin{pv}
  Voir [He] Th. 3.12 et 3.13. On peut donner une preuve (tr\`es rapide) du
  premier \'enonc\'e
  en disant
  qu'on connait l'espace tangent au sch\'ema de Hilbert des quotients $V$
en un point au dessus de $\Lambda$ (voir paragraphe
\ref{rappelconstruct}), et qu'on a un crit\`ere connu de lissit\'e
de $V$ en un tel point. Il suffit ensuite d'utiliser que sur
l'ouvert de $V^{ss}$ des points repr\'esentant des classes de
syst\`emes coh\'erents stables, l'action de $GL(H)$ est
{\bf{libre}}. \end{pv} \\

 Pour tout syst\`eme coh\'erent $\mathbf{\Lambda}$ {\bf{plat}} sur $S$, il
 existe en tout point $s \in S$ un morphisme dit de
 {\bf{Kodaira-Spencer}}
$$ \omega_s : T_s S \ri \Ext^1 \left( \mathbf{\Lambda}_s ,
\mathbf{\Lambda}_s \right) $$
 La construction de $\omega_s$ est classique (voir [L1]): un vecteur tangent
 dans $T_s S$ s'interpr\`ete comme un morphisme $t : D \ri S$ o\`u
 $D = \mathrm{Spec} \left( \mathbf{C} [ \epsilon ] / \epsilon^2
 \right)$ envoyant le point ferm\'e de $D$ sur $s$; soit $I$
 l'id\'eal du point ferm\'e de $D$, de la suite exacte canonique
 $0 \ri I \ri \mathbf{C} [ \epsilon ] / \epsilon^2 \ri \mathbf{C}
 \ri 0$ on d\'eduit la suite exacte image r\'eciproque
 $$ 0 \ri \mathbf{\Lambda}_s \ri t^* (\mathbf{\Lambda})
 \ri \mathbf{\Lambda}_s \ri 0 $$
 A la lumi\`ere du th\'eor\`eme pr\'ec\'edent, lorsque $\mathbf{\Lambda}$
 param\`etre en $s$ un syst\`eme
 coh\'erent stable, le morphisme $\omega_s$ n'est autre que la
 diff\'erentielle du morphisme modulaire $f_S$ en $s$.

\subsection{Valeurs critiques du param\`etre}

    Revenons au cas o\`u $X=\P_2$, $n \geq 4$ est un entier et $\alpha > 0$ est un
nombre rationnel. Dans ce paragraphe nous abordons l'\'etude des
espaces de modules $\mathcal{S}_{\alpha}$, consid\'er\'ees comme
formant une famille - au sens ensembliste - de vari\'et\'es
d\'ependantes de $\alpha$.

    \begin{Def}
    \label{defcrit}On dit que $\alpha$ est {\bf{critique}} si et seulement si pour
tout $\beta > \alpha$ rationnel, il existe un syst\`eme coh\'erent
$\alpha-$semi-stable mais non $\beta-$semi-stable.
    \end{Def}

    On pourrait objecter qu'on peut d\'efinir une valeur critique \`a l'aide des
valeurs inf\'erieures $\beta < \alpha$. Ceci est en fait
\'equivalent, ainsi que la prochaine proposition le montre. Notons
d'abord que si $\alpha$ est \'egal \`a $\fc{n}{2} - k$, avec $1 \leq
k < \fc{n}{2}$ un entier, alors $\alpha$ est critique au sens de la
d\'efinition \ref{defcrit}. En effet un syst\`eme coh\'erent
$(\Gamma, \O_{l^{'}}(\fc{n}{2} - \alpha)) \oplus (0,
\O_{l^{''}}(\fc{n}{2} + \alpha))$, avec $\Gamma$ de dimension 2, est
$\alpha-$semi-stable mais n'est pas $\beta-$semi-stable pour $\beta
> \alpha$.

    \begin{prop}
    \label{Valeurcrit}Les valeurs critiques de $\alpha$ sont les rationnels de la forme $\fc{n}{2} -
k$, avec $1 \leq k < \fc{n}{2}$ un entier.
    \end{prop}

    \begin{pv}
    Si $\Lambda=(\Gamma, \Theta)$ n'est pas $\beta-$semi-stable il existe un
sous-faisceau $\Theta^{'} \subset \Theta$ de multiplicit\'e 1, tel
que en posant $\Gamma^{'} = \Gamma \cap \H^0 (\Theta^{'})$ on ait:
    $$ \beta (1 - \dim \Gamma^{'}) < \chi (\Theta^{'}) - \fc{n}{2} - 1 $$
    On a toujours $ \alpha (1 - \dim \Gamma^{'}) \geq \chi (\Theta^{'}) - \fc{n}{2}
- 1 $, de sorte que si $\beta > \alpha$ on a $\dim \Gamma^{'}=2$.
Choisissant $\beta$ tel que $\beta - \alpha < \fc{1}{2}$, on en
d\'eduit l'existence d'un sous-syst\`eme coh\'erent
$\Lambda^{'}=(\Gamma^{'},\Theta^{'})$ de $\Lambda$ (c'est-\`a-dire
un sous-objet dans la cat\'egorie ab\'elienne des syst\`emes
alg\'ebriques) tel que: $\dim \Gamma^{'}=2$ et $\chi(\Theta^{'}) =
\fc{n}{2} - \alpha + 1$. Ceci impose que $\fc{n}{2} - \alpha$ est
entier et $\geq 1$. Avec l'observation pr\'ec\'edant la proposition
on peut conclure.
    \end{pv}
    \hspace{0.2cm} \\

     Les syst\`emes coh\'erents $\alpha-$semi-stables mais $\beta$-instables
exhib\'es dans la preuve sont aussi instables pour un param\`etre $0 < \gamma
< \alpha$. Ceci montre que la d\'efinition \ref{defcrit} avec des
valeurs inf\'erieures du param\`etre donnent les m\^emes valeurs
critiques. Il r\'esulte de ce qui pr\'ec\`ede le corollaire suivant.

     \begin{corollaire}
     Le nombre rationnel $\alpha$ ne prend qu'un nombre fini de valeurs critiques
$\alpha_i = \fc{n}{2} - i$ pour $1 \leq i < \fc{n}{2}$ entier, et pour
$\alpha$ compris entre deux valeurs critiques cons\'ecutives
$\alpha_{i}$ et $\alpha_{i+1}$ la classe d'isomorphisme de la
vari\'et\'e $\mathcal{S}_{\alpha}$ reste constante.

     \end{corollaire}

   \begin{pv}
   En effet entre deux valeurs critiques cons\'ecutives la condition de
$\alpha-$semi-stabilit\'e est ind\'ependante de $\alpha$.
   \end{pv}
  \hspace{0.2cm} \\

   D'autre part si $\alpha_{-}, \alpha_{+}$ sont des nombres rationnels positifs
tels que $\alpha -1 < \alpha_{-} < \alpha < \alpha_{+} < \alpha +
1$, avec $\alpha$ critique, on montre facilement que tout syst\`eme
coh\'erent $\alpha_{-}$ ou $\alpha_{+}-$semi-stable est
$\alpha-$semi-stable. On a deux morphismes birationnels de
vari\'et\'es alg\'ebriques
     $$\xymatrix{
     \mathcal{S}_{\alpha_{-}} \ar[dr]_{\pi_{-}} & & \mathcal{S}_{\alpha_{+}}
     \ar[dl]^{\pi_{+}} \\
      & \mathcal{S}_{\alpha} & \\
      } $$
   qui ont des sections rationnelles $\mathcal{S}_{\alpha} -->
\mathcal{S}_{\alpha_{+/-}}$ dont le lieu d'ind\'etermination est
celui des points de $\mathcal{S}_{\alpha}$ correspondant aux classes
de S-\'equivalence $\Lambda^{'} \oplus \Lambda^{''}$ o\`u
$\Lambda^{'} = (\Gamma, \O_{l^{'}}(\fc{n}{2} - \alpha))$ et
$\Lambda^{''} = (0, \O_{l^{''}}(\fc{n}{2} + \alpha))$. C'est donc un
sous-sch\'ema ferm\'e de $\mathcal{S}_{\alpha}$.

   \begin{Def}
   \label{defSigma}On note $\Sigma^{\alpha}$ le sous-sch\'ema ferm\'e r\'eduit de
$\mathcal{S}_{\alpha}$ des classes de syst\`emes coh\'erents
$\alpha-$semi-stables mais non $\beta-$ stables pour $\beta >
\alpha$.
   \end{Def}

    On a alors la proposition suivante. On pose $\Sigma = \Sigma^{\alpha}$, et $D$
la vari\'et\'e d'incidence
     $$ \xymatrix{
     D \ar[r]^{q} \ar[d]_{p} & \P_2 \\
     \P_{2}^* & \\
     } $$
     dont les points param\`etrent les couples $(l,x)$ o\`u $l$ est une droite et
$x$ un point sur $l$; lorsque l'on a affaire au produit $\P_{2}^*
\times \P_{2}^*$, on note $D_1$ (resp. $ D_2$) les sous-vari\'et\'es
d'incidence contenues dans $\P_{2}^* \times \P_{2}^* \times \P_2$
relatives au premier (resp. second) facteur.

      \begin{Def}
      \label{grassmanrelative}On note $G$ la grassmanienne relative $\Grass(2, p_*
\O_{D}(\fc{n}{2} - \alpha)) \ri \P_{2}^*$, dont les points
param\`etrent les sous-espaces $\Gamma \subset \H^0
(\O_{l}(\fc{n}{2} - \alpha))$ de rang 2, o\`u $l$ est une droite de
$\P_{2}$.
      \end{Def}

       On peut d\'ecrire explicitement $\Sigma$. Remarquons que la vari\'et\'e $G
\times \P_{2}^*$ param\`etre une famille plate de syst\`emes
coh\'erents $\alpha-$semi-stables, \`a savoir la famille
$\mathbf{\Lambda}^{'} \oplus \mathbf{\Lambda}^{''}$ o\`u
$\mathbf{\Lambda}^{'} = (\mathbf{\Gamma}, \O_{D_1}(\fc{n}{2} -
\alpha))$, $\mathbf{\Lambda}^{''}=(0, \O_{D_2}(\fc{n}{2} + \alpha))$
et $\mathbf{\Gamma}$ est le sous-fibr\'e tautologique de rang 2 sur
$G$. D'apr\`es ce qui pr\'ec\`ede il existe un morphisme $\phi: G
\times \P_{2}^* \ri \mathcal{S}_{\alpha}$ dont l'image est $\Sigma$,
et ce morphisme est bijectif.

    \begin{prop}
     \label{Sigma}Le sous-sch\'ema $\Sigma$ est isomorphe \`a $G \times \P_{2}^*$.
    \end{prop}

    \begin{pv}
    Nous allons montrer que $\phi$ est un plongement. Soit $V$ le
sch\'ema produit introduit au paragraphe \ref{rappelconstruct}, dont
un ouvert $V^{ss}$ donne par quotient l'espace de modules
$\mathcal{S}_{\alpha}$. Le sch\'ema $V^{ss}$
     param\`etre une famille plate de syst\`emes coh\'erents $\mathbf{\Lambda}$
     qui sont $\alpha-$semi-stables. Soit $d$ la classe dans
     $K(\P_2)$ \'egale \`a
$h + (\fc{n}{2} + \alpha) h^2$, autrement dit celle des faisceaux
$\O_{l}(\fc{n}{2} + \alpha)$, et Quot$^{d}(\Lambda)$ le sch\'ema
dont les points param\`etrent les morphismes surjectifs de
syst\`emes coh\'erents
     $$ (\Gamma, \Theta) \ri (0, \O_{l}(\fc{n}{2} + \alpha)) $$
     o\`u $(\Gamma, \Theta)$ est $\alpha-$semi-stable.
     En d'autres termes le morphisme de
faisceaux $\Theta \ri  \O_{l}(\fc{n}{2} + \alpha)$ est surjectif, et
l'application induite $\Gamma \ri \H^0 ( \O_{l}(\fc{n}{2} +
\alpha))$ est nulle. Il existe un morphisme naturel $i:$
Quot$^{d}(\mathbf{\Lambda}) \ri V^{ss}$, et on a \'egalement une
action naturelle de $GL(H)$ sur Quot$^{d}(\mathbf{\Lambda})$ pour
laquelle $i$ est \'equivariant. Un point de
Quot$^{d}(\mathbf{\Lambda})$ param\`etre une suite exacte au sens
des syst\`emes alg\'ebriques
     $$0 \ri \Lambda^{'} \ri \Lambda \ri \Lambda^{''} \ri 0 \hspace{2cm}(+)$$
     avec $\Lambda=(\Gamma, \Theta)$, $\Lambda^{'}= (\Gamma, \O_{l^{'}}(\fc{n}{2} -
\alpha))$ et $\Lambda^{''} = (0, \O_{l^{''}}(\fc{n}{2} + \alpha) )
$, et globalement il existe une suite exacte de syst\`emes
coh\'erents sur Quot$^{d}(\mathbf{\Lambda})$ donn\'ee par
$$ 0 \ri \left( \mathbf{\Gamma} , \mathbf{\Theta}^{'} \right)
\ri i^* \left( \mathbf{\Lambda} \right) \ri \left( 0 ,
\mathbf{\Theta}^{''} \right) \ri 0 \hspace{2cm} (*)$$
 Maintenant la fibre du morphisme $i$ au
dessus d'un point de $V^{ss}$ qui param\`etre un syst\`eme
coh\'erent $\Lambda$ a pour espace tangent
$\Hom(\Lambda^{'},\Lambda^{''})$ avec les notations de (+). Il est
donc nul, la fibre est finie et r\'eduite, et $i$ est une immersion
ferm\'ee car \'etant donn\'e un syst\`eme coh\'erent $\Lambda =
(\Gamma, \Theta)$, il existe au plus une suite exacte du type (+)
\`a scalaire pr\`es.

  Consid\'erons \`a pr\'esent le
sch\'ema Quot$_{red}^{d}(\mathbf{\Lambda})$ muni de sa structure
{\bf{r\'eduite}}, il se plonge donc dans $V^{ss}$. Soit $F$ une
quelconque de ses composantes irr\'eductibles.
     On a un diagramme commutatif de vari\'et\'es alg\'ebriques (sauf
\'eventuellement le sch\'ema en haut \`a droite)
     $$\xymatrix{
     F \ar@{^{(}->}[r]^i \ar[d]_j &
     V^{ss} \ar[d]^{//} \\
     G \times \P_{2}^* \ar[r]^{\phi} & \mathcal{S}_{\alpha} \\
     } $$
     o\`u la fl\`eche verticale de gauche associe \`a un quotient de $\Lambda$ les
deux droites supports de $\O_{l^{'}}$ et $\O_{l^{''}}$, et le
sous-espace $\Gamma \subset \H^0 (\O_{l^{'}}(\fc{n}{2} - \alpha))$.
\\

 Pour justifier rigoureusement l'existence de $j$ on proc\`ede comme
suit: le faisceau $\mathbf{\Theta}^{''}$ dans (*) est par
hypoth\`ese plat sur $F$, de dimension relative 1. Il est facile de
voir qu'il existe une r\'esolution de $\mathbf{\Theta}^{''}$ \`a
deux termes du type
$$ 0 \ri H_1 \ri H_0  \ri \mathbf{\Theta}^{''} \ri 0 $$
 o\`u $H_0 , H_1$ sont des fibr\'es vectoriels sur
 $F$ de m\^eme rang.
 En effet si $p:$
 $F \times \P_2
 \ri F$ est la premi\`ere projection le
 faisceau coh\'erent $p_* (\mathbf{\Theta}^{''})$ est localement
 libre d'apr\`es le th\'eor\`eme des images directes de Grauert (voir [Ha]
 III.12.9 qui s'applique car $F$ est int\`egre), et on a un morphisme surjectif d'\'evaluation
 $p^* \left( p_* (\mathbf{\Theta}^{''}) \right) \ri
 \mathbf{\Theta}^{''}$ dont le noyau est localement libre d'apr\`es
 le m\^eme th\'eor\`eme.

 On peut
 donc prendre le d\'eterminant $\mathrm{det} ( \mathbf{\Theta}^{''} )$, qui
 est un fibr\'e inversible sur $F \times \P_2$ dont la restriction \`a chaque fibre au dessus d'un
 point de $F$ d\'efinit un fibr\'e inversible
 de classe constante dans Pic($\mathbf{P}_2$), qui est $\O(1)$. On
 proc\`ede de m\^eme pour $\mathbf{\Theta}^{'}$ et par hypoth\`ese
 le fibr\'e $\mathbf{\Gamma}$ s'injecte dans $p_*
 ( \mathbf{\Theta}^{'} )$. La propri\'et\'e universelle de la grassmanienne
 ach\`eve de justifier l'existence de $j$. \\

 Maintenant quand on passe au quotient par l'action de $GL(H)$,
 l'image de $F$ se plonge dans $\mathcal{S}_{\alpha}$
et se factorise par le morphisme $\phi: G
\times \P_{2}^* \ri \mathcal{S}_{\alpha}$. Comme on peut choisir
$F$ dominant $G \times \P_{2}^*$ on peut conclure. \end{pv} \\

 Il n'est en fait pas difficile de voir que le sch\'ema
 Quot$^{d}(\mathbf{\Lambda})$ de la preuve ci-dessus est irr\'eductible. En
 effet les fibres du morphisme de passage au quotient $// : V^{ss}
 \ri \mathcal{S}_{\alpha}$ au dessus des points de $\Sigma$ sont
 toutes isomorphes \`a un m\^eme quotient de $GL(H)$ qui est
 irr\'eductible.

 On poursuit notre \'etude du comportement des espaces de modules au
voisinage d'une valeur critique avec la d\'efinition suivante.

      \begin{Def}
      Soit $\alpha$ critique. On note $\Sigma^{\alpha}_{-}$ le ferm\'e (muni de sa
structure de sch\'ema {\bf{r\'eduit}}) des points de
$\mathcal{S}_{\alpha_{-}}$ d'image contenue dans $\Sigma^{\alpha}$
par le morphisme naturel
      $\mathcal{S}_{\alpha_{-}} \ri \mathcal{S}_{\alpha}$.
      \end{Def}

    On fixe la valeur critique $\alpha$, et on pose $\Sigma^{\alpha} = \Sigma$ et
$\Sigma_{-} = \Sigma^{\alpha}_{-}$. On s'int\'eresse donc \`a des
classes de syst\`emes coh\'erents $\alpha_{-}-$semi-stables
$\Lambda=(\Gamma, \Theta)$ admettant une filtration dont les
gradu\'es sont du type $\Lambda^{'}=(\Gamma, \O_{l^{'}}(\fc{n}{2} -
\alpha))$ et $\Lambda^{''}=(0,\O_{l^{''}}(\fc{n}{2} + \alpha))$.
L'hypoth\`ese de $\alpha_{-}-$semi-stabilit\'e  entra\^ine
l'existence d'une extension {\bf{non scind\'ee}} de syst\`emes
coh\'erents
    $$0 \ri \Lambda^{'} \ri \Lambda \ri \Lambda^{''} \ri 0 \hspace{0.4cm} (*)$$
    R\'eciproquement, une telle extension non scind\'ee d\'efinit un syst\`eme
coh\'erent $\Lambda$ qui est $\alpha_{-}-$semi-stable, et dont la
classe est dans $\Sigma_{-}$. Il est facile de voir que l'on
d\'ecrit ainsi tous les syst\`emes coh\'erents
$\alpha_{-}-$semi-stables mais $\alpha_{+}-$instables.
    Les droites $l^{'}, l^{''}$ \'etant fix\'ees, ainsi que le sous-espace $\Gamma$,
les extensions de type (*) sont param\'etr\'ees par l'espace
projectif $\P (\Ext^1 (\O_{l^{''}} (2 \alpha), \O_{l^{'}}))$.

    \begin{Def}
    \label{idealdiagonale}On note $\mathcal{J}$ le faisceaux d'id\'eaux du ferm\'e
image r\'eciproque de la sous-vari\'et\'e diagonale $\P_{2}^*
\subset \P_{2}^* \times \P_{2}^*$ par la projection naturelle
$\Sigma \ri \P_{2}^* \times \P_{2}^*$.
    \end{Def}

    On d\'esigne par $\mathbb{P}(E)$ la fibration au sens de Grothendieck,
param\'etrant les quotients de rang 1 des fibres d'un faisceau
coh\'erent $E$. Rappelons que ce sch\'ema (projectif sur sa base),
est d\'efini comme Proj($\oplus_{l \geq 0} \hspace{0.1cm} \S^l E$)
lorsque $E$ est localement libre. Dans le cas g\'en\'eral pour avoir
une interpr\'etation g\'eom\'etrique concr\`ete de $\mathbb{P}(E)$
il faut \'ecrire une pr\'esentation de $E$ par des faisceaux
localement libres.

    \begin{prop}
     \label{sigma-}Pour $\alpha \geq 1$ critique on a un isomorphisme canonique
    $$\Sigma_{-} \simeq \mathbb{P} (\mathcal{J}^{2 \alpha - 2} (2 \alpha - 1, 2
\alpha - 2))$$
    Pour $\alpha=1/2$ on a un isomorphisme
    canonique $\Sigma_{-} \simeq \Sigma$, de m\^eme que pour $\alpha = 1$.
    Pour $\alpha = 3/2$ le sch\'ema $\Sigma_{-}$ est isomorphe \`a
    l'\'eclat\'e de
    $\Sigma$ le long de l'id\'eal $\mathcal{J}$.
    \end{prop}

     \begin{pv}
     Soient \`a nouveau $\mathbf{\Lambda}^{'}$ et $\mathbf{\Lambda}^{''}$ les
familles plates param\'etrees par $\Sigma= G \times \P_{2}^*$
introduites avant l'\'enonc\'e de la proposition \ref{Sigma}. On
note de plus $\tau_{1}, \tau_{2}$ les faisceaux sur $\P_{2}^* \times
\P_{2}^*$ images r\'eciproques de $\O_{\P_{2}^*}(1)$ par la
premi\`ere et seconde projection.

     Le faisceau coh\'erent $\underline{\Ext}^{1}_{p} (\mathbf{\Lambda}^{''},
\mathbf{\Lambda}^{'})$ sur $\Sigma$ est isomorphe \`a
$\mathcal{E}=\underline{\Ext}^{1}_{p} (\O_{D_{2}}( 2 \alpha),
\O_{D_{1}})$, et ses sous-espaces de rang 1 param\`etrent les
classes d'extensions $(*)$. Posons $\mathcal{F} =
\underline{\Ext}^{1}_{p} (\O_{D_{1}}, \O_{D_{2}} ( 2 \alpha - 3))$,
on a $\mathcal{F}^{\vee} \simeq \mathcal{E}$. Sur le sch\'ema
$\mathbb{P} (\mathcal{F})$, le quotient $\pi^* (\mathcal{F}) \ri
\O_{\mathbb{P}(\mathcal{F})}(1)$, o\`u $\pi: \mathbb{P}
(\mathcal{F}) \ri \Sigma$ est la projection naturelle, donne une
section partout non nulle du faisceau $\pi^* (\mathcal{E}) \otimes
\O_{\mathbb{P}(\mathcal{F})}(1)$, et par cons\'equent une extension
canonique de familles de syst\`emes coh\'erents
     $$ 0 \ri (\O_{\mathbb{P}(\mathcal{F})}(1) \otimes \mathbf{\Gamma},
\O_{\mathbb{P}(\mathcal{F})}(1) \otimes \O_{D_{1}} (\fc{n}{2} -
\alpha)) \ri \mathbf{\Lambda}
     \ri (0, \O_{D_{2}} (\fc{n}{2} + \alpha)) \ri 0 $$
      La famille $\mathbf{\Lambda}$ \'etant plate sur sa base $\mathbb{P}
(\mathcal{F})$ on obtient un morphisme
     $$ \sigma: \mathbb{P} (\mathcal{F}) \ri S_{\alpha_{-}}$$
     qui est bijectif sur son image, qui ensemblistement est $\Sigma_{-}$. En effet
tout point de $\Sigma_{-}$ correspond \`a une unique extension de
type (*) \`a multiplication par un scalaire pr\`es. On doit donc
tout d'abord identifier le
faisceau $\mathcal{F}$. \\

     De la pr\'esentation du $\O_{\P_{2}^* \times \P_{2}}-$ module coh\'erent
$\O_{D}$ donne\'e par
     $$ 0 \ri \O_{\P_{2}^*}(-1) \boxtimes \O_{\P_{2}}(-1) \ri \O \ri \O_{D} \ri 0 $$
     on d\'eduit la suite exacte longue:
     $$ p_{*} \O_{D_{2}}(2 \alpha - 3) \ri \tau_{1} \otimes p_{*}  \O_{D_{2}}(2 \alpha
- 2) \ri \mathcal{F} \ri R^1 p_* \O_{D_{2}}(2 \alpha - 3) \ri 0 $$
     Lorsque $\alpha=1/2$ on a $\mathcal{F} \simeq \tau_{2}^*$. Dans ce cas dans le
diagramme commutatif
     $$ \xymatrix{
     \mathbb{P}(\mathcal{F}) \ar[r] \ar[dr] & \Sigma_{-} \ar[d] \\
      & \Sigma \\
      } $$
      la fl\`eche oblique est un isomorphisme, la fl\`eche horizontale est
surjective, on peut donc conclure. Quand $\alpha \geq 1$, on obtient
la pr\'esentation suivante de $\mathcal{F}$
     $$ \O \boxtimes \S^{2 \alpha - 3} Q \ri \tau_{1} \boxtimes \S^{2 \alpha - 2} Q \ri
\mathcal{F} \ri 0 \hspace{2.5cm} (1)$$
     o\`u $Q$ est le fibr\'e tautologique de rang 2 sur $\P_{2}^*$. La fl\`eche de
gauche dans la suite exacte ci-dessus est \`a torsion pr\`es la
multiplication par la section canonique d\'efinissant la diagonale
dans $\P_{2}^* \times \P_{2}^*$. On v\'erifie donc facilement que
     $$\mathcal{F} \simeq \mathcal{J}^{2 \alpha - 2} (2 \alpha - 1, 2 \alpha - 2)$$
      Il reste \`a montrer que le morphisme $\mathbb{P} (\mathcal{F}) \ri
\mathcal{S}_{\alpha_{-}}$ est un plongement.

\begin{lem}
 \label{P(F)}Le sch\'ema $\mathbb{P}(\mathcal{F})$ est {\bf{r\'eduit}}.
 \end{lem}

 \begin{pv} Comme on l'a vu ci-dessus $\mathbb{P}(\mathcal{F})$ un
 sous-sch\'ema ferm\'e du fibr\'e
 projectif $\mathbb{P} ( \tau_1 \boxtimes \S^{2 \alpha - 2} Q )$ sur
 $\P_{2}^* \times \P_{2}^*$. On recouvre $\P_{2}^* \times \P_{2}^*$ par des
  ouverts isomorphes \`a $\mathbf{C}^2 \times \mathbf{C}^2$ tels que
  la restriction de l'espace total du fibr\'e projectif
  $\mathbb{P} ( \tau_1 \boxtimes \S^{2 \alpha - 2} Q )$ \`a chaque ouvert
  du recouvrement soit isomorphe \`a $\mathbf{C}^2 \times \mathbf{C}^2 \times
  \mathbf{C}^{2 \alpha - 1}$.
   Il suffit donc de montrer que
  la restriction du sous-sch\'ema ferm\'e $\mathbb{P}(\mathcal{F})$
  \`a chacun de ces ouverts affines est d\'efini par un id\'eal d'un anneau de
  polyn\^omes \`a $ 2 \alpha + 3$ variables sans \'el\'ements nilpotents.
  D'apr\`es la suite exacte (1), on sait qu'il existe une base de $\mathbf{C}^{2 \alpha
  -1}$ form\'ee de mon\^omes telle que si $(f_0 , ..., f_{2 \alpha - 2})$
  d\'esigne le syst\`eme de coordonn\'ees affines correspondant
  les \'equations qui d\'efinissent le sous-sch\'ema soient
  $$ (x_0 - x_1 ) f_{k - 1} + (y_0 - y_1) f_{k} \hspace{0.2cm},
  \hspace{0.5cm}
  1 \leq k \leq 2 \alpha - 2 $$
  o\`u $(x_0 , y_0)$ et $(x_1 , y_1)$ sont les deux syst\`emes de coordonn\'ees
  affines sur chaque facteur $\mathbf{C}^2$. Comme les formes lin\'eaires
  $x_0 - x_1$ et $y_0 - y_1$ sur $\mathbf{C}^2 \times \mathbf{C}^2$ sont
  ind\'ependantes, il revient au m\^eme de montrer que l'alg\`ebre quotient
  $$ \mathcal{A} := \mathbf{C} \left[ x,y,f_0 ,...,f_{k} \right] /
  \left< x f_{i-1} - y f_i \hspace{0.1cm} , \hspace{0.1cm} 1 \leq i
  \leq k \right> $$
  o\`u $k = 2 \alpha - 2$ est sans \'el\'ement nilpotent. Pour cela
  on ordonne les mon\^omes en $(f_0 , ..., f_k)$ selon l'ordre
  lexicographique avec de plus la convention $f_0 < f_1 < ... < f_k$.
  Chaque \'el\'ement de $\mathcal{A}$ s'\'ecrit
  de mani\`ere unique sous la forme
 $$ \rho (x,y) \hspace{0.1cm} + \hspace{0.1cm} \sum_{ \begin{array}{ll}
                        m_0 + ... + m_{k-1} \geq 1 \\
                        m_k \geq 0 \\
                        \end{array} }
                        \sigma_{m_0 , ..., m_k} (y) \hspace{0.1cm}
                        f_{0}^{m_0} \cdot \cdot \cdot f_{k}^{m_k}
                        \hspace{0.1cm} + \hspace{0.1cm} \sum_{l \geq
                        1} \tau_l (x,y) f_{k}^l $$
                        o\`u les $\rho, \sigma, \tau$ sont des
                        polyn\^omes. Si l'\'el\'ement est non nul
                        il existe un terme non nul de plus
                        petit ordre dans cette \'ecriture; s'il
                        est de plus nilpotent le polyn\^ome en $y$
                        ou en $(x,y)$ devant ce mon\^ome de plus petit ordre
                        est lui-m\^eme nilpotent, ce qui est absurde car on a
                        une inclusion canonique $\mathbf{C} \left[ x, y \right]
                        \subset \mathcal{A}$. \end{pv} \\

                        On ach\`eve maintenant la preuve de la proposition
                        \ref{sigma-}. D'apr\`es ce qui pr\'ec\`ede on a le
                        diagramme commutatif
$$ \xymatrix{
   \mathbb{P} \left( \mathcal{F} \right) \ar[r]^{\sigma} \ar[d] &
   \mathcal{S}_{\alpha_{-}} \ar[d] \\
   \Sigma \ar@{^{(}->}[r] & \mathcal{S}_{\alpha} \\
   } $$
 o\`u la fl\`eche du bas est un plongement d'apr\`es la proposition
 \ref{Sigma}. Soit $x$ un point de $\mathbb{P} \left( \mathcal{F} \right)$, et
 param\'etrant une extension $u: 0 \ri \Lambda^{'} \ri
 \Lambda \ri \Lambda^{''} \ri 0$. Soit $v$ un vecteur tangent en $x$.
  Supposons que $d \sigma_x (v)=0$;
 alors l'image de $v$ par la diff\'erentielle de la projection
 $\mathbb{P} \left( \mathcal{F} \right) \ri \Sigma$ est \'egalement
 nulle. Mais le noyau de cette diff\'erentielle est l'espace tangent
 \`a la fibre de cette projection, qui est isomorphe au quotient
 $\Ext^1 ( \Lambda^{''} , \Lambda^{'} ) / (u)$. La diff\'erentielle
 de $\sigma$ restreinte \`a l'espace tangent \`a la fibre
 s'identifie \`a la projection canonique
 $$ \Ext^1 ( \Lambda^{''} , \Lambda^{'} ) / (u) \ri \Ext^1 (\Lambda ,
 \Lambda) $$
 d\'eduite de la suite exacte $u$. En effet au paragraphe \ref{famillesrappels}
 on a rappel\'e que la diff\'erentielle de $\sigma$ s'identifiait au
 morphisme de d\'eformation de Kodaira-Spencer, dont on a donn\'e la construction.
  On a d'autre part la suite exacte
 $$ \Hom (\Lambda^{''} , \Lambda^{''} ) \ri \Ext^1 (\Lambda^{''} ,
 \Lambda^{'} ) \ri \Ext^1 (\Lambda^{''} , \Lambda)$$
 o\`u l'image de la fl\`eche de gauche est $(u)$. Enfin il est clair
 que le morphisme $\Ext^1 (\Lambda^{''} , \Lambda) \ri \Ext^1 (
 \Lambda , \Lambda )$ est injectif, ce qui entraine finalement que
 $v$ est nul. Le morphisme $\sigma$ est donc injectif, non
 ramifi\'e, d'image r\'eduite. La proposition s'ensuit. \end{pv} \\

 Le sous-sch\'ema $\Sigma_{-}$ n'est en g\'en\'eral pas int\`egre, ce
qui justifie que l'on ne puisse le d\'ecrire \`a l'aide d'un fibr\'e
associ\'e \`a un fibr\'e vectoriel sur $\Sigma$. Plus exactement on
a la proposition suivante.

     \begin{prop}
     \label{compsigma-}Lorsque $\alpha > 0$ est valeur critique et inf\'erieure \`a $3/2$, le
sous-sch\'ema $\Sigma_{-}$ est une vari\'et\'e irr\'eductible de
dimension $n - 2 \alpha + 2$. Lorsque $\alpha \geq 2$, la
vari\'et\'e $\Sigma_{-}$ poss\`ede deux composantes irr\'eductibles,
que nous notons $A$ et $B$. La vari\'et\'e $A$ est de dimension $n -
2 \alpha + 2$ et isomorphe \`a l'\'eclat\'e de $\Sigma$ le long du
ferm\'e d'id\'eal $\mathcal{J}^{2 \alpha - 2}$, la vari\'et\'e $B$
est de dimension $n-2$ et isomorphe au changement de base du fibr\'e
projectif $\P(\S^{2 \alpha - 2} Q) \ri \P_{2}^*$ par la projection
naturelle $G \ri \P_{2}^*$.
     \end{prop}

     \begin{pv}
     Lorsque $\alpha \geq 1$, le morphisme canonique de faisceaux quasi-coh\'erents
d'alg\`ebres
     $\oplus \hspace{0.1cm} \S^* \mathcal{J}^{2 \alpha - 2} \ri \oplus
\hspace{0.1cm} \mathcal{J}^{* (2 \alpha - 2)}$, qui est surjectif,
donne une immersion de la vari\'et\'e $A$ de l'\'enonc\'e dans
$\Sigma_{-}$. L'image r\'eciproque du ferm\'e d'id\'eal
$\mathcal{J}$ par le morphisme naturel $\Sigma_{-} \ri \Sigma$ est
isomorphe au fibr\'e en espaces projectifs sur $G$ donn\'e par
     $\mathbb{P}(\mathcal{J}^{2 \alpha - 2} / \mathcal{J}^{2 \alpha - 1})$. Comme on
a $\mathcal{J}^{2 \alpha - 2} / \mathcal{J}^{2 \alpha - 1} \simeq
\S^{2 \alpha - 2} (\mathcal{J} / \mathcal{J}^2)$. Lorsque $\alpha =
3 / 2$, on voit que $\Sigma_{-}$ est int\`egre isomorphe \`a
l'\'eclat\'e de $G$ le long de l'id\'eal $\mathcal{J}$. Lorsque
$\alpha \geq 2$, on obtient les deux composantes irr\'eductibles de
l'\'enonc\'e.
     \end{pv}

      \begin{Def}
      Soit $\alpha$ une valeur critique. On note $\Sigma^{\alpha}_{+}$ le ferm\'e des
points de $\mathcal{S}_{\alpha_{+}}$ dont l'image par la projection
$\mathcal{S}_{\alpha_{+}} \ri \mathcal{S}_{\alpha}$ appartient au
ferm\'e $\Sigma^{\alpha}$.
      \end{Def}

       La valeur critique $\alpha$ \'etant fix\'ee, on note dans la suite
$\Sigma_{+} = \Sigma^{\alpha}_{+}$. La propri\'et\'e de
$\alpha_{+}-$semi-stabilit\'e permet d'affirmer que de tels
syst\`emes coh\'erents $\Lambda$ s'ins\`erent dans une suite exacte
{\bf{non scind\'ee}} de syst\`emes coh\'erents
       $$0 \ri \Lambda^{''} \ri \Lambda \ri \Lambda^{'} \ri 0 \hspace{0.4cm}(+)$$
       avec $\Lambda^{'} = (\Gamma, \O_{l^{'}}(\fc{n}{2} - \alpha))$, $\Lambda^{''}
= (0, \O_{l^{''}}(\fc{n}{2} + \alpha))$. R\'eciproquement, toute
extension non scind\'ee de ce type est $\alpha_{+}$-semi-stable et
d\'efinit un point de $\Sigma_{+}$. De plus on d\'ecrit ainsi tous
les syst\`emes coh\'erents $\alpha_{+}-$semi-stables mais
$\alpha_{-}-$instables. Les droites $l^{'},l^{''}$ \'etant fix\'ees,
ainsi que le sous-espace $\Gamma$, les extensions de type (+) sont
param\'etr\'ees par l'espace vectoriel $\P (\Ext^1 (\Lambda^{'},
\Lambda^{''}))$.

       \begin{prop}
       \label{sigma+} Soient $\mathbf{\Lambda}^{'}$et $\mathbf{\Lambda}^{''}$ les familles de
syst\`emes coh\'erents param\'etr\'ees par $\Sigma$ utilis\'ees \`a
la proposition \ref{Sigma}. Alors le ferm\'e $\Sigma_{+}$, muni de
sa structure r\'eduite, est isomorphe au fibr\'e projectif $\P
(\underline{\Ext}^{1}_{p} (\mathbf{\Lambda}^{'},
\mathbf{\Lambda}^{''}))$, qui est localement libre de rang $n + 2
\alpha + 2$. La vari\'et\'e $\Sigma_{+}$ est donc int\`egre de
dimension $ 2 n + 4$.
       \end{prop}

      \begin{pv}
      La d\'emonstration consiste d'abord \`a montrer que le faisceau
$\underline{\Ext}^{1}_{p} (\mathbf{\Lambda}^{'},
\mathbf{\Lambda}^{''})$ est localement libre sur $\Sigma$ de rang
$n + 2 \alpha + 3$. C'est une application de la prop. 1.5 de [He].
Puis on montre par des arguments connus qu'il existe un morphisme
canonique $\P (\underline{\Ext}^{1}_{p} (\mathbf{\Lambda}^{'},
\mathbf{\Lambda}^{''})) \ri \Sigma_{+}$ qui est un plongement.
      \end{pv}

   \section{fibr\'es d\'eterminants et passages de valeurs critiques}

    Dans cette section, nous posons les notations suivantes: $n\geq 4$ est un
entier, $\alpha > 0$ est un nombre rationnel, qui est une valeur
critique au sens de la d\'efinition \ref{defcrit}. On fixe deux
valeurs critiques $\alpha_{-}, \alpha_{+}$, proches de $\alpha$ et
telles que $\alpha_{-} < \alpha < \alpha_{+}$. On pose $\Sigma =
\Sigma^{\alpha}$, $\Sigma_{-} = \Sigma_{-}^{\alpha}$, $\Sigma_{+} =
\Sigma_{+}^{\alpha}$.

     Nous cherchons \`a calculer pour $\alpha > \fc{n}{2} - 1$, le nombre
d'intersection
     $$ c_{1}(\gamma^* \O(1))^{2 n + 5} \cap [\mathcal{S}_{\alpha}]$$
     o\`u $\gamma: \mathcal{S}_{\alpha} \ri |\O_{\P_{2}^*}(n)|$ est la compos\'ee du
morphisme naturel $\pi:  \mathcal{S}_{\alpha} \ri \mathcal{M}_{n}$
est du morphisme de Barth $\beta:  \mathcal{M}_{n} \ri
|\O_{\P_{2}^*}(n)|$.

     \subsection{Fibr\'es d\'eterminants sur les espaces de modules}
 \label{fibredet}On conserve les notations du paragraphe pr\'ec\'edent.

    \begin{Def}
     On nomme le fibr\'e inversible $\beta^* (\O(1))$ sur $\mathcal{M}_n$ le fibr\'e
{\bf{d\'eterminant}} de Donaldson et on le note $\mathcal{D}$.
     \end{Def}

     Dans [L3] l'auteur \'enonce la  construction de $\mathcal{D}$ en termes de
K-th\'eorie. Le principe est le suivant: prenons la classe $- d \in
K(\P_{2})$ d\'efinie comme celle du faisceau $\O_{l}(-1)$, o\`u $l$
est une droite. Une telle classe n'est pas choisie arbitrairement:
elle est orthogonale \`a la classe $c=2 - n h^2$ pour la forme
quadratique $<,>$ introduite \`a la section \ref{prelim}.
     Soit $\mathcal{F}$ un faisceau coh\'erent {\bf{universel}} sur $W \times
\P_{2}$, o\`u $W$ est une vari\'et\'e quasi-projective qui par
quotient par l'action d'un groupe lin\'eaire donne
$\mathcal{M}_{n}$. Par universel on entend plat sur $W$, et on
suppose que les faisceaux param\'etr\'es par $\mathcal{F}$ aux
points de $W$ d\'ecrivent toutes les classes de S-\'equivalence de
faisceaux semi-stables de classe $c$. La construction de $W$ et $\mathcal{F}$ est classique et fait intervenir un sch\'ema des quotients. \\

     Notons $p,q$ les projections de $W \times \P_{2}$ sur les premier et second
facteurs. Alors on pose sur $W$ le fibr\'e inversible $\det \left( R
p_{!} (\mathcal{F} \cdot q^* d) \right)$.
      La classe $q^* d$ est l'image r\'eciproque de $d$ par le morphisme plat $q$,
le produit $\mathcal{F} \cdot q^* d$ est le produit de K-th\'eorie
de la classe de $\mathcal{F}$ avec $q^* d$, et $p_{!}$ associe \`a
la classe d'un faisceau coh\'erent sur $W \times \P_2$ la somme
altern\'ee des classes de ses images directes sup\'erieures $R^i
p_*$. Enfin $det$ est le d\'eterminant. Il doit \^etre pr\'ecis\'e
que ces d\'efinitions sont licites car on montre que l'on travaille
toujours avec des sommes de classes de faisceaux admettant des
r\'esolutions localement libres {\bf{finies}}. Maintenant la
propri\'et\'e d'orthogonalit\'e de $c$ et $d$ entra\^ine que le
fibr\'e ainsi d\'efini se descend sur le quotient $\mathcal{M}_{n}$;
en effet le stabilisateur en un point d'orbite {\bf{ferm\'ee}}
op\`ere trivialement sur la fibre, et le lemme de
Kempf-Drezet-Narasimhan permet de conclure (voir [DN]). Un
th\'eor\`eme d\^u \`a Drezet (voir [L] ou [L3] Th 3.10) montre que
tout fibr\'e inversible sur $\mathcal{M}_n$ s'obtient ainsi \`a
l'aide d'une classe $d$ orthogonale \`a $c$. Ces fibr\'es
inversibles v\'erifient tous une propri\'et\'e universelle
identique, et que nous \'ecrivons pour $\mathcal{D}$.

      \begin{theo}
       Soit un sch\'ema $S$ et un faisceau coh\'erent $\mathcal{F}$ sur $S \times
\P_2$, plat sur $S$ et param\'etrant des faisceaux semi-stables de
classe $c=2 - n h^2$, et soit $f: S \ri \mathcal{M}_n$ le morphisme
modulaire associ\'e; alors le fibr\'e inversible sur $S$ d\'efini
par $det \left( R p_{!} (\mathcal{F} \cdot q^* d) \right)$ est image
r\'eciproque du fibr\'e d\'eterminant $\mathcal{D}$ par $f$.
      \end{theo}

      \begin{pv}
      Voir [L3] sect. 2.13. \end{pv} \\

       On se propose d'effectuer une construction identique pour les espaces de
modules de syst\`emes coh\'erents. Soit donc $\alpha > 0$ un nombre
rationnel, et soit un espace de param\`etres $V$ sur lequel agit un
groupe $GL(H)$ et dont le bon quotient d'un ouvert $V^{ss}$ donne
l'espace de modules $\mathcal{S}_{\alpha}$ (voir section
\ref{rappelconstruct}). La vari\'et\'e $V^{ss}$ param\`etre une
famille plate de syst\`emes coh\'erents
$\mathbf{\Lambda}=(\mathbf{\Gamma}, \mathbf{\Theta})$. Utilisant les
m\^emes notations que ci-dessus, on pose alors sur $V^{ss}$ le
fibr\'e inversible donn\'e par
       $$det \left( R p_{!} (\mathbf{\Theta} \cdot q^* h) \right) \otimes
       det \left( \mathbf{\Gamma}^* \right) \hspace{0.3cm}(+)$$
        Ce choix non plus n'est pas anodin: on a en effet $< \Theta , h > - \dim
\Gamma = 0$ pour tout syst\`eme coh\'erent $(\Gamma, \Theta$ tel
que $c(\Theta)= 2 h + n h^2$ et $\dim \Gamma=2$. On peut ainsi
remarquer que l'expression (+) ne d\'epend pas du choix de la
famille $\mathbf{\Lambda}$ en vertu de [He], lemme 3.5, et du fait
que le ferm\'e compl\'ementaire de l'ouvert des points
$\alpha-$stables est de codimension $\geq 2$.

          \begin{prop}
       Pour $\alpha$ {\bf{non critique}}, le fibr\'e inversible d\'efini par (+) se
descend sur l'espace de modules, i.e il existe un fibr\'e inversible
not\'e $\mathcal{D}_{\alpha}$ sur $\mathcal{S}_{\alpha}$ tel que son
image r\'eciproque sur $V^{ss}$ par le morphisme de passage au
quotient ait l'expression (+).

        De plus, pour $\alpha > \fc{n}{2} - 1$, on a la relation
$\mathcal{D}_{\alpha} = \pi^* \mathcal{D}$, o\`u $\pi:
\mathcal{S}_{\alpha} \ri \mathcal{M}_n$ est le morphisme d\'efini au
th\'eor\`eme \ref{morphismesystcohfaisc}.
       \end{prop}

       \begin{pv}
       Le stabilisateur du groupe $GL(H)$ est isomorphe \`a $\mathbb{C}^*$ en chaque
point de $V^{ss}$ d\'efinissant un quotient $\alpha-$stable. En un
tel point $\Lambda = (\Gamma,\Theta)$ l'action d'un scalaire non nul
$l$ sur la fibre du fibr\'e inversible d\'efini par (+) est donn\'ee
par la multiplication par $l^{\chi (\Theta \cdot h) - \dim \Gamma} =
1$. Si maintenant $\Lambda$ est un point de $V^{ss}$ non
$\alpha-$stable mais d'orbite ferm\'ee sous $GL(H)$, on sait qu'il
s'\'ecrit
       $\Lambda = \Lambda^{'} \oplus \Lambda^{''}$ avec $\Lambda^{'} = (\Gamma^{'},
\O_{l^{'}}(\fc{n}{2}))$ et $\Lambda^{''} = (\Gamma^{''},
\O_{l^{''}}(\fc{n}{2}))$; les espaces $\Gamma^{'}$ et $\Gamma^{''}$
sont de dimension 1. Alors le stabilisateur en un tel point est soit
isomorphe \`a $\mathbb{C}^* \times  \mathbb{C}^*$ si $l^{'} \neq
l^{''}$, soit isomorphe \`a $GL(\mathbb{C}^2)$ si $l^{'}=l^{''}$.
L'action d'un endomorphisme inversible $g$ de rang 2 est alors:
       $det(g)^{<\Theta, h> / 2} \cdot det(g)^{-1} = 1$. Ceci prouve la premi\`ere
partie de la proposition, avec le lemme de Kempf-Drezet-Narasimhan. L'hypoth\`ese
$\alpha$ non critique est ici essentielle car si $\alpha$ est
critique on peut v\'erifier que les assertions pr\'ec\'edentes sont
fausses aux points correspondant aux classes $(\Gamma,
\O_{l^{'}}(\fc{n}{2} - \alpha)) \oplus (0,
\O_{l^{''}}(\fc{n}{2} + \alpha))$. \\

        Pour la seconde partie, supposons donc $\alpha > \fc{n}{2} - 1$, et notons
que le morphisme d'\'evaluation $p^* \mathbf{\Gamma} \ri
\mathbf{\Theta}$ fournit une suite exacte de faisceaux coh\'erents
sur $V^{ss} \times \P_2$
        $$ 0 \ri p^* \mathbf{\Gamma} \ri \mathcal{F}(1) \ri \check{\mathbf{\Theta}} \ri 0$$
        o\`u par d\'efinition $\check{\mathbf{\Theta}} = \underline{\Ext}^{1}_{p}
(\mathbf{\Theta}, \O)$
        et $\mathcal{F}$ est un faisceau plat sur $V^{ss}$ d\'ecrivant les faisceaux de
Poncelet de classe $c$. D'apr\`es la propri\'et\'e universelle du
fibr\'e d'eterminant $\mathcal{D}$,  l'image r\'eciproque sur
$V^{ss}$ de $\mathcal{D}$ par le morphisme compos\'e $V^{ss} \ri
\mathcal{S}_{\alpha} \ri \mathcal{M}_n$ est donn\'ee par
l'expression $det \hspace{0.1cm} p_{!} (\mathcal{F} \cdot q^* d) $
donn\'ee ci-dessus. Ce fibr\'e est donc \'egal \`a $det
\hspace{0.1cm} p_{!} (\check{\mathbf{\Theta}}(-1) \cdot q^* d)
\otimes det \hspace{0.1cm} \mathbf{\Gamma}^*$, qui est encore \'egal
d'apr\`es le th\'eor\`eme de dualit\'e relative (voir [H2]) et la
d\'efinition de $\check{\mathbf{\Theta}}$ \`a l'expression de
l'\'enonc\'e de la proposition. On peut conclure par passage au
quotient. \end{pv} \\

        Comme le fibr\'e $\mathcal{D}$, le fibr\'e $\mathcal{D}_{\alpha}$ v\'erifie une
propri\'et\'e universelle qui s'\'enonce en termes de familles
plates de syst\`emes coh\'erents. Il est clair \'egalement que
lorsque $\alpha$ est critique, en dehors des ferm\'es $\Sigma_{-}$
et $\Sigma_{+}$, les fibr\'es $\mathcal{D}_{\alpha_{-}}$ et
$\mathcal{D}_{\alpha_{+}}$ se correspondent par l'isomorphisme
birationnel $\mathcal{S}_{\alpha_{-}} --> \mathcal{S}_{\alpha_{+}}$.

\subsection{Eclatements d'espaces de modules}

        Soit $\alpha >0$ une valeur critique, et soient $\alpha_{-} < \alpha < \alpha_{+}$
des valeurs non critiques du param\`etres proches de $\alpha$. On
pose $\mathcal{S}=\mathcal{S}_{\alpha}$, $\mathcal{S}_{+/-} =
\mathcal{S}_{\alpha_{+/-}}$. Dans ce paragraphe on peut trouver une
vari\'et\'e $\widetilde{\mathcal{S}}$ munie de deux projections sur
$\mathcal{S}_{-}$ et $ \mathcal{S}_{+}$ tel que l'on ait un
diagramme commutatif
        $$ \xymatrix{
         & \widetilde{\mathcal{S}} \ar[dl] \ar[dr] & \\
        \mathcal{S}_{-} \ar@{-->}[rr] & &  \mathcal{S}_{+} \\
        } $$
      dont les deux fl\`eches obliques sont des morphismes birationnels. Comme on
s'y attend, on va \'eclater les deux vari\'et\'es
$\mathcal{S}_{+/-}$ le long des ferm\'es $\Sigma_{+/-}$.

      \begin{theo}
      \label{eclatements}Pour tout $\alpha > 0$ les \'eclat\'es de $\mathcal{S}_{-}$
et $\mathcal{S}_{+}$ le long de $\Sigma_{-}$ et $\Sigma_{+}$
respectivement sont isomorphes. Notons $\widetilde{\mathcal{S}}$ la
vari\'et\'e projective obtenue. Pour $\alpha \leq 3/2$ le diagramme
      $$ \xymatrix{
         & \widetilde{\mathcal{S}} \ar[dl] \ar[dr] & \\
        \mathcal{S}_{-} \ar[dr] & &  \mathcal{S}_{+} \ar[dl] \\
        & \mathcal{S}  & \\
        } $$
        est cart\'esien.
      \end{theo}

     La d\'emonstration utilise un lemme.

     \begin{lem}
     \label{espacetangentsigma-}Soit $x$ un point de $\Sigma_{-}$, correspondant \`a la classe d'un syst\`eme
coh\'erent $\Lambda$ s'ins\'erant dans une suite exacte (unique \`a
scalaire pr\`es)
     $$ 0 \ri \Lambda^{'} \ri \Lambda \ri \Lambda^{''} \ri 0
     \hspace{2cm}(+)$$
     o\`u $\Lambda^{'} = (\Gamma, \O_{l^{'}}(\fc{n}{2} - \alpha))$, $\Lambda^{''} =
(0, \O_{l^{''}}(\fc{n}{2} + \alpha))$. On a une suite exacte
d'espaces vectoriels
     $$0 \ri T_{x} \Sigma_{-} \ri T_{x} \mathcal{S}_{-} \ri \Ext^1 (\Lambda^{'},
\Lambda^{''}) $$
     o\`u les deux premiers termes sont des espaces tangents.
     Lorsque $\alpha \leq 3/2$ la fl\`eche de droite est surjective.
     \end{lem}

     \begin{pv} On sait que le
     morphisme de d\'eformation infinit\'esimale de Kodaira-Spencer
     $\omega: T_{x} \mathcal{S}_{-} \ri
     \Ext^1 \left( \Lambda , \Lambda \right)$ est un isomorphisme.
      La fl\`eche de droite
     dans la suite exacte en dessous de (+) s'identifie
     via $\omega$ \`a la projection
     canonique $\Ext^1 (\Lambda , \Lambda) \ri \Ext^1 (\Lambda^{'} , \Lambda^{''})$
     d\'eduite de la suite exacte
     $0 \ri \Lambda^{'} \ri \Lambda \ri \Lambda^{''} \ri 0$.
 On note $\Ext^{1}_{+/-} \left( \Lambda , \Lambda \right)$
     les groupes $\Ext$ filtr\'es relativement \`a la filtration
     d'Harder-Narasimhan $(\Lambda^{'} , \Lambda^{''})$ de
     $\Lambda$ (cf [L] (chap.15.3 ou [DL]).
     Il existe en particulier une suite exacte longue
     d' $\Ext$ filtr\'es dont nous n'\'ecrivons que la s\'equence qui nous int\'eresse
     $$  0 \ri \Ext_{-}^1 \left( \Lambda , \Lambda \right) \ri
     \Ext^1 \left( \Lambda , \Lambda \right) \ri \Ext^{1}_{+} \left(
     \Lambda , \Lambda \right) \hspace{2cm}(*) $$
     La fl\`eche $\Ext_{-}^1 ( \Lambda , \Lambda ) \ri
     \Ext^1 \left( \Lambda , \Lambda \right)$ est en effet injective car
      le groupe $\Ext_{+}^0 ( \Lambda , \Lambda ) =
      \Hom ( \Lambda^{'} , \Lambda^{''} )$ est nul. Le lemme
      \ref{espacetangentsigma-} est alors cons\'equence du lemme qui
      suit. \end{pv}

 \begin{lem}
  La restriction du morphisme de Kodaira-Spencer $\omega$ au
  sous-espace $T_x \Sigma_{-}$ de $T_x \mathcal{S}_{-}$ induit un
  isomorphisme
  $T_x \Sigma_{-} \simeq \Ext_{-}^1 (\Lambda , \Lambda)$.
  \end{lem}

    \begin{pv}  Comme $\Ext_{-}^1$ est foncteur d\'eriv\'e du sous-foncteur du foncteur
      $\Hom$
      respectant la filtration, on sait qu'on a le diagramme commutatif suivant
      $$ \xymatrix{
        T_x \Sigma_{-} \ar@{^{(}->}[r] \ar[d]_{\omega} & T_x \mathcal{S}_{-}
        \ar[d]^{\omega} \\
        \Ext_{-}^1 ( \Lambda , \Lambda ) \ar@{^{(}->}[r] & \Ext^1 (
        \Lambda , \Lambda ) \\
        } $$
      Il existe d'autre part un isomorphisme canonique $\Ext^{1}_{+} \left(
     \Lambda , \Lambda \right) \simeq \Ext^1 ( \Lambda^{'} ,
     \Lambda^{''} )$ et la fl\`eche de droite dans (*) s'identifie \`a
     l'application lin\'eaire $T_{x} \mathcal{S}_{-} \ri \Ext^1 (\Lambda^{'},
\Lambda^{''}) $. La suite spectrale convergeant vers $\Ext_{-}^1 (
\Lambda , \Lambda )$ en degr\'e 1 comporte les quatres termes
$E_{1}^{0,0} = \Hom ( \Lambda^{'} , \Lambda^{'} ) \oplus \Hom
(\Lambda^{''} , \Lambda^{''})$, $E_{1}^{1,0} = \Ext^1 ( \Lambda^{''}
, \Lambda^{'} )$, $E_{1}^{0,1} = \Ext^1 ( \Lambda^{'} , \Lambda^{'}
) \oplus \Ext^1 ( \Lambda^{''} , \Lambda^{''} )$, $E_{1}^{1,1} =
\Ext^2 ( \Lambda^{''} , \Lambda^{'} )$. Les autres termes sont nuls.
On a donc une suite exacte
$$ 0 \ri E_{2}^{1,0} \ri \Ext_{-}^1 ( \Lambda , \Lambda ) \ri
   E_{2}^{0,1} \ri 0 $$
   o\`u $E_{2}^{1,0}$ est le conoyau de l'application
   $d^{0,0}: \Hom(\Lambda^{'},\Lambda^{'}) \oplus \Hom(\Lambda^{''}
   ,  \Lambda^{''} ) \ri \Ext^1 (\Lambda^{''} , \Lambda^{'})$, qui
   est isomorphe \`a l'espace vectoriel quotient $\Ext^1 (\Lambda^{''} ,
   \Lambda^{'} ) / (u)$, o\`u $(u)$ est la classe de l'extension
   (+). L'espace vectoriel $E_{2}^{0,1}$ est un sous-espace de
   $\Ext^1 (\Lambda^{'} , \Lambda^{'}) \oplus
   \Ext^1 (\Lambda^{''} , \Lambda^{''})$, qui lui-m\^eme s'identifie
   \`a $T_{ \pi_{-} (x) } \Sigma$. La compos\'ee (2) suivante
   $$ \xymatrix{
   T_{x} \Sigma_{-} \ar[r] \ar@/_/[dr]^{(2)} &
      \Ext_{-}^1 (\Lambda , \Lambda) \ar[d]^{(1)} \\
      &  \Ext^1 (\Lambda^{'} ,
      \Lambda^{'} ) \oplus
   \Ext^1 (\Lambda^{''} , \Lambda^{''}) \\
 } $$
   s'identifie \`a l'application lin\'eaire
   tangente $\left( d \pi_{-} \right)_x : T_x \Sigma_{-} \ri
   T_{\pi_{-} (x) } \Sigma$. Pour v\'erifier que $\omega: T_{x} \Sigma_{-}
   \ri \Ext_{-}^1 (\Lambda , \Lambda)$ est un isomorphisme, il faut
   donc v\'erifier que les fl\`eches (1) et (2) dans le diagramme ci-dessus
   ont des noyaux isomorphes. Mais le noyau de $\left( d \pi_{-}
   \right)_x$ est l'espace tangent \`a la fibre de $\pi_{-}$ en $x$,
   isomorphe \`a l'espace projectif $\mathbf{P} \left( \Ext^1
   (\Lambda^{''} , \Lambda^{'} ) \right)$; son espace tangent en $x$
   est isomorphe \`a $\Ext^1 (\Lambda^{''} ,
   \Lambda^{'} ) / (u)$. Le lemme est montr\'e. \end{pv} \\

     \begin{pv} (du th\'eor\`eme \ref{eclatements}) On va faire la
     d\'emonstration quand les vari\'et\'es $\mathcal{S}_{-}$ et
     $\mathcal{S}_{+}$ param\`etrent des familles universelles de syst\`emes
     coh\'erents, la preuve dans le cas g\'en\'eral n'\'etant
     conceptuellement gu\`ere plus compliqu\'ee. Soit
     $\mathbf{\Lambda}_{-}$ la famille param\'etr\'ee par
     $\mathcal{S}_{-}$. On a une suite exacte de familles de
     syst\`emes coh\'erents param\'etr\'ees par $\Sigma_{-}$
     $$ 0 \ri \mathbf{\Lambda}^{'} \ri \mathbf{\Lambda}_{-}
     |_{\Sigma_{-}} \ri \mathbf{\Lambda}^{''} \ri 0 $$
 o\`u $\mathbf{\Lambda}^{'}$ param\`etre les syst\`emes coh\'erents
 de la forme $\left( \Gamma, \O_l \left( \fc{n}{2} - \alpha \right)
 \right)$ et $\mathbf{\Lambda}^{''}$ les syst\`emes coh\'erents de
 la forme $\left( 0 , \O_l \left( \fc{n}{2} + \alpha \right)
 \right)$. Soit $\widetilde{\mathcal{S}_{-}}$ la vari\'et\'e projective
 \'eclat\'ee de $\mathcal{S}_{-}$ le long du ferm\'e $\Sigma_{-}$.
 On note $p_{-}: \widetilde{\mathcal{S}_{-}} \ri \mathcal{S}_{-}$
 la projection canonique. On note $E_{-} = p_{-}^{-1} (\Sigma_{-})$
 le diviseur exceptionnel. Posons $\widetilde{\mathbf{\Lambda}_{-}}$,
 $\widetilde{\mathbf{\Lambda}^{'}} =
 p_{-}^* \left( \mathbf{\Lambda}^{'} \right)$, $\widetilde{\mathbf{\Lambda}^{''}} =
 p_{-}^* \left( \mathbf{\Lambda}^{''} \right)$. Ce sont des familles de syst\`emes
 coh\'erents param\'etr\'ees par $E_{-}$. On d\'efinit le syst\`eme alg\'ebrique
 $\overline{\mathbf{\Lambda}}_{-}$ sur $\widetilde{\mathcal{S}_{-}}$ comme
 noyau de la projection canonique $\widetilde{\mathbf{\Lambda}_{-}}
 \ri \widetilde{\mathbf{\Lambda}^{''}}$. On a alors un diagramme commutatif
$$ \xymatrix{
  & 0 \ar[d] & 0 \ar[d] &  & \\
  & \widetilde{\mathbf{\Lambda}_{-}} \left( - E_{-} \right)
  \ar@{=}[r] \ar[d] & \widetilde{\mathbf{\Lambda}_{-}} \left( - E_{-}
  \right) \ar[d] & & \\
  0 \ar[r] & \overline{\mathbf{\Lambda}}_{-} \ar[d] \ar[r] &
  \widetilde{\mathbf{\Lambda}_{-}} \ar[r] \ar[d] &
  \widetilde{\mathbf{\Lambda}^{''}} \ar@{=}[d] \ar[r] & 0 \\
  0 \ar[r] & \widetilde{\mathbf{\Lambda}^{'}} \ar[r] &
  \widetilde{\mathbf{\Lambda}_{-}} |_{E_{-}} \ar[r] &
  \widetilde{\mathbf{\Lambda}^{''}} \ar[r] & 0 \\
  } $$
 Montrons que $\overline{\mathbf{\Lambda}}_{-}$ est une famille de syst\`emes
 coh\'erents $\alpha_{+}-$semi-stables, plate sur $\mathcal{S}_{-}$. Pour tout
 point $x \notin E_{-}$, on a \'evidemment
 $\left( \overline{\mathbf{\Lambda}}_{-}
 \right)_{x} \simeq \left( \widetilde{\mathbf{\Lambda}_{-}} \right)_x $
 qui est un syst\`eme coh\'erent $\alpha_{+}-$semi-stable. Soit maintenant
 $x \in E_{-}$, $y = p_{-} (x)$, et $\left[ \Lambda \right]$ la classe
 param\'etr\'ee en $y$, ainsi que $\left( \Lambda^{'} , \Lambda^{''} \right)$
 la filtration d'Harder-Narasimhan de $\Lambda$ pour le param\`etre
 $\alpha_{+}$. Soit $v$ un vecteur tangent \`a
 $\widetilde{\mathcal{S}_{-}}$ en $x$ et transverse au diviseur $E_{-}$,
  ce qui implique l'image de $v$ par la diff\'erentielle $\left( d p_{-} \right)_{x}$
  est non contenue dans l'image de l'inclusion d'espaces tangents
  $\mathrm{T}_y \Sigma_{-} \subset \mathrm{T}_y \mathcal{S}_{-}$, c'est \`a dire
  que l'image de $\left( d p_{-} \right)_{x} (v)$ par la projection
  $\Ext^1 \left( \Lambda , \Lambda \right) \ri
  \Ext^1 \left( \Lambda^{'} , \Lambda^{''} \right)$ est non nulle
  (cf lemme \ref{espacetangentsigma-}). Le vecteur $v$ correspond \`a
  un morphisme $t: D \ri \widetilde{\mathcal{S}_{-}}$ o\`u
  $D = \mathrm{Spec} \left( \mathbf{C} [ \epsilon ] /
  \epsilon^2 \right)$. Soit $I = \mathbf{C} \epsilon$ l'id\'eal du point ferm\'e
  de $D$. De la suite exacte $0 \ri I \ri \mathbf{C} [ \epsilon ] /
  \epsilon^2 \ri \mathbf{C} \ri 0$ on d\'eduit la suite exacte
  $$ 0 \ri I \otimes \Lambda \ri t^* \left( \widetilde{\mathbf{\Lambda}_{-}} \right)
  \ri \Lambda \ri 0 $$
  o\`u $\Lambda$ est le syst\`eme coh\'erent param\'etr\'e en $y$.
  L'extension $u \in \Ext^1 (\Lambda , \Lambda \otimes I)$ obtenue n'est autre que
  l'image $\left( d p_{-} \right)_x (v)$ (cf paragraphe \ref{famillesrappels}).
  On applique maintenant le lemme du serpent au diagramme
$$ \xymatrix{
   0 \ar[r] & I \otimes t^* \left( \overline{\mathbf{\Lambda}}_{-}
   \right) \ar[d] \ar[r] &
   I \otimes t^* \left( \widetilde{\mathbf{\Lambda}_{-}} \right) \ar[r] \ar[d] &
   I \otimes t^* \left( \widetilde{\mathbf{\Lambda}^{''}} \right) \ar[r] \ar[d]^{0} & 0
   \\
   0 \ar[r] & t^* \left( \overline{\mathbf{\Lambda}}_{-} \right) \ar[r] &
   t^* \left( \widetilde{\mathbf{\Lambda}_{-}} \right) \ar[r] &
   t^* \left( \widetilde{\mathbf{\Lambda}^{''}} \right) \ar[r] & 0 \\
   } $$
   o\`u les deux lignes sont exactes et la derni\`ere fl\`eche \`a
   droite est nulle, car $t^* \left(
   \widetilde{\mathbf{\Lambda}^{''}} \right)$ a pour support le
   point ferm\'e de $D$. On obtient une suite exacte de syst\`emes
   coh\'erents
   $$ \xymatrix{
   0 \ar[r] & I \otimes \Lambda^{''} \ar[r] &
   \left( \overline{\mathbf{\Lambda}}_{-} \right)_x \ar[dr] \ar[rr]
   & & \left( \mathbf{\Lambda}_{-} \right)_y \ar[r] & \Lambda^{''}
   \ar[r] & 0 \\
   & & & \Lambda^{'} \ar[ur] & &  \\
   } $$
   o\`u la fl\`eche horizontale du milieu se factorise comme
   indiqu\'e par l'injection $\Lambda^{'} \ri \left(
   \mathbf{\Lambda}_{-} \right)_y$. La suite exacte obtenue
   $$ 0 \ri I \otimes \Lambda^{''} \ri
   \left( \overline{\mathbf{\Lambda}}_{-} \right)_x \ri \Lambda^{'}
   \ri 0 \hspace{2cm}(*)$$
 est par construction l'image de la classe $v \in \Ext^1 (\Lambda ,
 \Lambda)$ par la projection $\Ext^1 ( \Lambda , \Lambda ) \ri
 \Ext^1 ( \Lambda^{'} , \Lambda^{''} )$. L'extension (*) est non
 scind\'ee, c'est donc que $\left( \overline{\mathbf{\Lambda}}_{-}
 \right)_x$ est $\alpha_{+}$-semi-stable; il est clair ensuite
 d'apr\`es la caract\'erisation de la platitude par la constance du
 polyn\^ome de Hilbert sur les fibres d'un morphisme projectif que
 la famille de syst\`emes coh\'erents $\overline{\mathbf{\Lambda}}_{-}$ est
 plate sur $\widetilde{\mathcal{S}_{-}}$. Cela donne par la
 propri\'et\'e de module un morphisme $p_{+} : \widetilde{\mathcal{S}_{-}} \ri
 \mathcal{S}_{+}$. Les deux morphismes $\pi_{+} \circ p_{+}$ et
 $\pi_{-} \circ p_{-}$ envoient un point ferm\'e de
 $\widetilde{\mathcal{S}_{-}}$ sur le m\^eme point ferm\'e de
 $\mathcal{S}$; en effet avec les hypoth\`eses et notations qui pr\'ec\`edent
 un point $x$ de $E_{-}$ est envoy\'e sur le point de $\Sigma$
 qui est la classe du syst\`eme coh\'erent $\Lambda^{'} \oplus \Lambda^{''}$.
  On d\'eduit que $\pi_{+} \circ p_{+}$ et
 $\pi_{-} \circ p_{-}$ sont deux morphismes \'egaux car la vari\'et\'e
 $\widetilde{\mathcal{S}_{-}}$ est int\`egre.

  On a donc en
 particulier l'\'egalit\'e sch\'ematique
 $p_{-}^{-1} \left( \pi_{-}^{-1} \left( \Sigma \right) \right) =
  p_{+}^{-1} \left( \pi_{+}^{-1} \left( \Sigma \right) \right)$.
  Comme $\pi_{-}^{-1} \left( \Sigma \right) = \Sigma_{-}$ et
  $\pi_{+}^{-1} \left( \Sigma \right) = \Sigma_{+}$ au sens
  sch\'ematique et $p_{-}^{-1} \left( \Sigma_{-} \right)$ est
  le diviseur exceptionnel $E_{-}$, on en d\'eduit par la
  propri\'et\'e universelle de l'\'eclatement que le morphisme
  $p_{+}: \widetilde{\mathcal{S}_{-}} \ri \mathcal{S}_{+}$ se
  factorise par l'\'eclat\'e $\widetilde{\mathcal{S}_{+}}$ de
  $\mathcal{S}_{+}$ le long de $\Sigma_{+}$.

   En effectuant un raisonnement similaire \`a ce qui pr\'ec\`ede
   avec l\'eclat\'e $\widetilde{\mathcal{S}_{+}}$ en place de
   $\widetilde{\mathcal{S}_{-}}$, on d\'eduit qu'il existe
   finalement deux morphismes $f_{-}: \widetilde{\mathcal{S}_{-}} \ri
   \widetilde{\mathcal{S}_{+}}$ et $f_{+} : \widetilde{\mathcal{S}_{+}}
   \ri \widetilde{\mathcal{S}_{-}}$ qui sont inverses
   l'un de l'autre au sens ensembliste. Ils sont sont donc inverses
   au sens sch\'ematique car on manipule des vari\'et\'es
   int\`egres sur un corps alg\'ebriquement clos.

    Pour $\alpha \leq 3/2$ les ferm\'es $\Sigma_{-}$, $\Sigma_{+}$ sont lisses
    et irr\'eductibles et dans le lemme \ref{espacetangentsigma-}
    la suite est aussi exacte \`a droite, ce qui montre que le diagramme
    de l'\'enonc\'e est cart\'esien. Le th\'eor\`eme est prouv\'e. \end{pv} \\

 Il est important pour la suite de conna\^itre les composantes
 irr\'eductibles du diviseur exceptionnel commun $E$. D'apr\`es la
 proposition \ref{compsigma-}, lorsque $\alpha \geq 2$ ce
 ferm\'e est r\'eunion de deux composantes $A$ et $B$, il est
 certain que dans ce dernier cas $E$ contiendra au moins deux
 composantes irr\'eductibles. On est amen\'e \`a \'etudier
 l'intersection $\Sigma_{-}^{'}$ de $\Sigma_{-}$ avec le lieu singulier de
 $\mathcal{S}_{-}$, ce qui fait l'objet de ce qui suit.

 \begin{Def}
 Soit $E$ un fibr\'e vectoriel de rang 2 sur une vari\'et\'e
 alg\'ebrique $X$ suppos\'ee lisse. Pour tout couple $(i,d)$ d'entiers
 positifs tels que $1 \leq i \leq \fc{d + 1}{2}$, on note
 Sec$_{i}^{d}(E)$ la projection sur $\mathbb{P}(\S^d E)$
 du lieu des points de $\mathbb{P}(\S^d E) \times X$ o\`u le morphisme
 naturel de
 fibr\'es dit de "contraction"
 $$ \O_{\mathbb{P}(\S^d E)}(-1) \boxtimes \S^{i} E \ri
 \O_{\mathbb{P}(\S^d E)} \boxtimes \S^{d-i} E^* $$
 soit de rang $\leq i$.
 \end{Def}

 On peut \'etendre cette d\'efinition \`a $i=0$, en remarquant que
 Sec$_{0}^{d}(E) = \emptyset$. Il s'agit bien sur de l'adaptation \`a un
 contexte relatif de la
d\'efinition classique
 de {\bf{$i-$\`eme vari\'et\'e de s\'ecantes}} \`a la courbe rationnelle normale
 de  degr\'e $d$ (voir par exemple [Ha] prop. 9.7).
 On montre
 facilement que Sec$_{i}^{d}(E)$ est lisse irr\'eductible de
 dimension relative $2 i - 1$ au dessus de $X$. \\

  Revenons maintenant \`a l'\'etude de $\Sigma_{-}^{'}$. On suppose
  $\alpha \geq 2$. Alors la composante $B$ est isomorphe au fibr\'e
  projectif $\mathbb{P} (\S^{2 \alpha - 2} (\mathcal{J} /
  \mathcal{J}^2))$, o\'u $\mathcal{J}$ est le faisceau d'id\'eaux
  introduit \`a la d\'efinition \ref{idealdiagonale}. On peut encore
  \'ecrire ce fibr\'e comme changement de base du fibr\'e
  $\P (\S^{2 \alpha - 2} Q)$ sur $\P_{2}^*$, o\`u $Q$ est le fibr\'e
  tautologique quotient.

  \begin{prop}
  \label{fermesing}On suppose $\alpha \geq 2$.
  Le ferm\'e $\Sigma_{-}^{'}$ est non vide si et seulement
  si $\alpha \geq 5/2$, et dans ce cas est contenu dans
  la composante $B$ et est isomorphe \`a la vari\'et\'e des
  s\'ecantes relatives Sec$_{[\alpha - \fc{3}{2}]}^{2 \alpha - 2} Q \subset
  \P (\S^{2 \alpha - 2} Q)$. Il est donc int\`egre de dimension
  $n-4$ si $n$ est impair, et $n-5$ si n est pair.
  \end{prop}

 \begin{pv} Soit $\Lambda = (\Gamma, \Theta)$ un syst\`eme coh\'erent
 $\alpha_{-}-$stable d\'efinissant un point singulier de $\mathcal{S}_{-}$.
 Alors le th\'eor\`eme \ref{critlissite} \'enonce que $\Ext^2 (\Lambda, \Lambda) \neq 0$
 et donc que $\Ext^2 (\Theta, \Theta)^{*} = \Hom (\Theta, \Theta(-3)) \neq 0$.
 Par suite le faisceau $\Theta$ ne peut-\^etre semi-stable, et il admet une filtration
 d'Harder-Narasimhan
 $$ 0 \ri \O_{l^{'}}(a) \ri \Theta \ri \O_{l^{''}}(b) \ri 0 $$
 o\`u $l^{'}, l^{''}$ sont deux droites et $a > b$. On a de plus $l^{'} =
 l^{''}$ car il existe un morphisme non trivial de $u: \Theta \ri \Theta(-3)$.
Un tel morphisme $u$ se factorise n\'ecessairement comme compos\'ee
$$ u: \Theta \ri \O_{l}(b) \stackrel{v}{\ri} \O_{l}(a - 3) \ri
\Theta(-3) $$ On en tire la condition $a \geq b + 3$, et donc $a
\geq \fc{n + 3}{2}$. Maintenant posons $\Gamma^{'} = \Gamma \cap
\H^0 (\O_{l}(a))$, la condition de $\alpha-$semi-stabilit\'e donne
$a \leq \fc{n}{2} + \alpha (1 - \dim \Gamma^{'})$, ce qui entra\^ine
$\fc{3}{2} \leq \alpha (1 - \dim \Gamma^{'})$, et par suite
$\Gamma^{'}=0$. On a donc une suite exacte de syst\`emes coh\'erents
$$ 0 \ri (0, \O_{l} (a)) \ri \Lambda \ri (\Gamma, \O_{l}(b)) \ri 0
$$
avec $\fc{n + 3}{2} \leq a \leq \fc{n}{2} + \alpha$, $b \geq
\fc{n}{2} - \alpha$, et $a + b =n$. De plus si la classe de
$\Lambda$ appartient \`a $\Sigma_{-}$, elle appartient \`a $B$ car
le support de $\Theta$ est une conique double. Supposons que cela
soit v\'erifi\'e alors on a de plus une suite exacte
$$ 0 \ri (\Gamma, \O_{l}(\fc{n}{2} - \alpha) \ri (\Gamma, \Theta)
\ri (0, \O_{l}(\fc{n}{2} + \alpha)) \ri 0 \hspace{0.6cm}(*)$$ et
l'on obtient n\'ecessairement une inclusion de faisceaux $i :
\O_{l}(a) \subset \O_{l}(\fc{n}{2} + \alpha))$ qui se rel\`eve comme
suit
$$ \xymatrix{
   \Theta \ar[r] & \O_{l}(\fc{n}{2} + \alpha) \\
   & \O_{l}(a) \ar[ul]^{j} \ar[u]^{i} \\
    } $$
   L'inclusion $i$ d\'etermine une section non nulle de
   $\O_{l}(\fc{n}{2} + \alpha - a)$. Que l'inclusion $i$ admette un
   rel\`evement $j$ est \'equivalent au fait que le cup produit de cette
   section par la classe de (*) soit nul. Notons que l'espace vectoriel
   $\H^0 (\O_{l} (\fc{n}{2} + \alpha - a))$ s'identifie \`a la fibre
   de $\S^{\fc{n}{2} + \alpha - a} Q$ en $\check{l} \in \P_{2}^*$ et que
   $\Ext^1 (\O_{l}(a), \O_{l}(\fc{n}{2} - \alpha))$ s'identifie \`a la fibre
   de $\S^{2 \alpha - 2 - (\fc{n}{2} + \alpha - a)} Q^*$ au m\^eme point.
   Le ferm\'e $\Sigma_{-}^{'}$ est donc la r\'eunion croissante des ferm\'es
   $B_{i} =
\mathrm{Sec}_{i}^{2 \alpha - 2} (Q)$ pour $i$ compris entre 0 et
$[\alpha - \fc{3}{2}]$. \end{pv} \\

 Le ferm\'e $\Sigma_{-}^{'}$ est donc muni d'une filtration par les
 ferm\'es $B_{i}$ introduits dans la preuve pr\'ec\'edente. Pour
 $i=1$ on remarque que $B_{1} = \mathbb{P}(Q)$, le plongement $B_{1}
 \subset B = \mathbb{P}(\S^{2 \alpha - 2} Q)$ \'etant celui de
 Veronese. Notons que $B_{1}$ est aussi l'intersection de $A$ et
 $B$.

 \begin{prop}
 Le diviseur exceptionnel de $\widetilde{\mathcal{S}}$ poss\`ede au
 plus deux composantes irr\'eductibles. Pour $\alpha \leq 3/2$ il
 est irréductible.
 \end{prop}

 \begin{pv} Soit $\Lambda$ un syst\`eme coh\'erent $\alpha_{-}-$stable.
 D'apr\`es le th\'eor\`eme \ref{critlissite}, il existe un voisinage du
 sch\'ema point\'e $(\mathcal{S}_{-}, \left[ \Lambda \right])$ qui est isomorphe \`a un sous-sch\'ema
ferm\'e
 d'un voisinage lisse $(\mathcal{N}, \left[ \Lambda \right])$, o\`u $\mathcal{N}$ est
 une vari\'et\'e d'espace tangent $\Ext^1 (\Lambda, \Lambda)$ en ce point,
 le faisceau d'id\'eaux \'etant localement d\'efini par
 $\Ext^2 (\Lambda, \Lambda)$ \'equations.
 Si maintenant $\Lambda$ d\'efinit
 un point de $\Sigma_{-}$ appartenant \`a la strate $B_{i}$, alors en
 posant $a = \fc{n}{2} + \alpha - i$ on a
$$ \dim \Ext^1 (\Lambda, \Lambda) = n + 2 a + 3, \hspace{0.4cm} \dim
\Ext^2 (\Lambda, \Lambda) = a - b - 2 $$
 Notons $I$ le faisceau d'id\'eaux associ\'e au
 plongement compos\'e $(\Sigma_{-}, \Lambda) \subset (\mathcal{S}_{-},
 \Lambda) \subset (\mathcal{N}, \Lambda)$, et $J$ celui du premier plongement.
 On a un morphisme naturel surjectif de faisceaux conormaux
 $$ I / I^2 |_{B_{i}} \ri J / J^2 |_{B_{i}} \ri 0 $$
 qui induit un plongement naturel $\mathbb{P}(J / J^2 |_{B_{i}})
 \subseteq \mathbb{P} (I / I^2 |_{B_{i}})$. Le terme de gauche est
 l'image r\'eciproque de $B_{i}$ dans l'\'eclat\'e, le terme de droite
 s'identifie \`a l'image inverse de $B_{i}$ dans l'\'eclat\'e de
 $\mathcal{N}$ le long de $\Sigma_{-}$. Le lieu singulier du
 faisceau $I / I^2 |_{B_{i}}$ est de support contenu dans $A \cap B$
 (si le voisinage consid\'er\'e a une intersection non vide
 avec ce ferm\'e), et donc ce module coh\'erent est sans torsion car
 $B_{i}$ est lisse. Sur l'ouvert compl\'ementaire de $B \setminus A$ le rang du fibr\'e $I / I^{2} |_{B_i}$ est donc $n + 2 a + 3 - (n-2) = 2 a + 5$. L'image r\'eciproque de la strate
 $B_{i}$ priv\'ee de $B \setminus A$ dans l'\'eclatement de $\mathcal{N}$ le long de $\Sigma_{-}$ est donc de dimension $\leq 2 a + 4 + \dim B_i = 2 (a + i) + 3 + n - 2 \alpha = 2 n + 3$. Par cons\'equent l'image r\'eciproque d'une telle strate dans l'\'eclat\'e de $\mathcal{S}_{-}$ le long de $\Sigma_{-}$ ne peut \^etre une composante irr\'eductible.
  \end{pv}

  \subsection{Sauts au passage d'une valeur critique}

   \label{calculsauts}On conserve les notations du paragraphe pr\'ec\'edent. On pose de plus
$\mathcal{D}_{+/-} = \mathcal{D}_{\alpha_{+/-}}$ les fibr\'es
d\'eterminants sur les vari\'et\'es $\mathcal{S}_{+/-}$ introduits
au paragraphe \ref{fibredet}. On cherche \`a \'evaluer la
diff\'erence
   $$ \Delta = c_{1} (\mathcal{D}_{+})^{2 n + 5} \cap [\mathcal{S}_{+}] - c_{1}
(\mathcal{D}_{-})^{2 n + 5} \cap [\mathcal{S}_{-}]
\hspace{0.5cm}(+)$$
   Pour cela on utilise le diagramme du th\'eor\`eme \ref{eclatements}. La
quantit\'e $(+)$ est encore \'egale \`a $ (c_1 (\phi_{+}^*
\mathcal{D}_{+})^{2 n + 5} - c_1 (\phi_{-}^* \mathcal{D}_{-})^{2 n +
5}) \cap [\widetilde{\mathcal{S}}]$, o\`u $\varphi_{+/-} :
\widetilde{\mathcal{S}} \ri \mathcal{S}_{+/-}$ sont les
\'eclatements du th\'eor\`eme pr\'ec\'edent.

   \begin{lem}
   \label{sautfibdet}On a l'isomorphisme de fibr\'es inversibles
   $$ \varphi_{-}^* \mathcal{D}_{-} \simeq \O(E) \otimes \varphi_{+}^*
\mathcal{D}_{+} $$
   o\`u $E \subset \widetilde{S}$ est le diviseur exceptionnel de l'\'eclatement.
   \end{lem}

  \begin{pv} En utilisant les notations introduites au d\'ebut de la preuve du
th\'eor\`eme \ref{eclatements}, et d'apr\`es la propri\'et\'e des
fibr\'es d\'eterminants, on a les \'egalit\'es
   $$   \phi_{+/-}^* \mathcal{D}_{+/-} \simeq det \hspace{0.1cm}( p_{!}
\hspace{0.1cm}\phi_{+/-}^* (\mathbf{\Theta}_{+/-}) \cdot h  )
\otimes det \hspace{0.1cm} \mathbf{\Gamma}_{+/-}^* $$
   Or par hypoth\`ese nous avons la suite exacte de syst\`emes alg\'ebriques sur
$\widetilde{V_{+/-}}$
   $$ 0 \ri (\mathbf{\Gamma}_{+}, \mathbf{\Theta}_{+}) \ri (\mathbf{\Gamma}_{-},
\mathbf{\Theta}_{-}) \ri (0, \mathbf{\Theta}) \ri 0 $$
   o\`u $\mathbf{\Gamma}_{-} \simeq \mathbf{\Gamma}_{+}$ et $\mathbf{\Theta}$ est
param\`etre les faisceaux de la forme $\O_{l}(\fc{n}{2} + \alpha)$.
Comme $\mathbf{\Theta}$ est plat sur le diviseur exceptionnel, on en
d\'eduit que la classe $p_{!} \hspace{0.1cm} \mathbf{\Theta} \cdot
h$ est celle d'un faisceau inversible sur le diviseur exceptionnel.
   Le lemme suivant permet alors de conclure. \end{pv}

    \begin{lem}
    Soit $X$ une vari\'et\'e projective int\`egre et $E$ un diviseur de
    Cartier effectif. Soit $L$ un faisceau coh\'erent de support $E$, et
inversible sur $E$. Alors $det \hspace{0.1cm} L \simeq \O(E)$.
    \end{lem}

    \begin{pv} Tout vient du fait que $L$ est un faisceau de
     dimension homologique 1: en effet toute surjection $A \ri L \ri
     0$ o\`u $A$ est un fibr\'e vectoriel sur $X$ a pour noyau
     un fibr\'e $B$, car en tout point le germe de $L$ est
     isomorphe au germe de $\O_{E}$. L'existence de $A$ est donn\'ee
     par le th\'eor\`eme de Serre. Par d\'efinition le d\'eterminant
     de $L$ est donn\'e par
     det$(A) \hspace{0.1cm} \otimes \hspace{0.1cm}$det$(B)^{-1}$, qui est
     isomorphe \`a $\O(E)$. \end{pv} \\

    \fbox{Notations:} on pose
    \begin{itemize}
    \item $d_{+/-} = c_1 (\mathcal{D}_{+/-})$;
    \item $L_{-}$ et $L_{+}$ les fibr\'es inversibles relativement
    amples sur les fibr\'es en espaces projectifs
    $\Sigma_{-} \ri \Sigma$ et $\Sigma_{+} \ri \Sigma$;
    on note leurs classes de Chern
$l_{+/-}$;
    \item $\gamma_{i}$, pour $i=1,2$, la $i-$\`eme classe de Chern du fibr\'e
tautologique $\mathbf{\Gamma}$;
    \item $\tau_1$ et $\tau_2$ les deux fibr\'es inversibles provenant
    des fibr\'es $\O_{\P_{2}^*}(1)$ venant de chacun des deux facteurs
    du produit $\P_{2}^* \times \P_{2}^*$; on pose $v$ la classe de Chern du fibr\'e inversible $\tau_{1}^{\otimes
\fc{n}{2} - \alpha} \otimes \tau_{2}^{\otimes \fc{n}{2} + \alpha}
\otimes det \hspace{0.1cm} \mathbf{\Gamma}^*$; on pose enfin $u = (2
\alpha - 1) \tau_1 + (2 \alpha - 2) \tau_2$;
   \item $\mathfrak{e}$ la classe du diviseur exceptionnel $E$.
 \end{itemize}

  On a alors les expressions suivantes des restrictions des fibr\'es
d\'eterminants $\mathcal{D}_{+/-}$ aux ferm\'es $\Sigma_{+/-}$.
    \begin{lem}
    \label{restrfibdet}On a:
    $$ \mathcal{D}_{-} |_{\Sigma_{-}} \simeq L_{-}^{-1} \otimes
    \tau_{1}^{\otimes
\fc{n}{2} - \alpha} \otimes \tau_{2}^{\otimes \fc{n}{2} + \alpha}
\otimes det \hspace{0.1cm} \mathbf{\Gamma}^*, \hspace{0.2cm}
\mathcal{D}_{+} |_{\Sigma_{+}} \simeq L_{+} \otimes \tau_{1}^{\otimes
\fc{n}{2} - \alpha} \otimes \tau_{2}^{\otimes \fc{n}{2} + \alpha}
\otimes det \hspace{0.1cm} \mathbf{\Gamma}^* $$
    \end{lem}

     \begin{pv}
 D'apr\`es la preuve de la proposition \ref{sigma-} et les notations qui la pr\'ec\`edent on a une suite exacte de familles de syst\`emes coh\'erents param\'etr\'ees par
 $\Sigma_{-}$ donn\'ee par
 $$ 0 \ri L_{-} \otimes \mathbf{\Lambda}^{'} \ri \mathbf{\Lambda}
 \ri \mathbf{\Lambda}^{''} \ri 0 $$
 Par la propri\'et\'e universelle du fibr\'e $\mathcal{D}_{-}$ on obtient
 $$ \mathcal{D}_{-} |_{\Sigma_{-}} \simeq L^{\otimes - 2}
 \otimes \det(\mathbf{\Gamma})^* \otimes \det ( L_{-} \otimes Rp_{!} (\O_{D_1}(\fc{n}{2} - \alpha) \cdot h)) \otimes \det (Rp_{!} (\O_{D_2}(\fc{n}{2} + \alpha) \cdot h))$$
 Rappelons que $h$ est la classe du faisceau $\O_{l}$  o\`u $l$ est une droite, on peut donc \'ecrire $\O_{D_i} \cdot h = \O_{D_i} - \O_{D_i}(-1)$ dans le groupe de Grothendieck de $\Sigma_{-} \times \P_2$; on obtient donc
 $$ \det ( L_{-} \otimes Rp_{!} (\O_{D_1}(\fc{n}{2} - \alpha) \cdot h)) \otimes \det (Rp_{!} (\O_{D_2}(\fc{n}{2} + \alpha) \cdot h)) = L_{-}^{\otimes \fc{n}{2} - \alpha +1} \otimes
 (L_{-} \otimes \tau_{1}^{-1})^{\otimes \alpha - \fc{n}{2}} \otimes \tau_{2}^{\otimes \fc{n}{2} + \alpha} $$
 ce qui donne la premi\`ere partie du lemme. La seconde se d\'emontre \`a l'identique, la diff\'erence venant de ce que $\mathbf{\Gamma}$ n'est pas tensoris\'ee par $L_{+}$. \end{pv} \\

   On a alors
   $$ \Delta  = - \fc{1}{\mathfrak{e}} \cdot ( \varphi_{+}^* (d_{+})^{2 n + 5} |_{E}
- \varphi_{-}^* (d_{-})^{2 n + 5} |_{E} ) \cap [E] $$
 o\`u l'expression $\fc{1}{\mathfrak{e}} \cdot ( \varphi_{+}^* (d_{+})^{2 n + 5} |_{E}
- \varphi_{-}^* (d_{-})^{2 n + 5} |_{E} )$ d\'esigne la division de
la classe entre parenth\`eses par la classe $e$: on sait en effet
d'apr\`es le lemme \ref{sautfibdet} que la classe $\varphi_{+}^*
(d_{+})^{2 n + 5} - \varphi_{-}^* (d_{-})^{2 n + 5}$ est un multiple
de la classe du diviseur exceptionnel. D'autre part le lemme
\ref{restrfibdet} donne:
$$ \varphi_{+}^* (d_{+}) |_{E} = l_{+} + v, \hspace{0.3cm} \varphi_{-}^*
(d_{-}) |_{E} = - l_{-} + v $$
 et on peut \'ecrire
 $$ \begin{array}{llcc}
   (l_{+} + v)^{2 n + 5} - (- l_{-} + v)^{2 n + 5} & = \sum_{i=0}^{2
   n + 5} \binom{2 n + 5}{i} ( l_{+}^{i} - (-1)^i l_{-}^i ) \cdot
   v^{2 n + 5 - i} \\
   & = (l_{+} + l_{-}) \cdot \sum_{ \begin{array}{ll}
                                    0 \leq i \leq 2 n + 5 \\
                                    0 \leq j \leq i - 1 \\
                                    \end{array} }
                                    \binom{2 n + 5}{i} l_{+}^{j}
                                    (- l_{-})^{i - j - 1}
                                    \cdot v^{2 n + 5 - i} \\
   \end{array} $$
 Apr\`es r\'eindexation de la somme pr\'ec\'edente, on peut donc
 \'ecrire
 $$ \Delta = \sum_{ \begin{array}{ll}
                                    0 \leq i \leq 2 n + 4 \\
                                    0 \leq j \leq i \\
                                    \end{array} }
                                    \binom{2 n + 5}{i+1} (-1)^{i+j}
                                    \hspace{0.2cm} l_{+}^{j} \cdot
                                    l_{-}^{i - j}
                                    \cdot v^{2 n + 4 - i} \cap [E]
                                    \hspace{4cm} (\mathbf{\Delta})$$
Lorsque $\alpha \in \{ \fc{1}{2},1,\fc{3}{2} \}$, on sait que $\Sigma_{-}$ est irr\'eductible, ainsi que $E$. Lorsque $\alpha \geq 2$ le diviseur $E$ est r\'eunion des deux
 composantes $A^{'}$ et $B^{'}$, qui sont les images r\'eciproques
 dans l'\'eclat\'e des deux composantes $A$ et $B$. On peut donc \'ecrire
 $\Delta = \Delta_{A} + \Delta_{B}$, o\`u les deux termes de droite sont
 respectivement les contributions de $A^{'}$ et $B^{'}$.
 Dans les calculs qui suivent on ne tiendra compte que de $\Delta_{A}$ si $\alpha < 2$. \\

\fbox{{\bf{Calculons $\Delta_{A}.$}}}
 Il existe un morphisme birationnel
 $\rho: A^{'} \ri A \times_{\Sigma} \Sigma_{+}$, et d'apr\`es le
 diagramme commutatif
 $$ \xymatrix{
    & A^{'} \ar[dl]_{\varphi_{-}} \ar[dr]^{\varphi_{+}} & \\
    A \ar[dr]_{\pi_{-}} & & \Sigma_{+} \ar[dl]^{\pi_{+}} \\
    & \Sigma & \\
    } $$
    et la formule de projection, on peut \'ecrire
    $$ v^{2 n + 4 - i} \cdot l_{+}^{j} \cdot l_{-}^{i  - j} \cap [A^{'}] =
    v^{2 n + 4 - i} \cdot \pi_{+ \hspace{0.1cm}*} (l_{+}^{j}) \cdot
    \pi_{- \hspace{0.1cm}*} (l_{-}^{i - j}) \cap
    [\Sigma] $$
    L'image directe $\pi_{+ \hspace{0.1cm}*} (l_{+}^{j})$
    s'interpr\`ete comme une classe de Segr\'e d'un fibr\'e
    vectoriel, \`a savoir le dual celui d\'efini \`a la proposition
    \ref{sigma+} et qui est not\'e $\underline{\Ext}^{1}_{p}(\mathbf{\Lambda}^{'},
    \mathbf{\Lambda}^{''})$, et dont le projectif
    au sens de Grothendieck (c'est-\`a-dire param\'etrant les quotients
    de rang 1) donne $\Sigma_{+}$.
    Par l'\'ecriture $\pi_{- \hspace{0.1cm}*} (l_{-}^{i - j}) \cap
    [\Sigma]$, nous entendons en fait l'image directe propre $\pi_{-
    \hspace{0.1cm}*} (l_{-}^{i - j} \cap [A])$ au sens des
    cycles (le cycle entre parenth\`eses \'etant de dimension $\dim
    A - (i - j) = n + j - 2 \alpha - i + 2$). Le probl\`eme ici est
    que $\mathcal{J}^{2 \alpha - 2} (2 \alpha - 1, 2 \alpha -
    2)$ est un faisceau coh\'erent singulier.
    Calculer sa classe de Segr\'e, qui est par d\'efinition la classe de
    Segr\'e du c\^one qu'il d\'efinit, n\'ecessite quelques
    pr\'ecautions (voir [F] ex. 4.1.6-7). Aussi nous allons faire en
    sorte de nous ramener \`a des classes de Segr\'e de fibr\'es
    vectoriels. \\

 Rappelons que $A$ est isomorphe, quand $\alpha \geq 1$, \`a
l'\'eclat\'e de $\Sigma$ le
     long du ferm\'e d'id\'eal $\mathcal{J}^{2 \alpha - 2}$, c'est-\`a-dire
     \`a Proj\hspace{0.1cm}($\oplus \mathcal{J}^{(2 \alpha - 2) *}$), par la
     proposition \ref{compsigma-}. Soit $\sigma$ la classe du
     fibr\'e canonique $\O(1)$ sur cet \'eclat\'e, qui est encore le
     faisceau d'id\'eaux inversible du diviseur exceptionnel. Alors
     en suivant les conventions ci-dessus on peut \'ecrire
     $l_{-} |_{A} = \sigma + u $. On en d\'eduit alors
     $$ \pi_{- \hspace{0.1cm} *} (l_{-}^{i-j} |_{A} \cap [A]) = \sum_{k=0}^{i-j}
     \binom{i-j}{k} \hspace{0.1cm} u^{i-j-k} \cdot \pi_{- \hspace{0.1cm} *} (\sigma^{k}
     \cap [A]) \hspace{2cm} (*)$$
 \begin{lem}
 \label{imdirpropre}L'image directe propre $\pi_{- \hspace{0.1cm} *} (\sigma^{k}
     \cap [A])$ est \'egale \`a $[\Sigma]$ si $k=0$ et \`a
     $- (2 - 2 \alpha)^k \hspace{0.1cm}
     \fc{k (k-1)}{2} \hspace{0.1cm} \cdot \tau^{k-2}$ sinon, o\`u
     $\tau$ est la classe du faisceau
     $\O_{\P_{2}^{*}}(1)$ port\'e par le sous-sch\'ema ferm\'e de $\Sigma$
     d'id\'eal $\mathcal{J}$. Autrement dit $\tau$ est la
     restriction de $\tau_i$ pour $i=1,2$ \`a ce sous-sch\'ema.
 \end{lem}

 \begin{pv} Posons $\mathcal{J}^{'} = \mathcal{J}^{2 \alpha - 2}$. Soit
 $E_{A}$ le diviseur exceptionnel de $A$. Il est l'image r\'eciproque
 au sens sch\'ematique du sous-sch\'ema ferm\'e d'id\'eal $\mathcal{J}^{'}$
 par la projection $A \ri \Sigma$. On peut encore \'ecrire que $E_{A}$
 est isomorphe \`a
 Proj\hspace{0.1cm}($\oplus_{k} \mathcal{J}^{' k} / \mathcal{J}^{' k +
 1}$). Comme $\sigma$ est la classe du faisceau
 $\O (- E_{A})$, on peut \'ecrire
$$ \pi_{- \hspace{0.1cm} *} (\sigma^k \cap [A]) = - \pi_{-
\hspace{0.1cm} *} (\sigma^{k-1} \cap [E_{A}]) $$
 L'image inverse dans $A$ de la "diagonale r\'eduite" d'id\'eal
 $\mathcal{J}$ est isomorphe \`a Proj\hspace{0.1cm}($\oplus_{k}
 \mathcal{J}^{' k} / \mathcal{J}^{(2 \alpha - 2) k + 1}$),
 lui-m\^eme isomorphe \`a $\mathbb{P}(\mathcal{J} / \mathcal{J}^2 )$
 plong\'e dans $\mathbb{P}(\S^{2 \alpha - 2} (\mathcal{J} /
 \mathcal{J}^2 ))$ via le plongement de Veronese. On a donc
 un plongement $E_{A}^{'} = \mathbb{P}(\mathcal{J} / \mathcal{J}^2 ) \subset
 E_{A}$, qui est un hom\'eomorphisme topologique. L'id\'eal
 du morphisme surjectif d'alg\`ebres
$$ \oplus_k \mathcal{J}^{' k } / \mathcal{J}^{' k + 1} \ri
\oplus_{k} \mathcal{J}^{' k} / \mathcal{J}^{(2 \alpha - 2) k + 1} $$
est engendr\'e par la partie de degr\'e 0 qui est $\mathcal{J} /
\mathcal{J}^{'}$; c'est donc un id\'eal nilpotent d'indice $2 \alpha
- 2$ et on a l'\'egalit\'e $[E_{A}] \equiv (2 \alpha - 2)
\hspace{0.1cm} [E_{A}^{'}]$. La restriction de $\sigma$ \`a
$\mathbb{P}(\mathcal{J} / \mathcal{J}^2 )$ est \'egale \`a
$\O_{\mathbb{P}(\mathcal{J} / \mathcal{J}^{2})}(2 \alpha - 2)$,
c'est-\`a-dire en termes de classe \`a $(2 \alpha - 2) \sigma^{'}$,
avec $\sigma^{'}$ la classe du fibr\'e tautologique
$\O_{\mathbb{P}(\mathcal{J} / \mathcal{J}^{2})}(1)$. Par suite on a
$$ \begin{array}{llcc}
   \pi_{- \hspace{0.1cm} *} (\sigma^{k} \cap [A]) & =
   - (2 \alpha - 2)^k \hspace{0.1cm}
   \pi_{- \hspace{0.1cm} *} (\sigma^{' k-1} \cap [E_{A}^{'}]) \\
    & = - (2 \alpha - 2)^k \hspace{0.1cm} s_{k-2}(\mathcal{J} /
    \mathcal{J}^2 ) \\
  \end{array} $$
 Or le fibr\'e $\mathcal{J} / \mathcal{J}^2 $ a pour support
 la grassmanienne relative $G$ et est isomorphe \`a l'image
 r\'eciproque du fibr\'e $Q^* (-1)$ par la projection $G \ri
 \P_{2}^*$. On \'ecrit ensuite
$$ \begin{array}{llcc}
  s_{k-2} (Q^* (-1)) & = [ c (Q(1))^{-1} ]_{k-2} \\
   & = [ (1 + h)^{-3} ]_{k-2} \\
   & = (-1)^k \hspace{0.1cm} \fc{k (k-1)}{2} \hspace{0.1cm} \tau^{k-2} \\
   \end{array} $$
 ce qui montre le lemme. \end{pv} \\

 Le lemme suivant est une cons\'equence du lemme \ref{imdirpropre}.

\begin{lem}
 \label{premiereimagedirecte}Utilisant les notations pr\'ec\'edentes, lorsque $\alpha \geq 1$ l'image directe $\pi_{- \hspace{0.1cm} *} (l_{-}^{i-j} |_{A} \cap [A])$
 est alors \'egale \`a
 $$ \epsilon_{i-j} = u^{i-j} - a_{i-j} \hspace{0.1cm}\tau^{i-j-2}$$
 o\`u $a_{i-j} = 2 (i-j) (i-j-1) (2 \alpha - 1)^{i-j-2} (\alpha - 1)^2$.
 Notons que pour $i-j \geq 5$ cette classe est nulle.
 Lorsque $\alpha = \fc{1}{2}$ on a \hspace{0.1cm} $l_{-} \simeq \tau_{2}^{\vee}$ et
 $\pi_{- \hspace{0.1cm} *} (l_{-}^{i-j}) = (-1)^{i + j}
 \hspace{0.1cm} \tau_{2}^{i-j}$.
\end{lem}

 \begin{pv}
 La seconde assertion est imm\'ediate car pour $\alpha = 1/2$ on a $\Sigma_{-}
 \simeq \Sigma$ et $l_{-} \simeq \tau_{2}^{\vee}$. Pour $\alpha \geq 1$, d'apr\`es ce qui pr\'ec\`ede la classe de l'\'enonc\'e s'exprime sous la forme
 $u^{i-j} - a_{i-j} \hspace{0.1cm}\tau^{i-j-2}$ avec $a_{i-j} = \fc{1}{2} \sum_{k=2}^{i-j}
 k (k-1) \binom{i-j}{k} (4 \alpha - 3)^{i-j-k} (2 - 2 \alpha)^k$. En vertu de l'identit\'e polyn\^omiale valable dans $\mathbb{C}[X,Y]$ et pour tout entier $m > 0$, et qui consiste en une diff\'erentiation de l'identit\'e du bin\^ome:
 $$ m (m-1) (Y + X)^{m-2} = \sum_{k=2}^{m} \binom{m}{k}
 \hspace{0.1cm} k (k-1) \hspace{0.1cm} Y^{m-k} X^{k-2} $$
 on en d\'eduit l'\'enonc\'e.  \end{pv} \\

 On peut donc \'ecrire $\Delta_A = \Delta_{A}^* +
 \Delta_{A}^{**}$, o\`u $\Delta_{A}^*$ est la contribution apport\'ee
 par les termes $u^{i-j}$ dans la classe $(\pi_{-})_* (l_{-}
 |_{A}^{i-j})$, et $\Delta_{A}^{**}$ est celle apport\'ee par
 $a_{i-j} \hspace{0.1cm}\tau^{i-j-2}$ (classe de support inclus dans la
 diagonale).

 Explicitons maintenant l'image directe $\pi_{+ \hspace{0.1cm} *} (l_{+}^j
 )$. Posons $\mathcal{E}$ le fibr\'e vectoriel sur $\Sigma$ d\'efini par
$\underline{\Ext}^{1}_{p}(\mathbf{\Lambda}^{'},
    \mathbf{\Lambda}^{''})$. On a par d\'efinition
    $\pi_{+ \hspace{0.1cm} *} (l_{+}^j ) = s_{j - n - 2 \alpha - 2}
    (\mathcal{E}^{\vee})$ car on a $\Sigma_{+} = \mathbb{P}
    (\mathcal{E}^{\vee})$ et le fibr\'e projectif est de rang $n + 2
    \alpha + 2$. On sait que $\mathcal{E}$ s'ins\`ere dans une suite
    exacte
    $$ 0 \ri \mathbf{\Gamma}^* \boxtimes \S^{\fc{n}{2} + \alpha} Q
    \ri \mathcal{E} \ri \underline{\Ext}^{1}_{p} (\O_{D_{1}},
    \O_{D_{2}}(2 \alpha)) \ri 0 $$
    Nous avons vu au cours de la preuve de la proposition
    \ref{sigma-} que le terme de droite est isomorphe au faisceau
    $\mathcal{J}^{2 \alpha + 1} ( 2 \alpha + 2, 2 \alpha + 1)$, qui
    admet \'egalement la r\'esolution suivante
    $$ 0 \ri \O \boxtimes \S^{2 \alpha} Q \ri \tau_1 \boxtimes \S^{2
    \alpha + 1} Q \ri \mathcal{J}^{2 \alpha + 1} ( 2 \alpha + 2, 2 \alpha + 1)
    \ri 0 $$
    on a donc l'identit\'e suivante des classes de Chern totales
    $$ c( \mathcal{E}) \hspace{0.1cm}
    c (\mathbf{\Gamma}^* \boxtimes \S^{\fc{n}{2}
    + \alpha} Q)^{-1} = c(\tau_1 \boxtimes \S^{2
    \alpha + 1} Q) \hspace{0.1cm} c(\O \boxtimes \S^{2 \alpha}
    Q)^{-1} $$
    qui entra\^ine l'identit\'e suivante entre classes de Segr\'e
    $$ \begin{array}{llccc}
    s (\mathcal{E}^{\vee}) & = & c (\mathcal{E})^{-1}
     & = & s(\mathbf{\Gamma} \boxtimes \S^{\fc{n}{2}
    + \alpha} Q^*) \hspace{0.1cm} c (\O \boxtimes \S^{2 \alpha} Q )
    \hspace{0.1cm} c (\tau_1 \boxtimes \S^{2
    \alpha + 1} Q)^{-1} \\
    \end{array} $$
    On pose $\delta =  \sum_{*} \delta_{*}$ la classe totale \'egale \`a
    $c (\O \boxtimes \S^{2 \alpha} Q )
    \hspace{0.1cm} c (\tau_1 \boxtimes \S^{2
    \alpha + 1} Q)^{-1}$.
 On a l'\'egalit\'e
 $$ \delta = \fc{(1 + \tau_1 - \tau_2)^{\binom{2 \alpha + 2}{2}}}{(1 +
 \tau_1)^{\binom{2 \alpha + 3}{2}} \cdot (1 - \tau_2)^{\binom{2
 \alpha + 1}{2}}}$$
 On a donc
 $$ \Delta_{A}^* = \sum_{n + 2 \alpha + 3 \leq j \leq i \leq 2 n +
 5} \hspace{0.1cm} (-1)^{i-j} \hspace{0.1cm} \binom{2 n + 5}{i}
 \hspace{0.1cm} v^{2 n + 5 - i} \cdot u^{i-j} \cdot s_{j-n-2 \alpha
 - 3}(\mathcal{E}^{\vee}) \cap [G \times \P_{2}^*] $$
 o\`u les indices $i$ et $j$ d\'esignent les indices $i+1$ et $j+1$
 dans l'\'ecriture $(\mathbf{\Delta})$ qui pr\'ec\`ede le paragraphe
 pr\'esent consacr\'e au calcul de $\Delta_A$. En d'autres termes on
 peut \'ecrire que
 $$ \Delta_{A}^* = \left[ (1 + v)^{2 n + 5} \cdot \fc{1}{1 + u} \cdot
 s(\mathcal{E}^{\vee}) \right]_{n - 2 \alpha + 2} \cap [G \times
 \P_{2}^*] $$
 o\`u la notation $[\hspace{0.5cm}]_{n - 2 \alpha + 2} $ d\'esigne
 la partie de degr\'e $n - 2 \alpha + 2$ de l'expression entre
 crochets. Par ailleurs on a
 $$ \begin{array}{ll}
 (1 + v)^{2 n + 5} \cdot \fc{1}{1 + u} \cdot
 s(\mathcal{E}^{\vee}) = (1 + v)^{2 n + 5} \cdot \fc{1}{1 + u} \cdot
 [\sum_{l} s_{l} (\mathbf{\Gamma} \boxtimes \S^{\fc{n}{2} + \alpha}
 Q^*)] \cdot \delta \\
  = \sum_{0 \leq m \leq 2 n + 5} \hspace{0.1cm}
  \binom{2 n + 5}{m} \hspace{0.1cm} (-1)^m \hspace{0.1cm} \fc{\left[ 1 +
  (\fc{n}{2} - \alpha) \tau_1 + (\fc{n}{2} + \alpha) \tau_2 \right]^m
  \cdot \delta}{1 + u} \cdot \gamma_{1}^{2 n + 5 - m} \cdot
  \left[ \sum_l s_l (\mathbf{\Gamma} \boxtimes \S^{\fc{n}{2} +
  \alpha} Q^*) \right] \\
  \end{array} $$
 Pour tout r\'eel $x$ on pose $(x)_{+} = \mathrm{max} (x,0)$. On a donc
 $$ \begin{array}{llcc}
 \Delta_{A}^* & = \sum_{ \begin{array}{ll}
                          n + 2 \alpha + 3 \leq m \leq 2 n + 5 \\
                          (m - n - 2 \alpha - 7)_+ \leq l \leq m - n
                          - 2 \alpha - 3 \\
                          \end{array} }
                          (-1)^m \binom{2 n + 5}{m}
                          \fc{\left[ 1 +
  (\fc{n}{2} - \alpha) \tau_1 + (\fc{n}{2} + \alpha) \tau_2 \right]^m
  \cdot \delta}{1 + u} \\
  & \hspace{8cm} \cdot \gamma_{1}^{2 n + 5 - m} \cdot
  s_{l}(\mathbf{\Gamma} \boxtimes \S^{\fc{n}{2} +
  \alpha} Q^*) \cap [G \times \P_{2}^*] \\
  \end{array} $$

 Le lemme suivant va nous permettre d'achever ce calcul.

    \begin{lem}
 \label{secondeimagedirecte} Soit $\rho: G \times \P_{2}^* \ri \P_{2}^* \times \P_{2}^*$ la
 projection canonique. On pose les entiers $\zeta = \binom{\fc{n}{2}
 + \alpha + 1}{2}$ et $\zeta^{'} = \binom{\fc{n}{2} - \alpha +
 1}{2}$. Alors l'image directe
 $\rho_* (c_{1}(\mathbf{\Gamma})^k \cdot s_l (\mathbf{\Gamma}
 \boxtimes \S^{\fc{n}{2} + \alpha} Q^*))$ est \'eventuellement non nulle
 si et seulement si en posant $l + k = n - 2 \alpha - 2 + \lambda$, on a $0 \leq \lambda
\leq 4$
 , et dans ce cas est \'egale \`a
 $$ \fbox{$\displaystyle \fc{1}{\zeta^{'} + 1} \hspace{0.1cm}
 \sum_{ \begin{array}{llcc}
  0 \leq l_0 \leq l_1 \leq l_2 \leq l_3 \leq l \\
  \lambda -2 \leq l_0 + (l_3 - l_2) \leq \mathrm{min}(\lambda, 2) \\
  \end{array} } f_{\alpha,n}(\underline{{\bf{l}}},k) \hspace{0.1cm}
 \tau_{1}^{\lambda - l_0 - (l_3 - l_2)} \cdot \tau_{2}^{l_0 + (l_3 -
 l_2)} $} $$
 o\`u la notation $\underline{{\bf{l}}}$ d\'esigne le vecteur $(l_0, l_1, l_2, l_3)$
 et o\`u la fonction num\'erique $f_{\alpha,n}(\underline{{\bf{l}}},k)$ d\'ependant des
 param\`etres $l,k,\alpha,n$ est \'egale \`a
 $$ \fbox{$\displaystyle \begin{array}{llll}
   f_{\alpha,n}(\underline{{\bf{l}}},k) = (-1)^{l_0 + l_1 + l_3} \hspace{0.1cm} \binom{\fc{n}{2} + \alpha +
 l_2}{l_2 - l_1} \hspace{0.1cm} \binom{\fc{n}{2} + \alpha +
 l_1}{l_1 - l_0} \hspace{0.1cm} \binom{\fc{n}{2} + \alpha
 + l - l_2}{l - l_3} \hspace{0.1cm} \binom{\zeta}{l_0}
 \hspace{0.1cm} \binom{\zeta}{l_3 - l_2} \\
  \sum_{ \begin{array}{ll}
                    0 \leq a \leq l_1 - l_0 + k \\
                    0 \leq b \leq l + l_2 - l_1 - l_3 \\
                    \end{array} } \hspace{0.1cm} \binom{l_1 + k -
                    l_0}{a} \hspace{0.1cm}
                    \binom{\zeta^{'} + 1}{\fc{n}{2} - \alpha + \lambda + l_2
                    - l_0 - l_3 - a - b}
                    \hspace{0.1cm} \binom{\zeta^{'} + 1}{a + b + \alpha + 1 -
                    \fc{n}{2}} \\
 \hspace{5.5cm} (n - 2 \alpha - 1 + l_2 - l_0 - l_3 - 2 a - 2 b + \lambda) \hspace{0.1cm} \\
 \end{array} $} $$
 On d\'eduit de cette formule que si $\rho$ d\'esigne cette fois la
 projection $G \ri \P_{2}^*$, l'image directe $\rho_* (c_{1}(\mathbf{\Gamma})^k \cdot s_l (\mathbf{\Gamma}
 \otimes \S^{\fc{n}{2} + \alpha} Q^*))$ (le produit tensoriel
 remplace le produit tensoriel externe dans le second facteur) est
 \'egale \`a
 $$ \fc{1}{\zeta^{'} + 1} \hspace{0.1cm} g_{\alpha,n}(l,k)
 \hspace{0.1cm} \tau^{\lambda} $$
 o\`u $g_{\alpha,n}(l,k) = \sum_{\underline{{\bf{l}}}} f_{\alpha,n}(\underline{{\bf{l}}}, k)$, la somme portant sur les
 quadruplets $\underline{{\bf{l}}}$ v\'erifiant les
 in\'egalit\'es ci-dessus. Naturellement si $\lambda > 2$ cette classe est nulle.
  \end{lem}

Pour montrer le lemme \ref{secondeimagedirecte} nous montrerons deux
r\'esultats p\'eliminaires. Le premier est une cons\'equence assez
facile de l'exemple 3.1.1 de [F]. Le second est un calcul d'image
directe qui se m\`ene de mani\`ere classique en faisant appel \`a
une vari\'et\'e de drapeaux relative.

\begin{lem}
\label{segretensor}Soit $E$ un fibr\'e vectoriel de rang 2 et $F$ un
fibr\'e vectoriel de rang $f + 1 > 0$ sur une vari\'et\'e projective
$X$. Soit $k \geq 0$ un entier. Alors $s_{k} (E \otimes F)$ est
\'egal \`a
$$ \begin{array}{ll}
 \sum_{0 \leq k_0 \leq k_1 \leq k_2 \leq k_3 \leq k}
\hspace{0.1cm} (-1)^{k_1 + k_2} \binom{f + k_2}{k_2 - k_1} \binom{f
+ k -
k_2}{k - k_3} \binom{f + k_1}{k_1 - k_0} \\
 \hspace{4cm} s_{k_0}(F) \cdot
s_{k_3 - k_2}(F) \cdot c_1 (E)^{k_1 - k_0} \cdot s_{k + k_2 - k_1
- k_3}(E) \\
 \end{array} $$
\end{lem}

\begin{pv}
On se place sur le fibr\'e projectif $\P (E) \stackrel{\pi}{\ri}
X$ qui param\`etre les droites vectorielles des fibres de $E$.
Soit $L$ le fibr\'e inversible quotient de $\pi^*(E)$, on note $l$
sa classe de Chern. On a la suite exacte de fibr\'es
$$ 0 \ri \pi^* (\Lambda^2 E) \otimes L^* \ri \pi^* (E) \ri
L \ri 0 $$ On a $\pi_* (s_{k} (\pi^*(E) \otimes \pi^*(F)) \cdot l) =
s_{k} (E \otimes F)$. On a alors $s_{k} (\pi^* (E \otimes F)) =
\sum_{0 \leq k_2 \leq k} s_{k_2} (\pi^* (\Lambda^2 E \otimes F)
\otimes L^*) \cdot s_{k - k_2} (\pi^* (F) \otimes L)$, ceci \'etant
\'egal d'apr\`es [F] ex.3.1.1 \`a
$$ \sum_{0 \leq k_1 \leq k_2 \leq k_3 \leq k} (-1)^{k_1 + k_2}
\hspace{0.1cm} \binom{f + k_2}{k_2 - k_1} \hspace{0.1cm} \binom{f +
k - k_2}{k - k_3} \hspace{0.1cm} s_{k_1} (\pi^* (\Lambda^2 E \otimes
F)) \cdot s_{k_3 - k_2} (\pi^* (F)) \cdot l^{k + k_2 - k_1 - k_3}
$$
 En \'ecrivant $s_{k_1} (
\Lambda^2 E \otimes F ) = \sum_{0 \leq k_0 \leq k_1} \binom{f +
k_1}{k_1 - k_0} \hspace{0.1cm} s_{k_0}(F) \cdot c_{1}(\Lambda^2
E)^{k_1 - k_0}$. On conclut gr\^ace \`a l'identit\'e $s_{k + k_2 -
k_1 - k_3} (E) = \pi_* (l^{k + k_2 - k_1 - k_3 + 1})$.
\end{pv} \\

 Soit \`a pr\'esent $\mathcal{E}$ un fibr\'e vectoriel de rang $e$
 sur une vari\'et\'e projective $X$. Soit $G(2, \mathcal{E})$ le
 fibr\'e en grassmaniennes au dessus de $X$ param\'etrant les
 plans vectoriels des fibres de $\mathcal{E}$. Soit $\mathbf{\Gamma}$ le sous-fibr\'e
 tautologique. On note
 $\mathcal{D}^{1,2}(\mathcal{E})$ la vari\'et\'e de drapeaux
 param\'etrant les couples $(\Gamma^{'}, \Gamma^{''} = \Gamma /
 \Gamma^{'})$ avec $\dim \Gamma^{'} = \dim \Gamma^{''} = 1$, et $\Gamma^{'}
 \subset \mathcal{E}_x$ pour $x \in X$
 et $\Gamma^{''} \subset \mathcal{E}_x / \Gamma^{'}$. On a un diagramme commutatif
$$ \xymatrix{
   \mathcal{D}^{1,2}(\mathcal{E}) \ar[d]_{\tau^{'}} \ar[r]^{\tau} & G(2,
   \mathcal{E}) \ar[d]^{\pi} \\
   \mathbf{P}(\mathcal{E}) \ar[r]_{\pi^{'}} & X \\
   } $$
o\`u les morphismes sont les projections naturelles, et
$\mathbf{P}(\mathcal{E})$ le fibr\'e projectif param\'etrant les
droites vectorielles des fibres de $\mathcal{E}$. Soient $\delta,
\epsilon$ deux entiers positifs.

\begin{prop}
\label{imdirclasseschernsegre}Avec les notations qui pr\'ec\`edent
$$\pi_* (c_{1}(\mathbf{\Gamma})^{\delta} \cdot
s_{\epsilon}(\mathbf{\Gamma})) = \sum_{ \begin{array}{ll}
                                     0 \leq d \leq \delta \\
                                     0 \leq e \leq \epsilon \\
                                     \end{array} }
                                     \binom{\delta}{d} \hspace{0.1cm}
                                     \left| \begin{array}{llcc}
                                            s_{\delta + \epsilon - d- e + 2
                                            - rg(\mathcal{E})}(\mathcal{E})
                                            & s_{d + e + 1 - rg(\mathcal{E})}(\mathcal{E}) \\
                                            s_{\delta + \epsilon -
                                            d- e +
                                            3 - rg(\mathcal{E})}(\mathcal{E})
                                            & s_{d + e + 2 - rg(\mathcal{E})}(\mathcal{E}) \\
                                            \end{array} \right| $$
\end{prop}

\begin{pv}
Soit $L^{'}$ le fibr\'e $\pi^{'}-$ample sur
$\mathbf{P}(\mathcal{E})$. On a une suite exacte
$$ 0 \ri L^{' \vee} \ri \pi^{' *} (\mathcal{E}) \ri
\mathcal{E}^{'} \ri 0 $$ ainsi que
$$ 0 \ri L^{'' \vee} \ri \tau^{' *} (\mathcal{E}^{'}) \ri
\mathcal{E}^{''} \ri 0$$ et enfin l'extension
$$ 0 \ri \tau^{' *} (L^{' \vee}) \ri \tau^* (\mathbf{\Gamma}) \ri
L^{'' \vee} \ri 0 $$ Notant $l^{'}, l^{''}$ les classes de Chern
de $L^{'}$ et $L^{''}$ on obtient ainsi les relations sur les
classes de Chern et de Segr\'e
$$ c_1 (\tau^* (\mathbf{\Gamma})) = - \tau^{' *} (l^{'}) - l^{''}$$
On a de plus
$$ s_{\epsilon} (\tau^* (\mathbf{\Gamma})) = (-1)^{\epsilon} \hspace{0.1cm}
\sum_{e=0}^{\epsilon} \tau^{' *} (l^{' e}) \cdot l^{'' \epsilon -
e}$$ La classe $c_1(\mathbf{\Gamma})^{\delta} \cdot s_{\epsilon}
(\mathbf{\Gamma})$ est de plus \'egale \`a $- \tau_* (c_1 (\tau^*
(\mathbf{\Gamma}))^{\delta} \cdot s_{\epsilon} (\tau^*
(\mathbf{\Gamma})) \cdot l^{''})$. La classe de l'\'enonc\'e est
donc \'egale \`a
$$ - \sum_{ \begin{array}{ll}
          0 \leq d \leq \delta \\
          0 \leq e \leq \epsilon \\
          \end{array} }
          (-1)^{\delta + \epsilon} \binom{\delta}{d} \hspace{0.1cm}
          \pi^{'}_* \tau^{'}_* (\tau^{' *}(l^{' d + e})
          \cdot l^{'' \delta + \epsilon - d - e + 1})$$
du fait de la commutativit\'e du diagramme ci-dessus. La formule
de projection donne alors
$$ \begin{array}{lllcc}
   \pi^{'}_* \tau^{'}_* (\tau^{' *}(l^{' d + e})
          \cdot l^{'' \delta + \epsilon - d- e + 1}) & =
          \pi^{'}_* (l^{' d + e } \cdot
          \tau^{'}_* (l^{'' \delta + \epsilon - d - e +
          1})) \\
          & = \pi^{'}_* (l^{' d + e} \cdot s_{\delta + \epsilon
          - d - e + 1
          - rg(\mathcal{E}^{''})} (\mathcal{E}^{' \vee}) \\
          & = s_{d + e - rg(\mathcal{E}^{'})} (\mathcal{E}^{\vee})
          \cdot s_{\delta + \epsilon - d - e + 1 -
          rg(\mathcal{E}^{''})}(\mathcal{E}^{\vee}) - \\
          & \hspace{0.5cm} s_{d + e + 1 - rg(\mathcal{E}^{'})} (\mathcal{E}^{\vee})
          \cdot s_{\delta + \epsilon - d - e -
          rg(\mathcal{E}^{''})}(\mathcal{E}^{\vee}) \\
   \end{array}$$
   On a par d\'efinition $rg(\mathcal{E}^{'}) = rg(\mathcal{E})-1$
   et $rg(\mathcal{E}^{''}) = rg(\mathcal{E})-2$. En tenant compte
   du fait que $s_{k}(\mathcal{E}^{\vee}) = (-1)^k \hspace{0.1cm}
   s_{k}(\mathcal{E})$ on peut conclure. \end{pv} \\

 \begin{pv} (du lemme \ref{secondeimagedirecte}) Dans ce cas le fibr\'e $\mathcal{E}$ est le fibr\'e
 $\S^{\fc{n}{2} - \alpha} Q$ sur $\P_{2}^*$. Posant $\zeta^{'} =
 \binom{\fc{n}{2} - \alpha + 1}{2}$ il vient
 $s_k (\S^{\fc{n}{2} - \alpha} Q) = \binom{\zeta^{'}}{k} \tau^k$, et
 donc les mineurs intervenant dans la formule de la proposition
 \ref{imdirclasseschernsegre} sont \'egaux \`a
 $$ \left| \begin{array}{llcc}
   \binom{\zeta^{'}}{\delta + \epsilon - d - e + 1 - \fc{n}{2} +
 \alpha} & \binom{\zeta^{'}}{d + e - \fc{n}{2} + \alpha} \\
 \binom{\zeta^{'}}{\delta + \epsilon - d - e + 2 - \fc{n}{2} +
 \alpha} & \binom{\zeta^{'}}{d + e + 1 - \fc{n}{2} + \alpha} \\
 \end{array} \right|
 \hspace{0.1cm} \tau^{\delta + \epsilon + 2 - n + 2 \alpha} $$
c'est-\`a-dire, compte tenu de l'identit\'e (calcul direct)
$$ \binom{\zeta^{'}}{a} \binom{\zeta^{'}}{b+1} -
\binom{\zeta^{'}}{a+ 1} \binom{\zeta^{'}}{b} = \fc{a-b}{\zeta^{'} +
1} \hspace{0.1cm} \binom{\zeta^{'} +1}{a+1} \binom{\zeta^{'} +
1}{b+1} $$ \'egaux \`a
$$ \fc{\delta + \epsilon - 2 d - 2 e + 1}{\zeta^{'} + 1}
\hspace{0.1cm} \binom{\zeta^{'} + 1}{\delta + \epsilon - d - e + 2 -
\fc{n}{2} + \alpha} \hspace{0.1cm} \binom{\zeta^{'} + 1}{d + e + 1 -
\fc{n}{2} + \alpha} \hspace{0.1cm} \tau^{\delta + \epsilon + 2 - n +
2 \alpha}$$

 D'apr\`es le lemme \ref{segretensor} on peut \'ecrire
 $$ \begin{array}{llcc}
 s_{l}(\mathbf{\Gamma} \boxtimes \S^{\fc{n}{2} + \alpha} Q^*) & =
 \sum_{ \begin{array}{l}
        0 \leq l_0 \leq l_1 \leq l_2 \leq l_3 \leq l \\
        \end{array} } (-1)^{l_1 + l_2} \hspace{0.1cm}
        \binom{\fc{n}{2} + \alpha + l_2}{l_2 - l_1} \hspace{0.1cm}
        \binom{\fc{n}{2} + \alpha + l_1}{l_1 - l_0} \hspace{0.1cm}
        \binom{\fc{n}{2} + \alpha + l - l_2}{l - l_3} \\
        & \hspace{5cm} s_{l_0} (\S^{\fc{n}{2} + \alpha} Q^*)  \cdot
        s_{l_3 - l_2} (\S^{\fc{n}{2} + \alpha} Q^*) \cdot c_1
        (\mathbf{\Gamma})^{l_1 - l_0} \cdot s_{l + l_2 - l_1 - l_3}
        (\mathbf{\Gamma}) \\
        \end{array} $$
 On a $s_k (\S^{\fc{n}{2} + \alpha} Q^*) = (-1)^k \binom{\zeta}{k}
 \tau_{2}^{k} $ où $\zeta$ est le nombre $\binom{\fc{n}{2} + \alpha
 + 1}{2}$. On pose de m\^eme $\zeta^{'} = \binom{\fc{n}{2} - \alpha + 1}{2}$.
 D'apr\`es le lemme \ref{imdirclasseschernsegre}, l'image directe de la classe
 $c_1 (\mathbf{\Gamma})^{l_1 - l_0 + k} \cdot s_{l + l_2 - l_1 - l_3} (\mathbf{\Gamma})$ par
 le morphisme de projection naturel $G \ri \P_{2}^{*}$ est donn\'e par
$$ \sum_{ \begin{array}{ll}
          0 \leq a \leq l_1 - l_0 + k \\
          0 \leq b \leq l + l_2 - l_1 - l_3 \\
          \end{array} }
          \binom{l_1 + k - l_0}{a} \hspace{0.1cm} \left|
 \begin{array}{llcc}
 \binom{\zeta^{'}}{l + l_2 + k + \alpha + 1 - a - b - l_0 - l_3 -
 \fc{n}{2}} & \binom{\zeta^{'}}{a + b + \alpha - \fc{n}{2}} \\
 \binom{\zeta^{'}}{l + l_2 + k + \alpha + 2 - a - b - l_0 - l_3 -
 \fc{n}{2}} & \binom{\zeta^{'}}{a + b + \alpha + 1 - \fc{n}{2}} \\
 \end{array} \right| \hspace{0.1cm} \tau_{1}^{l + l_2 + k + 2 \alpha + 2
 -  l_0 - l_3 - n} $$
 Les d\'eterminants $2 \times 2$ pr\'ec\'edents sont \'egaux \`a
 $$ \fc{l + l_2 + k + 1 - l_0 - l_3 - 2 a - 2 b}{\zeta^{'} + 1}
 \hspace{0.1cm} \binom{\zeta^{'} + 1}{l + l_2 + k + \alpha + 2 - a - b -
 l_0 - l_3 - \fc{n}{2}} \binom{\zeta^{'} + 1}{a + b + \alpha + 1 -
 \fc{n}{2}}$$
 On arrive donc au r\'esultat \'enonc\'e dans le lemme. \end{pv} \\

 Apr\`es r\'eindexation $l:= l$, $m := m - (n + 2 \alpha + 3)$ on obtient donc
 $$ \fbox{$\displaystyle \begin{array}{ll}
 \Delta_{A}^{*} = \fc{(-1)^{n + 2 \alpha + 3} }{\zeta^{'} + 1} \hspace{0.1cm}\sum_{ \begin{array}{llll}
                                     0 \leq l_0 \leq l_1 \leq l_2
                                     \leq l_3 \leq l \leq m
                                     \leq n - 2 \alpha + 2 \\
                                     m-l \leq 4 \\
                                     2 - (m-l) \leq l_0 + (l_3 - l_2) \leq \mathrm{min}(4 -
                                     (m-l) , 2) \\
                                     \end{array} } \\
 (-1)^m \binom{2 n + 5}{n - 2 \alpha + 2 - m} \hspace{0.1cm}
 f_{\alpha,n}(l_0 , l_1 , l_2 , l_3, n - 2 \alpha + 2 - m)
 \left[ \fc{\left[ 1 + (\fc{n}{2} - \alpha) \tau_1 +
                                     (\fc{n}{2} + \alpha) \tau_2
                                     \right]^m \cdot \delta}{1 + u}
                                     \right]^{l_0 + (l_3 - l_2) +
                                     (m-l) - 2}_{2 - l_0 - (l_3 -
                                     l_2)}  \\
                                     \end{array} $} $$
 o\`u la notation $\left[ \gamma(\tau_1 , \tau_2) \right]^{x}_{y}$
 d\'esigne le coefficient se trouvant devant le terme $\tau_{1}^x
 \cdot \tau_{2}^y$ dans le d\'evelopement en s\'erie formelle d'une
 fraction rationelle $\gamma(\tau_1 , \tau_2)$. Notons que pour $\alpha = 1/2, 1$
 on a $\Delta_A = \Delta_{A}^*$, et on pose dans l'expression ci-dessus
 $u = - \tau_2$, $u= \tau_1$ respectivement. \\

  Calculons maintenant $\Delta_{A}^{**}$ (on suppose donc $\alpha \geq 3 / 2)$.
  Ceci correspond \`a la
  contribution \`a $\Delta_{A}$ des termes $- a_{i-j} \tau^{i-j-2}$
  apparaissant dans les classes $(\pi_{-})_{*}(l_{-}^{i-j})$ (cf
  lemme \ref{premiereimagedirecte}). Plus pr\'ecis\'ement on a
$$ \begin{array}{lllcc}
 \Delta_{A}^{**} &  = - 2 (\alpha - 1)^2 \hspace{0.1cm} \sum_{
 \begin{array}{l}
 0 \leq j \leq i \leq 2 n + 4 \\
 \end{array} }
  \binom{2 n + 5}{i+1} \hspace{0.1cm}(-1)^{i+j} \hspace{0.1cm}(i-j)
  \hspace{0.1cm}(i-j-1) \hspace{0.1cm} (2 \alpha -1)^{i-j-2} \\
  & \hspace{6cm} \tau^{i-j-2} \cdot (n \tau - \gamma_{1})^{2 n + 4 - i} \cdot
  s_{j - n - 2 \alpha - 2} (\mathcal{E}^{\vee}) \cap [G] \\
  & = - 2 (\alpha - 1)^2 \hspace{0.1cm} \sum_{ \begin{array}{l}
                                              n + 2 \alpha + 5 \leq j \leq i \leq 2 n  + 5 \\
                                              \end{array} }
                                              \hspace{0.1cm} \binom{2 n + 5}{i}
                                              \hspace{0.1cm}(-1)^{i-j} \hspace{0.1cm} (i-j+2) \hspace{0.1cm}(i-j+1)
                                              \\
         & \hspace{6cm} \left[ (2 \alpha - 1) \tau \right]^{i-j} \cdot (n \tau - \gamma_{1})^{2 n + 5
 - i} \cdot s_{j - n - 2 \alpha - 5} (\mathcal{E}^{\vee}) \cap [G] \\
 & = - 2 (\alpha - 1)^2 \hspace{0.1cm} \left[ (1 + n \tau - \gamma_{1})^{2 n + 5} \cdot \left(
 \sum_k (k+1)(k+2) \left[ (1 - 2 \alpha) \tau \right]^k \right)
 \cdot s (\mathcal{E}^{\vee}) \right]_{n - 2 \alpha} \cap [G] \\
 \end{array} $$
 Posons $\zeta = (1 - 2 \alpha) \tau$, on a
 $$ \begin{array}{l}
 \sum_{k \geq 0} (k+1) (k+2) \zeta^k  = \left( \sum_{k} \zeta^k
 \right)^{(2)} = \left( \fc{1}{1 - \zeta} \right)^{(2)} =
 \fc{2}{(1 - \zeta)^3} \\
 \end{array} $$
o\`u la notation $(\hspace{0.1cm})^{(2)}$ d\'esigne la d\'eriv\'ee
seconde de l'expression entre parenth\`eses. On pose $\delta$ la
classe $(1 + \tau)^{- \binom{2 \alpha + 3}{2}} \cdot (1 - \tau)^{-
\binom{2 \alpha + 1}{2}}$. On a donc
$$ \begin{array}{lllcc}
   - 2 (\alpha - 1)^2 \hspace{0.1cm} (1 + n \tau - \gamma_{1})^{2 n + 5}
   \cdot \fc{2}{\left( 1 + (2 \alpha - 1) \tau \right)^3}
 \cdot s (\mathcal{E}^{\vee}) & = 4 (\alpha - 1)^2 \hspace{0.1cm} \sum_{0 \leq
   m \leq 2 n + 5} \binom{2 n + 5}{m} \hspace{0.1cm} (-1)^m
   \hspace{0.1cm} \left[ \fc{(1 + n \tau)^m \cdot \delta}{\left( 1
   + (2 \alpha - 1) \tau \right)^3 } \right] \\
   & \hspace{3cm} \cdot \gamma_{1}^{2 n +
   5 - m} \cdot \left( \sum_l s_l (\mathbf{\Gamma} \otimes
   \S^{\fc{n}{2} + \alpha} Q^*) \right) \\
   = \lefteqn{ \fc{4 (\alpha - 1)^2}{\zeta^{'} + 1} \hspace{0.1cm}
   \sum_{ \begin{array}{ll}
          0 \leq m \leq 2 n + 5 \\
          \hspace{0.8cm} l \\
          \end{array} } \binom{2 n + 5}{m} \hspace{0.1cm} (-1)^m
   \hspace{0.1cm} g_{\alpha, n}(l , 2 n + 5 - m) \hspace{0.1cm}
   \left[ \fc{(1 + n \tau)^m \cdot \delta}{\left( 1
   + (2 \alpha - 1) \tau \right)^3 } \right] \hspace{0.1cm} \tau^{n + 2 \alpha + 7 - m + l} } \\
 \end{array} $$
 On a introduit la fonction $g_{\alpha, n}$ au lemme
 \ref{secondeimagedirecte}. Pour calculer $\Delta_{A}^{**}$ on ne retient que les couples
 $(m,l)$ tels que $m - (n + 2 \alpha + 7) \leq l \leq m - (n + 2
 \alpha + 5)$; finalement en r\'eindexant la somme ci-dessus par $m
 := m - (n + 2 \alpha + 5)$ et $p := l-m + (n + 2 \alpha + 7)$ on
 obtient
 $$ \fbox{$ \begin{array}{lll}
 \Delta_{A}^{**} = \\
   \\
  \fc{4 (\alpha - 1)^2 \hspace{0.1cm}(-1)^{n + 2 \alpha + 5}}{\zeta^{'} + 1} \hspace{0.1cm}
   \sum_{ \begin{array}{ll}
          0 \leq m \leq n - 2 \alpha \\
          (2 - m)_{+} \leq p \leq 2 \\
          \end{array} } \binom{2 n + 5}{n - 2 \alpha - m} \hspace{0.1cm}
          (-1)^{m} \hspace{0.1cm} g_{\alpha, n}(m + p - 2 , n - 2 \alpha - m)
\hspace{0.1cm}
   \left[ \fc{(1 + n \tau)^m \cdot \delta}{\left( 1
   + (2 \alpha - 1) \tau \right)^3 } \right]_{2 - p} \\
   \end{array} $} $$
 o\`u la notation $\left[ \gamma(\tau) \right]_{x}$ d\'esigne le
 coefficient devant le terme $\tau^x$ dans le d\'eveloppement de la
 fraction rationnelle $\gamma (\tau)$. \\

 \fbox{{\bf{Calculons $\Delta_B$.}}} On a encore un morphisme
 canonique $B^{'} \ri B \times_{\Sigma} \Sigma_{+}$, qui est
 birationnel sur son image $B^{''} \subseteq B \times_{\Sigma} \Sigma_{+}$
 qui est un ferm\'e int\`egre de codimension $2 \alpha - 4$. En
 effet la dimension du produit fibr\'e est $n - 2 + n + 2 \alpha + 2
 = 2 n + 2 \alpha$. Lorsque $\alpha = 2$, le ferm\'e $B^{''}$ est tout
 le produit fibr\'e. Lorsque $\alpha = 5 / 2$ c'est une hypersurface.
 Nous d\'esirons conna\^itre en toute g\'en\'eralit\'e la classe
 fondamentale $[B^{''}]$ du ferm\'e $B^{''}$: celle-ci est bien d\'efinie car
 le produit fibr\'e est lisse.

 \begin{prop}
 \label{classefond}Avec les notations ci-dessus, on a
 $$ [B^{''}] = [\fc{c (L_{-}^{\otimes 2}
 \otimes L_{+} \otimes \S^{2 \alpha - 3} Q^*)}{c (L_{-} \otimes L_{+}
 \otimes Q(1))}]_{2 \alpha - 4} $$
 \end{prop}

 \begin{pv} Soient $\mathbf{\Lambda}^{'}$ et $\mathbf{\Lambda}^{''}$
 les familles de syst\`emes coh\'erents param\'etr\'ees par la
 vari\'et\'e $\Sigma$ (cf notations de la proposition \ref{Sigma}).
 On a une suite exacte de familles
$$ 0 \ri L_{-} \otimes \mathbf{\Lambda}^{'} \ri \mathbf{\Lambda}
\ri \mathbf{\Lambda}^{''} \ri 0 \hspace{2cm} (+)$$
 o\`u $\mathbf{\Lambda}$ est une famille param\'etrant les
 syst\`emes coh\'erents dont la classe appartient \`a $\Sigma_{-}$.
 On restreint cette suite exacte \`a la composante $B$.
 La suite exacte (+) d\'efinit une filtration de $\mathbf{\Lambda}$
 et il en r\'esulte une suite exacte longue de faisceaux
 $\underline{\Ext}$ relatifs filtr\'es (voir [L] chap. 15.3 ou [HL]
 annexe du chap. 2 pour les d\'efinitions)
 $$ \begin{array}{ll}
 0 \ri \underline{\Ext}_{p,\hspace{0.1cm}-}^{1} (\mathbf{\Lambda},
 \mathbf{\Lambda}) \ri \underline{\Ext}_{p}^{1} (\mathbf{\Lambda},
 \mathbf{\Lambda}) \ri \underline{\Ext}_{p,\hspace{0.1cm}+}^{1} (\mathbf{\Lambda},
 \mathbf{\Lambda}) \ri \\
  \hspace{1cm} \underline{\Ext}_{p,\hspace{0.1cm}-}^{2} (\mathbf{\Lambda},
 \mathbf{\Lambda}) \ri \underline{\Ext}_{p}^{2} (\mathbf{\Lambda},
 \mathbf{\Lambda}) \ri \underline{\Ext}_{p,\hspace{0.1cm}+}^{2} (\mathbf{\Lambda},
 \mathbf{\Lambda}) \ri 0 \hspace{1cm} (++) \\
 \end{array}$$
 Le morphisme $p$ dans la suite exacte qui pr\'ec\`ede est
 la premi\`ere projection $B \times \P_2 \ri B$. On a $\underline{\Ext}_{p,\hspace{0.1cm}+}^{2} (\mathbf{\Lambda},
 \mathbf{\Lambda}) =0$, d'apr\`es la suite spectrale classique aboutissant vers
 un groupe $\Ext$ filtr\'e.
 Le faisceau $\underline{\Ext}^{2}_{p} (\mathbf{\Lambda},
 \mathbf{\Lambda})$ est de support le ferm\'e $\Sigma_{-}^{'}$,
 en utilisant les notations de la proposition \ref{fermesing}.
 On a d'autre part $\underline{\Ext}^{1}_{p, \hspace{0.1cm} +} (\mathbf{\Lambda},
 \mathbf{\Lambda}) \simeq L_{-}^{\vee} \otimes
 \underline{\Ext}^{1}_{p}(\mathbf{\Lambda}^{'},
 \mathbf{\Lambda}^{''})$, et le morphisme $\underline{\Ext}_{p,\hspace{0.1cm}-}^{1} (\mathbf{\Lambda},
 \mathbf{\Lambda}) \ri \underline{\Ext}_{p}^{1} (\mathbf{\Lambda},
 \mathbf{\Lambda})$ est celui de l'inclusion du faisceau tangent \`a
 $\Sigma_{-}$ dans le faisceau tangent de $\mathcal{S}_{-}$
 restreint \`a $B$ (voir proposition \ref{espacetangentsigma-}).
 D'apr\`es le th\'eor\`eme \ref{critlissite} et la proposition
 \ref{fermesing}, le faisceau
 $\underline{\Ext}_{p}^{2} (\mathbf{\Lambda},
 \mathbf{\Lambda})$ est de support le ferm\'e $\Sigma_{-}^{'}$ qui
 est de codimension $\geq 2$ dans $B$. Notons $W$ le faisceau
 image du morphisme $\underline{\Ext}_{p,\hspace{0.1cm}+}^{1} (\mathbf{\Lambda},
 \mathbf{\Lambda}) \ri \underline{\Ext}_{p,\hspace{0.1cm}-}^{2} (\mathbf{\Lambda},
 \mathbf{\Lambda})$; on a donc une suite exacte
$$ 0 \ri W \ri \underline{\Ext}_{p, \hspace{0.1cm} -}^{2} (\mathbf{\Lambda},
 \mathbf{\Lambda})  \ri \underline{\Ext}_{p}^{2} (\mathbf{\Lambda},
 \mathbf{\Lambda}) \ri 0 $$
 Appliquons alors le foncteur
 $\underline{\Hom}_{\O_{B}}(\Box , \O_{B})$ \`a cette suite exacte,
 comme le faisceau $\underline{\Ext}_{p}^{2} (\mathbf{\Lambda},
 \mathbf{\Lambda})$ a un support de codimension $\geq 2$ dans $B$
 qui est lisse, on a $W^{\vee} \simeq
 \underline{\Ext}_{p, \hspace{0.1cm} -}^{2} (\mathbf{\Lambda},
 \mathbf{\Lambda})^{\vee}$. On se place sur le projectif au sens de
 Grothendieck $\mathbb{P} (\underline{\Ext}_{p,\hspace{0.1cm}+}^{1} (\mathbf{\Lambda},
 \mathbf{\Lambda})^{\vee})$, qui est isomorphe \`a $B
 \times_{\Sigma} \Sigma_{+}$; sur cette vari\'et\'e on a d'apr\`es
 ce qui pr\'ec\`ede un morphisme $W^{\vee} \ri L_{-} \otimes L_{+}$,
 autrement dit une section du faisceau
 $\underline{\Ext}_{p, \hspace{0.1cm} -}^{2} (\mathbf{\Lambda},
 \mathbf{\Lambda})^{\vee \vee} \otimes L_{-} \otimes L_{+}$ dont le
 sch\'ema des z\'eros est pr\'ecis\'ement le ferm\'e $B^{''}$.

 \begin{lem}
 \label{faisceauExt-}Le faisceau $\underline{\Ext}_{p, \hspace{0.1cm} -}^{2} (\mathbf{\Lambda},
 \mathbf{\Lambda})$ sur $B$ est conoyau d'un morphisme de fibr\'es
 vectoriels
 $$\imath : Q(1) \ri \S^{2 \alpha - 3} Q^* \otimes L_{-}$$
  qui est
 g\'en\'eriquement injectif. Son lieu singulier est l'intersection
 $B \cap A$ et il est localement libre de rang $2 \alpha - 4$ sur $B
 \setminus A$. De plus ce faisceau est sans torsion lorsque $\alpha
 = 5 / 2$ et r\'eflexif lorsque $\alpha > 5 / 2$.
 \end{lem}

 \begin{pv}
 Ceci vient de la suite spectrale convergeant vers
 $\underline{\Ext}_{p, \hspace{0.1cm} -}^{2} (\mathbf{\Lambda},
 \mathbf{\Lambda})$: en effet ce dernier est conoyau de la
 diff\'erentielle
 $$ d^{0,1}: E_{1}^{0,1} =
 \underline{\Ext}_{p}^{1} (\mathbf{\Lambda}^{'},
 \mathbf{\Lambda}^{'}) \oplus
 \underline{\Ext}_{p}^{1} (\mathbf{\Lambda}^{''},
 \mathbf{\Lambda}^{''}) \ri E_{1}^{1,1} = L_{-} \otimes
 \underline{\Ext}_{p}^{2} (\mathbf{\Lambda}^{''},
 \mathbf{\Lambda}^{'})$$
 Notons que le noyau de $d^{0,1}$ est le faisceau
 $\underline{\Ext}_{p, \hspace{0.1cm} -}^{1} (\mathbf{\Lambda},
 \mathbf{\Lambda})$, qui est isomorphe au faisceau tangent de $B$.
 L'inclusion $Ker(d^{0,1}) \subseteq E_{1}^{0,1}$ s'identifie \`a la
 diff\'erentielle $d \pi_{-}$ de la projection $\pi_{-}: \Sigma_{-}
 \ri \Sigma$, restreinte \`a la composante $B$. L'image de $d
 \pi_{-}$ est \'egale en tout point de $B \setminus A$ \`a l'espace
 tangent de la sous-vari\'et\'e de $\Sigma$ d'id\'eal $\mathcal{J}$
 (voir d\'efinition \ref{idealdiagonale}), et en chaque point de $A
 \cap B$ contient un tel espace tangent.

  C'est donc que $d^{0,1}$ se factorise en un morphisme injectif de faisceaux
  $\imath : Q(1) \ri \S^{2 \alpha - 3} Q^* \otimes L_{-}$,
  car $Q(1)$ est le fibr\'e normal $\mathcal{J} / \mathcal{J}^2$, et
  $\underline{\Ext}_{p, \hspace{0.1cm} -}^{2}
  (\mathbf{\Lambda}, \mathbf{\Lambda})$ est le conoyau de $\imath$.
  De plus $\imath$ est un morphisme injectif au sens des fibr\'es sur
  $B \setminus A$, ce qui prouve le lemme: en effet quand $\alpha = 5 / 2$ le ferm\'e $A \cap B$
  est de codimension 2 dans $B$ qui est lisse, et quand $\alpha > 5 / 2$
  cette codimension est $\geq 3$. \end{pv} \\

  Revenant \`a la preuve de la proposition \ref{classefond},
 le faisceau $\underline{\Ext}_{p, \hspace{0.1cm} -}^{2} (\mathbf{\Lambda},
 \mathbf{\Lambda})^{\vee \vee}$ co\"incide toujours avec
 $\underline{\Ext}_{p, \hspace{0.1cm} -}^{2} (\mathbf{\Lambda},
 \mathbf{\Lambda})$ sur $B \setminus A$ dont le compl\'ementaire
 dans $B$ est de codimension $2 \alpha - 3 > 2 \alpha - 4$,
 c'est donc que la classe fondamentale de $B^{''}$ est \'egale \`a
 $$ [B^{''}] = c_{2 \alpha - 4} (\underline{\Ext}_{p, \hspace{0.1cm} -}^{2} (\mathbf{\Lambda},
 \mathbf{\Lambda}) \otimes L_{-} \otimes L_{+})$$
 ce qui termine la preuve de la proposition compte tenu du lemme. \end{pv} \\

  Revenons au calcul de $\Delta_B$. On se doit de calculer pour tout
  couple $(i,j)$ d'entiers tels que $0 \leq i \leq 2 n + 4$, $0 \leq
  j \leq i$, le nombre d'intersection $v^{2 n + 4 - i} \cdot l_{+}^j
  \cdot l_{-}^{i-j} \cap [B^{'}]$. On a un diagramme
  $$ \xymatrix{
     B^{'} \ar[r] & B^{''} \ar@{^{(}->}[d] & \\
     & B \times_{\Sigma} \Sigma_{+} \ar[dl]_{\varphi_{-}}
     \ar[dr]^{\varphi_{+}} & \\
     B \ar[dr]_{\pi_{-}} & & \Sigma_{+} \ar[dl]^{\pi_{+}} \\
     & \Sigma & \\
     } $$
 o\`u la fl\`eche horizontale est birationnelle. Du fait que sur
 $B$, $B^{'}$ et $B^{''}$ on a $\tau_1 = \tau_2 = \tau$, on peut \'ecrire
 en notation additive $v = n \tau - \gamma_1$. Posons
 $\mathcal{F}$ le faisceau coh\'erent $\underline{\Ext}_{p, \hspace{0.1cm} -}^{2} (\mathbf{\Lambda},
 \mathbf{\Lambda})$ intervenant dans le lemme \ref{faisceauExt-}.
 On a donc
 $$ \Delta_B = \sum_{0 \leq j \leq i \leq 2 n + 4}
 \binom{2 n + 5}{i+1} \hspace{0.1cm} (-1)^{i-j} \hspace{0.1cm}
 v^{2 n + 4 - i} \cdot l_{-}^{i-j} \cdot l_{+}^j \cdot
 c_{2 \alpha - 4}(\mathcal{F} \otimes L_{-} \otimes L_{+}) \cap
 [B \times_{\Sigma} \Sigma_{+}] $$
 Sur $B \setminus A$ le faisceau coh\'erent $\mathcal{F}$ est localement libre
 de rang $2 \alpha - 4$, on peut donc \'ecrire d'apr\`es [F] ex.3.2.2
 $$ c_{2 \alpha - 4}(\mathcal{F} \otimes L_{-} \otimes L_{+}) =
    \sum_{0 \leq k \leq 2 \alpha - 4} c_{2 \alpha - 4 -
    k}(\mathcal{F} \otimes L_{-}) \cdot l_{+}^k $$
 Comme au paragraphe pr\'ec\'edent on r\'eindexe la somme en posant
 $i := i+1$, $j :=j+1$, ce qui donne
 $$ \Delta_B = \sum_{ \begin{array}{ll}
                                 1 \leq j \leq i \leq 2 n + 5 \\
                                 0 \leq k \leq 2 \alpha - 4 \\
                                 \end{array} }
 \binom{2 n + 5}{i} \hspace{0.1cm} v^{2 n + 5 - i} \cdot (-1)^{i-j} \hspace{0.1cm}
l_{-}^{i-j} \cdot c_{2 \alpha - 4 - k}(\mathcal{F} \otimes L_{-})
\cdot l_{+}^{j + k - 1} \cap [B \times_{\Sigma} \Sigma_{+}] $$
 Les classes $v^{2 n + 5 - i} \cdot l_{-}^{i-j} \cdot
 c_{2 \alpha - 4 - k}(\mathcal{F} \otimes L_{-})$ proviennent du
 ferm\'e $B$ qui est de dimension $n-2$; elles ne sont non nulles
 que lorsque l'in\'egalit\'e suivante sur les indices est
 v\'erifi\'ee
 $$ (2 n + 5 - i) + (i-j) + (2 \alpha - 4 - k) =
 2 n + 2 \alpha + 1 - (j + k) \leq n - 2 $$
 c'est-\`a-dire $j+k \geq n + 2 \alpha + 3$. L'exclusion des indices $(j \leq i)$
 tels que $i=0$ ou $i > 0$ et $j=0$ dans la somme pr\'ec\'edente n'est donc que
 factice, et on peut \'ecrire, apr\`es avoir appliqu\'e
 {\bf{seulement}} l'image directe $(\pi_{+})_{*}$ aux termes $l_{+}^{j+ k - 1}$
 $$ \Delta_{B} = \left[ (1 + v)^{2 n + 5} \cdot \fc{c (\mathcal{F} \otimes L_{-})}{1 + l_{-}} \cdot \delta \cdot
 s (\Gamma \otimes \S^{\fc{n}{2} + \alpha} Q^*) \right]_{n-2} \cap
 [B] $$
 o\`u $\delta$ d\'esigne la classe $(1 + \tau)^{- \binom{2 \alpha +
 3}{2}} \cdot (1 - \tau)^{- \binom{2 \alpha + 1}{2}}$. De la suite
 exacte suivante sur $B \setminus A$
$$ 0 \ri Q(1) \otimes L_{-} \ri \S^{2 \alpha - 3} Q^* \otimes
L_{-}^{\otimes 2} \ri \mathcal{F} \otimes L_{-} \ri 0 $$
 On peut donc \'ecrire
 $$ \begin{array}{lllcc}
  c (\mathcal{F} \otimes L_{-}) & = \fc{c (\S^{2 \alpha - 3} Q^*
  \otimes L_{-}^{\otimes 2})}{c(Q(1) \otimes L_{-})} \\
  & = \fc{(1 + 2 l_{-})^{\binom{2 \alpha - 1}{2}}}{(1 + \tau + 2
  l_{-})^{\binom{2 \alpha - 2}{2}} \cdot c (Q(1) \otimes L_{-})} \\
  & = \fc{(1 + 2 l_{-})^{\binom{2 \alpha - 1}{2}} \cdot (1 +
  l_{-})}{(1 + \tau + 2 l_{-})^{\binom{2 \alpha - 2}{2}} \cdot
  (1 + \tau + l_{-})^{3}} \\
  \end{array} $$
  On cherche \`a pr\'esent \`a d\'evelopper la classe
$\fc{(1 + 2 l_{-})^{\binom{2 \alpha - 1}{2}}}{(1 + \tau + 2
l_{-})^{\binom{2 \alpha - 2}{2}} \cdot
  (1 + \tau + l_{-})^{3}}$ sous forme d'une s\'erie en $l_{-}$; on
  obtient
  $$ \fc{(1 + 2 l_{-})^{\binom{2 \alpha - 1}{2}}}{(1 + \tau + 2
l_{-})^{\binom{2 \alpha - 2}{2}} \cdot
  (1 + \tau + l_{-})^{3}} = \fc{1}{(1 + \tau)^{3 + \binom{2 \alpha
  - 2}{2}}} \cdot \fc{(1 + 2 l_{-})^{\binom{2 \alpha -
  1}{2}}}{\left[ 1 + \left( \fc{l_{-}}{1 + \tau} \right) \right]^3
  \cdot \left[ 1 + \left( \fc{2 l_{-}}{1 + \tau} \right)
  \right]^{\binom{2 \alpha - 2}{2}}} $$
  ceci \'etant encore \'egal \`a
$$ \fc{(1 + 2 l_{-})^{\binom{2 \alpha -
  1}{2}}}{(1 + \tau)^{3 + \binom{2 \alpha
  - 2}{2}}} \cdot \left( \sum_{p \geq 0} (-1)^p \hspace{0.1cm}
  \binom{2 + p}{2} \hspace{0.1cm} \fc{l_{-}^{p}}{(1 + \tau)^p}
  \right) \hspace{0.1cm} \left( \sum_{m \geq 0} (-1)^m
  \hspace{0.1cm} \binom{\binom{2 \alpha - 2}{2} - 1 + m}{\binom{2
  \alpha - 2}{2} -1} \hspace{0.1cm} \fc{2^m \hspace{0.1cm}
  l_{-}^m}{(1 + \tau)^m} \right)  $$
c'est-\`a-dire, en posant $k := m+p$ et $l := m$ dans les sommes
ci-dessus
$$ \begin{array}{ll}
 \fc{(1 + 2 l_{-})^{\binom{2 \alpha - 1}{2}}}{(1 + \tau)^{3 + \binom{2 \alpha
  - 2}{2}}} \cdot \sum_{0 \leq l \leq k } (-1)^k \hspace{0.1cm}
  \binom{2 + k - l}{2}
\hspace{0.1cm} \binom{\binom{2 \alpha - 2}{2} - 1 + l}{\binom{2
\alpha - 2}{2} - 1} \hspace{0.1cm} 2^l \hspace{0.1cm}
\fc{l_{-}^k}{(1 + \tau)^k} = \\
 \fc{1}{(1 + \tau)^{3 + \binom{2 \alpha
  - 2}{2}}} \cdot \sum_{ 0 \leq l \leq k \leq m}
  (-1)^k \hspace{0.1cm} \binom{2 + k - l}{2} \hspace{0.1cm} \binom{\binom{2 \alpha - 2}{2} - 1 + l}{\binom{2
\alpha - 2}{2} - 1} \hspace{0.1cm} \binom{\binom{2 \alpha -
1}{2}}{m-k} \hspace{0.1cm} 2^{m + l-k}  \hspace{0.1cm}
\fc{l_{-}^m}{(1 + \tau)^k} \hspace{2cm} {\bf{(1)}} \\
\end{array} $$
 On a d'autre part
 $$ (\pi_{-})_* (l_{-}^m) = s_{m - 2 \alpha + 2}
 (\S^{2 \alpha - 2} Q(1)) = \binom{\binom{2 \alpha}{2} + m - 2 \alpha + 1}{\binom{2 \alpha}{2}
- 1} \hspace{0.1cm} \tau^{m - 2 \alpha + 2} $$
 Pour que cette classe soit non nulle, on doit donc se restreindre
 \`a $2 \alpha - 2 \leq m \leq 2 \alpha$. Enfin en posant $\rho: G
 \ri \P_{2}^*$ le morphisme structural du fibr\'e en grassmaniennes
 (cf lemme \ref{secondeimagedirecte}) on a
 $$ \begin{array}{lll}
  \rho_* \left( (1 + n \tau - \gamma_1)^{2 n + 5} \cdot s(\Gamma
  \otimes \S^{\fc{n}{2} + \alpha} Q^*) \right) = \\
  \sum_{ \begin{array}{ll}
 n - 2 \alpha - 2 \leq p + l \leq n - 2 \alpha \\
 p \hspace{0.1cm}, \hspace{0.1cm} l \geq 0
  \end{array} } (-1)^p \hspace{0.1cm} \binom{2 n + 5}{p} \hspace{0.1cm}
  (1 + n \tau)^{2 n + 5 - p} \cdot \rho_* \left( \gamma_{1}^p \cdot s_l
  (\Gamma \otimes \S^{\fc{n}{2} + \alpha} Q^*) \right) = \\
   \hspace{0.1cm} \\
  \sum_{ \begin{array}{ll}
 n - 2 \alpha - 2 \leq p + l \leq n - 2 \alpha \\
  p \hspace{0.1cm}, \hspace{0.1cm}l \geq 0 \\
  \end{array} } \fc{(-1)^p}{\zeta^{'} + 1} \hspace{0.1cm} \binom{2 n + 5}{p} \hspace{0.1cm}
  g_{\alpha, n}(l,p) \hspace{0.1cm} (1 + n \tau)^{2 n + 5 - p} \cdot
 \tau^{p + l + 2 \alpha + 2 - n} \\
 \end{array} $$
  L'expression pr\'ec\'edente se r\'eecrit, avec la r\'eindexation
  $p:=p$; \hspace{0.2cm}$q:=p + l - (n - 2 \alpha - 2)$
  $$ \fc{1}{\zeta^{'} + 1}
  \hspace{0.1cm} \sum_{ \begin{array}{ll}
                        0 \leq q \leq 2 \\
                        0 \leq p \leq q + (n - 2 \alpha - 2) \\
                        \end{array} } \hspace{0.1cm} (-1)^p \hspace{0.1cm}
                        \binom{2 n + 5}{p} \hspace{0.1cm} g_{\alpha,n}(
                        q + n - 2 \alpha - 2 - p, p)
  \hspace{0.1cm}(1 + n \tau)^{2 n + 5 - p} \cdot \tau^q
  \hspace{2cm}{\bf{(2)}}$$
 On r\'eindexe la somme dans l'expression {\bf{(1)}} en posant $m:=m
 - (2 \alpha - 2)$; en combinant avec l'expression {\bf{(2)}} on
 obtient
 $$ \fbox{$\displaystyle \begin{array}{llcc}
 \Delta_B & = \fc{2^{2 \alpha - 2}}{\zeta^{'} + 1} \hspace{0.1cm}
 \sum_{ \begin{array}{lll}
        0 \leq m,q \hspace{0.1cm}; \hspace{0.1cm}m + q \leq 2 \\
        0 \leq l \leq k \leq m + 2 \alpha - 2 \\
        0 \leq p \leq q + n - 2 \alpha - 2 \\
        \end{array} } \hspace{0.2cm} (-1)^{k + p}
        \hspace{0.1cm}\binom{2 + k - l}{2}
        \hspace{0.1cm}\binom{\binom{2 \alpha - 2}{2} - 1 + l}{\binom{2 \alpha - 2}{2} -
        1} \hspace{0.1cm}\binom{\binom{2 \alpha - 1}{2}}{m - k + 2
        \alpha - 2} \\
        & \hspace{6cm} \binom{\binom{2 \alpha}{2} + m - 1}{\binom{2 \alpha}{2}
        -1} \hspace{0.1cm}\binom{2 n + 5}{p} \hspace{0.1cm}2^{m + l
        - k} \hspace{0.1cm} g_{\alpha, n}(q - p + n - 2 \alpha - 2,
        p) \\
        & \hspace{6cm} \left[ \fc{(1 + n \tau)^{2 n + 5 - p}}{(1 -
        \tau)^{\binom{2 \alpha + 1}{2}}} \cdot \fc{1}{( 1 + \tau)^{k
        + 3 + 2 (2 \alpha^2 + 3)}} \right]_{2 - m - q} \\
        \end{array} $} $$
o\`u la notation $\left[ f(\tau) \right]_{x}$ d\'esigne la partie de
degr\'e $x$ dans le d\'eveloppement en s\'erie formelle d'une
fraction rationnelle $f(\tau)$.

\subsection{Quelques applications num\'eriques}
\label{appnumDelta}
 Dans ce paragraphe nous appliquons les
 formules pr\'ec\'edentes. Les r\'esultats sont not\'es dans le
 tableau ci-dessous.
 Nous n'avons bien sur pas pu nous
 passer d'un logiciel de calcul, en l'occurrence Maple.
 Nous faisons d'abord quelques commentaires.
 \\

  Pour ${\bf{n=4}}$: il n'y a qu'un seul saut \`a calculer
  \`a la valeur $\alpha = 1$. Le ferm\'e $\Sigma_{-}$
  est lisse ir\'eductible, isomorphe \`a $\P_{2}^{*} \times \P_{2}^{*}$,
  et le saut $\Delta$ est \'egal \`a $\Delta_{A} = \Delta_{A}^*$.
  D'apr\`es ce qui pr\'ec\`ede, on a
  $$ \begin{array}{cc}
  \Delta_{A}^* & = \frac{-1}{2} \,
  \sum_{ \begin{array}{ll}
         0 \leq l_0 \leq l_1 \leq l_2 \leq l_3 \leq
         l \leq m \leq 4 \\
     2 - (m-l) \leq l_0 + (l_3 - l_2) \leq min \{
         4 - (m-l) , 2 \} \\
         \end{array} }
         (-1)^m \, \binom{13}{4 - m} \,
         f_{1,4} (l_0 , l_1 , l_2 , l_3 , 4 - m) \\
 & \hspace{10cm} \left[ \frac{(1 + \tau_1 + 3 \tau_2 )^m \cdot
         (1 + \tau_1 - \tau_2)^3 }{(1 + \tau_1)^{11} \cdot
         (1 - \tau_2)^3 } \right]_{2 - l_0 - (l_3 - l_2)}^{l_0
         + (l_3 - l_2) + (m-l) - 2} \\
         \end{array} $$
 o\`u la fonction $f_{1,4} (l_0 , l_1 , l_2 , l_3, 4 - m)$
 d\'epend du param\`etre entier $l$ et a pour expression celle
 donn\'ee au lemme \ref{secondeimagedirecte}, et est une somme
 de $(l_1 - l_0 + 5 -m) \, (l - l_3 + l_2 - l_1 + 1)$ termes.
 Le calcul de $\Delta$ requiert ici environ une centaine
 de sommations.
 \\

 Pour ${\bf{n=5}}$: il y a deux sauts \`a calculer aux valeurs
 $\alpha = \frac{1}{2}$ et $\alpha = \frac{3}{2}$;
 \begin{itemize}
 \item pour $\alpha
 = \frac{1}{2}$, le ferm\'e $\Sigma_{-}$ est isomorphe \`a
 $\Sigma = \mathrm{Grass}(2 , \S^2 Q) \times \P_{2}^{*}$, qui
 est une fibration en plans projectifs sur $\P_{2}^{*}
 \times \P_{2}^{*}$, et on a encore
 $\Delta = \Delta_{A} = \Delta_{A}^*$ avec $u = - \tau_2$;
 \item pour
 $\alpha = \frac{3}{2}$, le ferm\'e $\Sigma_{-}$ est isomorphe
 \`a l'\'eclat\'e de $\P_{2}^* \times \P_{2}^*$ le long de la
 diagonale, et on a
 $\Delta = \Delta_A = \Delta_{A}^* + \Delta_{A}^{**}$, avec
 $$ \begin{array}{llcc}
\Delta_{A}^* & = - \frac{1}{4}
 \sum_{ \begin{array}{lll}
        0 \leq l_0 \leq l_1
        \leq l_2 \leq l_3 \leq l \leq m \leq 6 \\
        m - l \leq 4 \\
        2 - (m-l) \leq l_0 + (l_3 - l_2 ) \leq
        min \{ 2, 4 - (m-l) \} \\
        \end{array} }
        (-1)^m \, \binom{15}{6-m} \,
        f_{\frac{1}{2} , 5} (l_0 , l_1 , l_2 , l_3 , 6 - m)
        \, \\
        & \hspace{10cm} \left[
        \frac{ (1 + 2 \tau_1 + 3 \tau_2 )^m \cdot
        (1 + \tau_1 - \tau_2)^3 }{(1 + \tau_1 )^6 \cdot
        (1 - \tau_2 )^2 }
        \right]_{2 - l_0 - (l_3 - l_2)}^{l_0 + (l_3 - l_2) + (m-l)
        -2} \\
 \Delta_{A}^{**} & = - \frac{1}{2} \sum_{
   \begin{array}{ll}
   0 \leq m \leq 2 \\
   2 - m \leq p \leq 2 \\
   \end{array}  } \binom{15}{2-m} \, (-1)^m \,
   g_{\frac{3}{2}, 5} (m + p - 2 , 2 -m) \,
   \left[ \frac{(1 + 5 \tau)^m }{(1 + 2 \tau)^3 (1 + \tau)^{15}
   (1-\tau)^6 } \right]_{2-p} \\
   \end{array} $$
 \end{itemize}
 Nous n'avons pas calcule l\`a encore le nombre exact de sommations
 qui apparaissent dans $\Delta_A$,
n\'eanmoins nous l'avons \'evalu\'e \`a $300$ pour $\Delta_{A}^*$
et moins de $10$ pour $\Delta_{A}^{**}$, en raison du fait
que c'est le premier cas o\`u ce terme suppl\'ementaire intervient. \\

Pour ${\bf{n=6}}$: il y a \'egalement deux sauts \`a calculer aux
valeurs $\alpha=1$ et $\alpha = 2$;
\begin{itemize}
\item pour $\alpha = 1$ le ferm\'e $\Sigma_{-}$ est encore
isomorphe \`a $\mathrm{Grass}(2 , \S^2 Q) \times \P_{2}^{*}$ et
\`a cette valeur critique on a $\Delta = \Delta_{A} =
\Delta_{A}^*$ avec $u = \tau_1$; \item pour $\alpha=2$ le ferm\'e
$\Sigma_{-}$ n'est plus irr\'eductible (c'est le premier cas
ainsi), et il contient l'\'eclat\'e $A$ de $\P_{2}^* \times
\P_{2}^*$ le long du carr\'e de l'id\'eal de la diagonale, ainsi
que la composante $B$ isomorphe au fibr\'e projectif $\P (\S^2 Q)
\ri \P_{2}^*$, tous deux de dimension $4$. On a dans ce cas
 $\Delta = \Delta_A + \Delta_B$, avec
 $\Delta_{A} = \Delta_{A}^* + \Delta_{A}^{**} $:
 \begin{enumerate}
 \item on pose
 $\zeta^{'}=1$, $u= 3 \tau_1 + 2 \tau_2$, $\delta =
 \frac{(1 + \tau_1 - \tau_2 )^{15} }{(1 + \tau_1 )^{21}
 \cdot (1 - \tau_2 )^{10} }$ dans le calcul de $\Delta_{A}^*$,
 et $\delta = (1 + \tau_1 )^{- 21}
 \cdot (1 - \tau_2 )^{- 10}$ dans celui de $\Delta_{A}^{**}$;
 \item pour $\Delta_B$ en appliquant la formule pr\'ec\'edente
 \`a $n=6$ et $\alpha = 2$ on obtient
 $$ \begin{array}{llcc}
  \Delta_B & = 2 \sum_{ \begin{array}{ll}
   m + q \leq 2 \\
   l \leq k \leq m + 2 \\
   p \leq q \\
   \end{array} } \, (-1)^{k + p} \, \binom{2 + k - l}{2}
   \, \binom{3}{m - k + 2} \, \binom{5 + m}{5} \,
   \binom{17}{p} \, 2^{m + l - k} \,
   g_{2,6} (q-p , p) \, \\
   & \hspace{5cm} \left[
   \frac{(1 + 5 \tau)^{17-p} }{(1-\tau)^{10} } \cdot
   \frac{1}{(1 + \tau)^{k + 25} } \right]_{2 - m - q} \\
   \end{array} $$
 \end{enumerate}
 Les indices intervenant dans les sommes pr\'ec\'edentes
 sont tous positifs. Dans ce cas $n=6$ le nombre total
 de sommations n\'ecessaires
 au calcul des sauts avoisine le millier.
\end{itemize}

\begin{center}
 \begin{table} \caption{Sauts pour $n \leq 6$, $\alpha \leq 2$}
 \begin{tabular}{|r|l|l|l|} \hline
 $\alpha$ & $n=4$ & $n=5$ & $n=6$ \\
 \hline\hline
 $\frac{1}{2}$ & $\times$ &
 $13347$ & $\times$ \\
 \hline
 $1$ & $- 222$ & $\times$ &
 $- 444 $ \\
 \hline
 $\frac{3}{2}$ & $\times$ &
 $- 1332$ & $\times$
 \\
 \hline
 $2$ & $\times$ &
 $\times$ & $213 $ \\
 \hline
 \end{tabular}
 \end{table}
 \end{center}

\section{Calcul en dessous de la plus petite valeur critique}

 Soit toujours un entier $n \geq 4$.
 Il nous reste \`a pr\'esent \`a calculer le nombre
 $c_{1} (\mathcal{D}_{-})^{2 n + 5} \cap [\mathcal{S}_{-}]$,
 o\`u $\mathcal{D}_{-} = \mathcal{D}_{\alpha_{-}}$ et
 $\mathcal{S}_{-} = \mathcal{S}_{\alpha_{-}}$, avec $\alpha_{-}$ un
 nombre rationnel strictement positif juste en dessous de la
 plus petite valeur critique $\alpha$. Rappelons que celle-ci
 est \'egale \`a 1/2 lorsque $n$ est impair, et \`a 1 sinon (voir
 proposition \ref{Valeurcrit}).
  Dans le cas $n$ impair on
 peut d\'ecrire explicitement $\mathcal{S}_{-}$ comme grassmanienne
 relative sur l'espace des coniques $\P_5$, mais pour $n$ pair
 ce n'est plus possible. Il existe seulement une grassmanienne
 relative de dimension $2 n + 6$ {\bf{dominant}} $\mathcal{S}_{-}$,
 qui est rappelons le de dimension $ 2n + 5$. Les fibres de ce morphisme
 dominant sont
 g\'en\'eriquement des coniques. \\

  Avant d'entamer ces deux parties nous commen\c{c}ons par un lemme.

 \begin{lem}
 \label{epsilonsstab}Soit $(\Gamma, \Theta)$ un syst\`eme coh\'erent
 $\alpha_{-}-$semi-stable, alors $\Theta$ est un faisceau
 semi-stable.
 \end{lem}

 \begin{pv}
 En effet si $\Theta^{'} \subset \Theta$ est un sous-faisceau coh\'erent
 de multiplicit\'e 1 v\'erifiant
 $\chi(\Theta^{'}) > \fc{\chi(\Theta)}{2}$, alors $\chi(\Theta^{'})
 \geq \fc{\chi(\Theta)}{2} + \fc{1}{2}$. Posant $\Gamma^{'} =
 \Gamma \cap \H^0 (\Theta^{'})$, pour un param\`etre $0 < \epsilon
 < \fc{1}{2}$, on a
 $$ \epsilon \hspace{0.1cm} \dim \Gamma^{'} + \chi(\Theta^{'}) >
 \epsilon + \fc{\chi(\Theta)}{2} $$
 ce qui contredit la $\epsilon-$semi-stabilit\'e. \end{pv} \\

  L'implication inverse est vraie lorsque $n$ est impair,
  mais fausse quand $n$ est pair comme en t\'emoigne ce qui suit.

 \subsection{Le cas pair}

 Soit donc $\alpha_{-} < \alpha$ une valeur positive inf\'erieure \`a la
 plus petite valeur critique.
 Soit $\mathcal{C} \subset \P_5 \times \P_2$ la conique universelle
  param\'etrant les couples $(C,x)$ avec $C$ une conique et $x$ un
  point de $C$. Le produit fibr\'e $\mathcal{C} \times_{\P_5}
  \mathcal{C}$ param\`etre les coniques avec deux points marqu\'es, et
  on note $\mathcal{I}$ le faisceau d'id\'eaux du plongement diagonal
  $\mathcal{C} \subset \mathcal{C} \times_{\P_5}
  \mathcal{C}$. On note $p: \mathcal{C} \times_{\P_5}
  \mathcal{C} \ri \mathcal{C}$ la projection sur le premier
  facteur, et on pose $\mathcal{G}$ la grassmanienne relative
  $Grass (2, p_{*} \mathcal{I} (\fc{n}{2} + 1))$, qui est une
  vari\'et\'e de dimension $2 n + 6$ munie d'une projection canonique
  $q$ sur $\mathcal{C}$. D'apr\`es le lemme pr\'ec\'edent on a un
  morphisme dominant $\mathcal{G} --> \mathcal{S}_{-}$, qui n'est que
  rationnel.

  \begin{lem}
  \label{instab}Un syst\`eme coh\'erent $(\Gamma, \I_{ \{ x \} } (\fc{n}{2} +
  1))$ est $\alpha_{-}-$instable si et seulement si il s'ins�e dans
  une suite exacte
  $$ 0 \ri (\Gamma, \O_{l^{'}}(\fc{n}{2})) \ri (\Gamma,
  \I_{ \{ x \} } (\fc{n}{2} +
  1)) \ri (0, \O_{l^{''}}(\fc{n}{2})) \ri 0 \hspace{1.5cm}(+)$$
  \end{lem}

  \begin{pv} Supposons que le syst\`eme coh\'erent de l'\'enonc\'e soit
  instable, alors il existe un sous-faisceau $\Theta^{'} \subset \I_{ \{ x \} } (\fc{n}{2} +
  1)$ de multiplicit\'e 1 tel que
  $$ \epsilon \hspace{0.1cm} \dim \Gamma^{'} + \chi(\Theta^{'}) > \epsilon +
  \fc{n + 2}{2} \hspace{2cm}(*)$$
  Or comme $0 < \epsilon < 1$ et $\dim \Gamma^{'} \leq 2$,
  on a n\'ecessairement $\chi(\Theta^{'}) \geq \fc{n+2}{2}$. On a d'autre part
  l'inclusion $\Theta^{'} \subset \O_{C}(\fc{n}{2} + 1)$, o\`u $C$
  est la conique support; comme le faisceau $\O_{C}$ est stable, on
  en d\'eduit que $\chi(\Theta^{'}) = \fc{n}{2} + 1$, et donc que
  $\Theta^{'} \simeq \O_{l^{'}}(\fc{n}{2})$ pour une droite $l^{'}$.
  D\`es lors $(*)$ impose $\dim \Gamma^{'} = 2$, et la condition
  n\'ecessaire du lemme s'en d\'eduit. La condition suffisante est
  imm\'ediate. \end{pv} \\

   On note dans cette partie $\Sigma \subset \mathcal{G}$ le ferm\'e
   (muni de sa structure de sch\'ema {\bf{r\'eduit}})
   des points de $\mathcal{G}$ qui repr\'esentent des syst\`emes
   coh\'erents $\epsilon-$instables. Soit $D$ la vari\'et\'e
   d'incidence introduite juste apr\`es la d\'efinition
   \ref{defSigma}, et posons $G := Grass (2, p_*
   \O_{D}(\fc{n}{2}))$.  \\

  \fbox{Notations:} on pose de plus
    \begin{itemize}
    \item $\mathbf{D}_{i}$, pour $i \in \{ 1,2 \}$, la
 sous-vari\'et\'e de $G \times D \times \P_2$ image r\'eciproque de
 la vari\'et\'e d'incidence $D_i$ (voir notations apr\`es la
 d\'efinition \ref{defSigma}) par la projection naturelle
 $G \times D \times \P_2 \ri \P_{2}^* \times \P_{2}^* \times \P_2$,

  \item $\mathbf{\Gamma}$ le fibr\'e tautologique de rang 2 sur $G$;
  $\tau_{i}$ pour $i=1,2$ les images r\'eciproque des fibr\'es
  $\O_{\P_{2}^*}(1)$ sur chacun des deux facteurs $\P_{2}^* \times
  \P_{2}^*$ par la projection $G \times D \ri \P_{2}^* \times
  \P_{2}^*$ et $\zeta$ l'image r\'eciproque sur $G \times D$ du
  fibr\'e sur $D$ provenant de $\O_{\P_2}(1)$.
\end{itemize}
 \hspace{0.1cm} \\
 On a maintenant la proposition suivante, dont la preuve suit de
 pr\`es le sch\'ema de la proposition \ref{Sigma}, bien qu'\'etant un peu
 plus simple du fait que $\mathcal{G}$ est une vari\'et\'e de
 "param\`etres" et non une vari\'et\'e de modules.

\begin{prop}
 On a un isomorphisme $\Sigma \simeq G \times D$.
\end{prop}

 \begin{pv} On construit tout d'abord un morphisme modulaire
 $f: G \times D \ri \mathcal{G}$. L'existence de $f$ est assur\'ee
 par le lemme suivant.

 \begin{lem}
 \label{filtratHNenfamille} Il existe une suite exacte de syst\`emes
 coh\'erents sur $G \times D$
 $$ 0 \ri (\mathbf{\Gamma} \boxtimes \tau_{2}^{-1},
 \O_{\mathbf{D}_{1}}(\fc{n}{2}) \boxtimes \tau_{2}^{-1})
 \ri \mathbf{\Lambda} \ri
 (0, \tau_{2} \boxtimes \zeta^{-1} \otimes \O_{\mathbf{D}_{2}}(\fc{n}{2})) \ri 0 $$
  La famille $\mathbf{\Lambda}$ param\`etre tous les syst\`emes
  coh\'erents du lemme \ref{instab}.
 \end{lem}

 \begin{pv}
 Il suffit de montrer que l'on a l'extension de faisceaux
 correspondante sur $G \times D \times \P_{2}$. On a d\'ej\`a
 l'extension
 $$ 0 \ri \O_{\mathbf{D}_1}(\fc{n}{2}) \boxtimes \tau_{2}^{-1}
 \ri \mathbf{\mathcal{U}} \ri \O_{\mathbf{D}_2}(\fc{n}{2} + 1)
 \ri 0 $$
 o\`u le faisceau $\mathbf{\mathcal{U}}$ param\`etre en chaque point de $G \times
 D$ le faisceau structural, tordu par $\O(\fc{n}{2} + 1)$, de la conique
 image par la projection
 $G \times D \ri \P_{2}^* \times \P_{2}^* \ri \P_5$, la deuxi\`eme
 fl\`eche associant \`a une paire de droites la conique singuli\`ere
 correspondante. Tenant compte \`a pr\'esent du fait qu'un point de
 $D$ param\`etre une droite et un point sur cette droite, on \'ecrit
 la suite exacte sur $G \times D \times \P_{2}$
 $$ 0 \ri \tau_{2} \otimes \zeta^{-1} \otimes \O_{\mathbf{D}_2}(-1)
 \ri \O_{\mathbf{D}_2} \ri \O_{\Delta} \ri 0 $$
 o\`u $\Delta$ est l'image r\'eciproque de la diagonale de $\P_{2}
 \times \P_2$ par la projection $G \times D \times \P_2 \ri D \times
 \P_2 \ri \P_2 \times \P_2$. On v\'erifie en effet facilement que
 l'on a l'isomorphisme
 $$ \underline{\Ext}_{p}^{1} (\O_{\Delta}, \O_{\mathbf{D}_2}(-1))
 \simeq \tau_{2}^{-1} \otimes \zeta $$
 en consid\'erant par exemple la r\'esolution bien connue de la diagonale
 $\Delta$. On obtient donc une suite exacte de faisceaux coh\'erents
 sur $G \times D \times \P_2$
 $$ 0 \ri \O_{\mathbf{D}_1}(\fc{n}{2}) \boxtimes \tau_{2}^{-1}
 \ri \mathbf{\Theta} \ri \tau_2 \boxtimes \zeta^{-1} \otimes
 \O_{\mathbf{D}_2}(\fc{n}{2}) \ri 0 $$
 et une autre suite exacte
 $$ 0 \ri \mathbf{\Theta} \ri \mathcal{U} \ri \O_{\Delta}(\fc{n}{2} + 1)
 \ri 0 $$
 Le faisceau $\mathbf{\Theta}$ param\`etre donc des faisceaux de la
 forme $\I_{ \{ x \} } (\fc{n}{2} + 1)$, o\`u $x$ est un point
 choisi sur la droite $l^{''}$. Maintenant l'extension de faisceaux
 donne naissance \`a une extension de syst\`emes coh\'erents de fa\c{c}on
 naturelle. Il est clair que tous les syst\`emes coh\'erents du lemme
 \ref{instab} sont d\'ecrits de cette fa\c{c}on. Ceci prouve le lemme.
 \end{pv} \\

 Il est clair que $f$ est injectif. Soit $x$ un point de $G \times D$; nous
 notons \'egalement $x$ son image par $f$. Supposons que $x$ corresponde \`a un
 syst\`eme coh\'erent $(\Gamma, \Theta)$ de filtration
 $(\Gamma , \Theta^{'})$, $(0 , \Theta^{''})$. Soit $p: G \times D
 \ri \P_{2}^* \times D$ la projection naturelle. On a un diagramme
 commutatif dont les colonnes sont exactes
$$ \xymatrix{
   \Hom \left( \Gamma , \H^0 (\Theta^{'}) / \Gamma \right)
   \ar@{^{(}->}[r]
   \ar[d] &
   \Hom \left( \Gamma , \H^0 (\Theta) / \Gamma \right) \ar[d] \\
   T_x \left( G \times D \right) \ar[r]^{d_x f} \ar[d]^{d_x p} & T_x
   \left( \mathcal{G} \right) \ar[d]^{d_x q} \\
   T_{p(x)} \left( \P_{2}^* \times D \right) \ar@{^{(}->}[r] & T_{q(x)}
   \mathcal{C} \\
   } $$
 Ceci entraine que $d_x f$ est injectif et par suite $f$ est un plongement.
 \end{pv} \\

 Le r\'esultat suivant d\'etermine le fibr\'e normal
 $N_{\Sigma / \mathcal{G}}$ associ\'e au plongement $\Sigma
 \subset \mathcal{G}$. Nous utilisons une version du Th\'eor\`eme
 3.2 de [L1], lui-m\^eme provenant de [DL]. Les groupes $\Ext$ filtr\'es
 dont il est question dans le th\'eor\`eme proviennent de la
 filtration de Harder-Narasimhan relative \`a un param\`etre
 $\alpha$.

 \begin{theo}
 Soit $S$ une vari\'et\'e alg\'ebrique lisse, et $\mathbf{\Lambda} =
 (\Lambda_{s})_{s \in S}$ une famille plate de syst\`emes
 coh\'erents purs. On suppose que les conditions suivantes sont
 satisfaites $\forall s \in S$:

 \begin{enumerate}
 \item le morphisme de d\'eformation
 infinit\'esimale de Kodaira et Spencer
 $$ \omega_s : \T_{s} S \ri \Ext^1 (\Lambda_s , \Lambda_s) $$
 est surjectif (la famille $\mathbf{\Lambda}$ est dite
 {\bf{compl\`ete}}); \\

 \item $\Ext_{-}^{2}(\Lambda_s , \Lambda_s ) = 0$ et
 $\Ext_{+}^{2}(\Lambda_s , \Lambda_s ) = 0$.
\end{enumerate}

 Alors la strate de Harder-Narasimhan est lisse en $s$, d'espace
 normal $\Ext_{+}^{1}(\Lambda_s , \Lambda_s)$.
 \end{theo}

 Appliquant le th\'eor\`eme pr\'ec\'edent \`a $S = \mathcal{G}$,
 la strate de Harder-Narasimhan est alors $\Sigma$; la filtration
 de la famille $\mathbf{\Lambda}$ au dessus de
 $\Sigma$, d\'etermin\'ee par la suite exacte du lemme
 \ref{filtratHNenfamille}, param\`etre
 en $s \in \Sigma$ la filtration de Harder-Narasimhan du syst\`eme
 coh\'erent
 $\Lambda_s$. Une l\'eg\`ere adaptation du th\'eor\`eme
 ci-dessus donne alors le fibr\'e normal \`a $\Sigma$, moyennant la
 v\'erification de la compl\'etude de $\mathbf{\Lambda}$, voir
 condition (1).

\begin{lem}
 Soit $p \in \mathcal{G}$, et $\Lambda$ le syst\`eme coh\'erent
 fibre de la famille $\mathbf{\Lambda}_{\mathcal{G}}$ en $p$.
 Le morphisme
 de d\'eformation de Kodaira-Spencer
 $$ \omega: \T_{p} \mathcal{G} \ri \Ext^1 (\Lambda, \Lambda) $$
 est surjectif.
\end{lem}

 \begin{pv} Posons $\Lambda = (\Gamma, \Theta)$ avec
 $\Theta = \I_{ \{ x \} }
 (\fc{n}{2} + 1)$. On a la suite exacte
 $$ 0 \ri \Hom (\Gamma, \H^0 (\Theta) / \Gamma)
 \ri \Ext^1 (\Lambda, \Lambda) \ri \Ext^1 (\I_{ \{ x \} },
 \I_{ \{ x \} }) \ri 0 $$
 On a d'autre part la suite exacte d'espaces tangents associ\'ee \`a
 la submersion $\rho: \mathcal{G} \ri \mathcal{C}$
 $$ 0 \ri \Hom (\Gamma, \H^0 (\Theta) / \Gamma)
 \ri \T_{p} \mathcal{G} \ri \T_{\rho(p)} \mathcal{C} \ri 0 $$
 Le morphisme $\omega$ induit donc une application lin\'eaire
 $\omega^{'} : \T_{\rho(p)} \mathcal{C} \ri \Ext^1 (\I_{ \{ x \} },
 \I_{ \{ x \} })$. Il reste \`a v\'erifier que $\omega^{'}$ est
 surjectif. Mais le faisceau $\I_{ \{ x \} }$, plac\'e en
 degr\'e  0, est quasi-isomorphe au complexe $\I^{*} :=
 0 \ri \O_{C} \ri \mathbb{C}_x \ri 0 $, o\`u $C$ est la conique
 support. Une suite spectrale montre alors que $\Ext^1 (\I_{ \{ x \} },
 \I_{ \{ x \} })$ est le groupe de cohomologie en degr\'e 0 du complexe
 $$ \Ext^1 (\mathbb{C}_x , \O_{C}) \stackrel{(1)}{\ri} \Ext^1 (\O_C , \O_C)
 \oplus
 \Ext^1 (\mathbb{C}_x , \mathbb{C}_x ) \stackrel{(2)}{\ri} \Ext^1 (\O_C ,
 \mathbb{C}_x ) $$
 Mais la fl\`eche (2) s'identifie \`a la projection
 $ \T (\P_5 \times \P_2)_{\rho(p)} \ri (N_{\mathcal{C} / \P_5
 \times \P_2})_{\rho(p)}$. On a donc une surjection $\T_{\rho(p)}
 \mathcal{C} \ri Ker (2) / Im (1) \simeq \Ext^1 (\I_{ \{ x \} },
 \I_{ \{ x \} })$. \end{pv} \\

 On a d\`es lors la proposition suivante.

 \begin{prop}
 \label{fibrenormal}Avec les notations du lemme \ref{filtratHNenfamille}, si
 l'on pose $\mathbf{\Lambda}^{'} = (\mathbf{\Gamma},
 \O_{\mathbf{D}_{1}}(\fc{n}{2}))$ et
 $\mathbf{\Lambda}^{''} = (0,
 \O_{\mathbf{D}_{2}}(\fc{n}{2}))$, le fibr\'e normal $N_{\Sigma /
 \mathcal{G}}$ est isomorphe \`a
 $\tau_{2}^{\otimes 2} \otimes \zeta^{-1}
 \otimes \underline{\Ext}_{p}^1 (\mathbf{\Lambda}^{'},
 \mathbf{\Lambda}^{''})$.
 \end{prop}

 En effet la suite exacte du lemme \ref{filtratHNenfamille}
 s'\'ecrit
 $$ 0 \ri \tau_{2}^{-1} \otimes \mathbf{\Lambda}^{'}
 \ri \mathbf{\Lambda} \ri \tau_{2} \otimes \zeta^{-1} \otimes
 \mathbf{\Lambda}^{''} \ri 0 $$
 et sa restriction en chaque point de $\Sigma$ donne la filtration
 d'Harder-Narasimhan du syst\`eme coh\'erent $\Lambda$. On applique
 alors le th\'eor\`eme ci-dessus.

 Soit $\widetilde{\mathcal{G}}$ l'\'eclat\'e de $\mathcal{G}$
 le long du ferm\'e $\Sigma$. D'apr\`es la proposition qui
 pr\'ec\`ede un point du diviseur exceptionnel $\mathcal{E}$
 param\`etre une suite exacte non scind\'ee
 $$ 0 \ri (0, \O_{l^{'}}(\fc{n}{2})) \ri (\Gamma, \Theta)
 \ri (\Gamma, \O_{l^{''}}(\fc{n}{2})) \ri 0 \hspace{1.5cm} (+)$$
 Celle-ci d\'efinit un syst\`eme coh\'erent $\epsilon-$semi-stable.
 Ceci justifie la proposition suivante.

 \begin{prop}
  \label{diagcart}Il existe un morphisme canonique $\widetilde{\mathcal{G}}
  \ri \mathcal{S}_{-}$ rendant commutatif le diagramme suivant
  $$ \xymatrix{
   &  \widetilde{\mathcal{G}} \ar[dl] \ar[dr] & \\
   \mathcal{G} \ar@{-->}[rr] & & \mathcal{S}_{-} \\
   } $$
   o\`u le morphisme en pointill\'e est le morphisme
   rationnel canonique de lieu d'ind\'etermination $\Sigma$. De plus
   on a un diagramme {\bf{cart\'esien}}
   $$ \xymatrix{
   & \mathcal{E} \ar[dl] \ar[rd] & \\
   \Sigma \ar[dr] & & \mathcal{P} \ar[dl] \\
   & G \times \P_{2}^* & \\
   } $$
   o\`u le morphisme $\Sigma \ri G \times \P_{2}^*$ est celui de
   l'oubli du point marqu\'e sur la conique, et $\mathcal{P} \ri G
   \times \P_{2}^*$ est le fibr\'e projectif (au sens usuel)
   $\P (\underline{\Ext}_{p}^1 (\mathbf{\Lambda}^{'},
   \mathbf{\Lambda}^{''}))$ (voir proposition \ref{fibrenormal}); il
   existe un morphisme canonique $\mathcal{P} \ri \mathcal{S}_{-}$
   birationnel sur son image form\'e des classes de syst\`emes
   coh\'erents s'ins\'erant dans une suite exacte de type (+).
 \end{prop}

 \begin{pv} La preuve est identique \`a celle du th\'eor\`eme
 \ref{eclatements}. La cart\'esianit\'e du diagramme est cons\'equence
 de la \end{pv} \\

  On peut imiter sur $\mathcal{G}$ la construction des fibr\'es d\'eterminants
  sur les espaces de modules $\mathcal{S}_{\alpha}$. On reprend les
  notations du d\'ebut de cette section. On note
  $\mathbf{\Gamma}_{\mathcal{G}}$ le fibr\'e tautologique de rang 2
  sur $\mathcal{G}$, et $p,q$ les
  projections
  $$ \xymatrix{
   \mathcal{G} \times \P_2 \ar[d]^{p} \ar[r]^{q} & \P_2 \\
   \mathcal{G} & \\
   } $$
 On pose le fibr\'e inversible $\mathcal{D}_{\mathcal{G}}$
 sur $\mathcal{G}$ d\'efini par
 $$ \mathcal{D}_{\mathcal{G}} =
 det \hspace{0.1cm} p_{\hspace{0.08cm}!} (\mathcal{J} (\fc{n}{2} + 1) \cdot q^* h ) \otimes
 det \hspace{0.1cm} \mathbf{\Gamma}_{\mathcal{G}}^{-1} $$
 C'est l'analogue du fibr\'e $\mathcal{D}_{-}$ sur
 $\mathcal{S}_{-}$, et il est plus simple \`a construire du fait de
 l'existence d'une famille de syst\`emes coh\'erents
 param\'etr\'ee par $\mathcal{G}$.

 \begin{lem}
  \label{expressDG}On a $\mathcal{D}_{\mathcal{G}} \simeq
   det \hspace{0.1cm} \mathbf{\Gamma}_{\mathcal{G}}^{-1}
   \otimes \O_{\mathcal{C}}(\fc{n}{2}, -1)$; utilisant les notations
   du lemme \ref{filtratHNenfamille},
   on a de plus
   $\mathcal{D}_{\mathcal{G}} |_{\Sigma} \simeq det \hspace{0.1cm}
   \mathbf{\Gamma}^{-1} \otimes \tau_{1}^{\otimes \fc{n}{2}}
   \otimes \tau_{2}^{\otimes \fc{n}{2} + 2} \otimes \zeta^{-1}$.
\end{lem}

 La notation $\O_{\mathcal{C}}(.,.)$ fait r\'ef\'erence au produit
 externe des fibr\'es inversibles provenant des projections de
 $\mathcal{C}$ sur chacun des deux facteurs. \\

 \begin{pv} De la suite exacte
 $0 \ri \mathcal{J} \ri \O_{\mathcal{C} \times_{\P_5} \mathcal{C}}
 \ri \O_{\Delta} \ri 0 $ et de la suite exacte
 $$ 0 \ri \O_{\mathcal{C}}(-1, 0) \boxtimes \O_{\P_2}(-2)
 \ri \O_{\mathcal{C} \times \P_2} \ri \O_{\mathcal{C} \times_{\P_5}
 \mathcal{C}} \ri 0 $$
 on d\'eduit que
 $$ \mathcal{D}_{\mathcal{G}} \simeq
   det \hspace{0.1cm} \mathbf{\Gamma}_{\mathcal{G}}^{-1} \otimes
 \O_{\mathcal{C}}(\fc{n}{2}, 0) \otimes
 det \hspace{0.1cm} p_{\hspace{0.08cm}!} (\O_{\Delta}(\fc{n}{2} + 1) \cdot q^* h)^{-1} $$
     Maintenant la r\'esolution du faisceau $\O_{\Delta}$ donne la
     suite exacte
$$ 0 \ri \O_{\mathcal{C}}(0,-2) \boxtimes \O_{\P_2}(\fc{n}{2})
\ri \O_{\mathcal{C}}(0,-1) \boxtimes Q (\fc{n}{2}) \ri \O_{\mathcal{C}
\times \P_2} (\fc{n}{2} + 1) \ri \O_{\Delta}(\fc{n}{2} + 1) \ri 0 $$
 d'o\`u l'on d\'eduit
 $$ det \hspace{0.1cm} p_{\hspace{0.08cm}!} (\O_{\Delta}(\fc{n}{2} + 1) \cdot q^* h)
 \simeq \O_{\mathcal{C}} (0, \chi (Q |_{l} (\fc{n}{2})) -
 2 \chi (\O_{l}(\fc{n}{2}))) \hspace{1.5cm} (*) $$
 On rappelle en effet que $h$ est la classe d'un faisceau structural
 $\O_{l}$, avec $l$ une droite. Comme on a $Q |_{l} \simeq \O_l
 \oplus \O_l (1)$, on peut dire que $(*) \simeq
 \O_{\mathcal{C}}(0,1)$. Nous savons que la construction
 d'un tel fibr\'e $\mathcal{D}_{\mathcal{G}}$ est fonctorielle (voir
 [L3]), et utilisant l'isomorphisme $G \times D \stackrel{\simeq}{\ri}
 \Sigma$ on a $\mathbf{\Gamma}_{\mathcal{G}} |_{\Sigma}
 \simeq \mathbf{\Gamma} \boxtimes \tau_{2}^{-1}$. Le lemme s'ensuit. \end{pv} \\

 Notons $\varphi: \widetilde{\mathcal{G}} \ri \mathcal{G}$ et
 $\psi: \widetilde{\mathcal{G}} \ri \mathcal{S}_{-}$ les
 morphismes pr\'ec\'edents. On a une projection compos\'ee
 $\widetilde{\mathcal{G}} \stackrel{\varphi}{\ri} \mathcal{G} \ri
 \mathcal{C} \ri \P_2$, et notons encore $\zeta$ l'image r\'eciproque
 sur $\mathcal{G}$ et $\widetilde{\mathcal{G}}$ de la classe de
 $\O_{\P_2}(1)$. Comme une droite coupe une conique g\'en\'erale en
 deux points et que la fibre g\'en\'erale de la projection
 $\widetilde{\mathcal{G}} \stackrel{\psi}{\ri} \mathcal{S}_{-}$
 est la conique support en un point, on a d'apr\`es la formule
 de projection
 $$ c_{1} (\psi^* \mathcal{D}_{-})^{2 n + 5} \cdot \zeta \cap
 [\widetilde{\mathcal{G}}] =
 2 \hspace{0.1cm} c_{1} (\mathcal{D}_{-})^{2 n + 5}
 \cap [\mathcal{S}_{-}] $$
 D'autre part on a \'evidemment
 $$ c_{1}(\varphi^* \mathcal{D}_{\mathcal{G}})^{2 n + 5} \cdot
 \zeta \cap [\widetilde{\mathcal{G}}] =
 c_{1} (\mathcal{D}_{\mathcal{G}})^{2 n + 5} \cdot \zeta
 \cap [\mathcal{G}] $$
 On pose donc naturellement
 $$ \Delta = 2 \hspace{0.1cm} c_{1} (\mathcal{D}_{-})^{2 n + 5}
               \cap [\mathcal{S}_{-}] -
               c_{1} (\mathcal{D}_{\mathcal{G}})^{2 n + 5} \cdot \zeta
               \cap [\mathcal{G}] $$
 qui est \'egal d'apr\`es ce qui pr\'ec\`ede \`a
 $(c_{1} (\psi^* \mathcal{D}_{-})^{2 n + 5}  -
 c_{1} (\varphi^* \mathcal{D}_{\mathcal{G}})^{2 n + 5}) \cdot \zeta \cap
 [\widetilde{\mathcal{G}}]$.

 \begin{lem}
 Notant $e$ la classe du diviseur exceptionnel
 $\mathcal{E}$, on a la relation
 $$ c_1 (\varphi^* \mathcal{D}_{\mathcal{G}}) =
 c_1 (\psi^* \mathcal{D}_{-}) + e $$
 \end{lem}

 \begin{pv} La d\'emonstration de ce lemme est identique \`a celle
 du lemme \ref{sautfibdet}. Elle s'appuie sur l'isomorphisme
 $\varphi^* \mathcal{D}_{\mathcal{G}} \simeq
  \psi^* \mathcal{D}_{-} \otimes \O(\mathcal{E})$, qui se d\'emontre \`a l'aide
  de la suite exacte dans la preuve de la proposition
  \ref{diagcart}. \end{pv} \\

  \fbox{{\bf{Calculons $\Delta$}}.}
  On pose $d_{\mathcal{G}} = c_1 (\mathcal{D}_{\mathcal{G}})$ et
  $d_{-} = c_1 (\mathcal{D}_{-})$. Par un raisonnement identique
  \`a celui effectu\'e \`a la section \ref{calculsauts}, on peut
  \'ecrire
  $$ \Delta = - \fc{1}{e} \cdot (\psi^* (d_{-})^{2 n + 5} |_{\mathcal{E}}
  - \varphi^* (d_{\mathcal{G}})^{2 n + 5} |_{\mathcal{E}}) \cdot \zeta
  \cap [\mathcal{E}] $$
 On pose $\gamma_i$, pour $i=1,2$ les classes de Chern du fibr\'e
 tautologique $\mathbf{\Gamma}$ sur $G = Grass (2, p_* \O_D
 (\fc{n}{2}))$. Rappelons alors que $\mathcal{E}$ s'ins\`ere dans un
 diagramme cart\'esien (voir proposition \ref{diagcart}) faisant
 intervenir $\Sigma$, $\mathcal{P}$ et $G$; soit $l_{\mathcal{P}}$
 le fibr\'e relativement ample sur le fibr\'e projectif
 $\mathcal{P}$.

 \begin{lem}
 On a $\psi^* (d_{-}) |_{\mathcal{E}} = l_{\mathcal{P}}
 + \fc{n}{2} (\tau_{1} + \tau_2 ) - \gamma_1$.
 \end{lem}

 On a par ailleurs $d_{\mathcal{G}} |_{\Sigma} = - \gamma_1 +
 \fc{n}{2} \tau_1 + (\fc{n}{2} + 2) \tau_2 - \zeta$ (voir
 le deuxi\`eme point
 du lemme \ref{expressDG}). Posons
 $v = - \gamma_1 + \fc{n}{2} (\tau_1 + \tau_2)$, on en tire
 $$ \Delta = \sum_{0 \leq j \leq i \leq 2 n + 4} \binom{2 n + 5}{i+1}
 \hspace{0.1cm} (2 \tau_2 - \zeta)^{i-j} \cdot v^{2 n + 4 - i} \cdot
 l_{\mathcal{P}}^j \cdot \zeta \cap [\mathcal{E}] \hspace{2cm}(*)$$
  Notons $\rho: \mathcal{P} \ri G \times \P_{2}^*$ la projection
  canonique, on a donc $\forall j$ \hspace{0.2cm}$\rho_* (l_{\mathcal{P}}^j)
  = s_{j - n - 2} (\mathcal{F}^{\vee})$, o\`u $\mathcal{F}$
  est le fibr\'e $\underline{\Ext}_{p}^1 (\mathbf{\Lambda}^{'},
  \mathbf{\Lambda}^{''})$ de la proposition \ref{fibrenormal}. En
  effet le fibr\'e projectif $\mathcal{P}$ est isomorphe au fibr\'e
  au sens de Grothendieck $\mathbb{P}(\mathcal{F}^{\vee})$. Comme
  pr\'ec\'edemment dans une situation identique, on a deux r\'esolutions
  $$ \xymatrix{
   0 \ar[r] & \mathbf{\Gamma}^* \boxtimes \S^{\fc{n}{2}} Q \ar[r] &
   \mathcal{F} \ar[r] & \mathcal{J}(2,1) \ar[r] \ar@{=}[d] & 0 \\
   0 \ar[r] & \O \ar[r] & \tau_1 \boxtimes Q \ar[r] &  \mathcal{J}(2,1)
   \ar[r] & 0 \\
   } $$
 On a d\`es lors
 $$ \begin{array}{llcc}
  s(\mathcal{F}^{\vee}) = c(\mathcal{F})^{-1} & =
 s(\mathbf{\Gamma} \boxtimes \S^{\fc{n}{2}} Q^{*}) \cdot c(\tau_1
 \boxtimes Q)^{-1} \\
  & = s(\mathbf{\Gamma} \boxtimes \S^{\fc{n}{2}} Q^{*}) \cdot
  (1 + \tau_1 - \tau_2) \cdot (1 + \tau_1)^{-3} \\
  \end{array} $$
 On r\'eindexe la somme $(*)$ en posant $i:=i+1$, $j:=j+1$
 $$ \Delta = \sum_{n + 3 \leq j \leq i \leq 2 n + 5} \binom{2 n +
 5}{i} \hspace{0.1cm} v^{2 n + 5 - i} \cdot (2 \tau_2 - \zeta)^{i-j}
 \cdot s_{j-n-3}(\mathcal{F}^{\vee}) \cdot \zeta \cap \left[ G \times
 D \right] $$
 Ceci revient \`a \'ecrire
 $$ \begin{array}{lllcc}
 \Delta & = \left[ (1 + v)^{2 n + 5} \cdot \fc{\zeta \cdot (1 +
 \tau_1 - \tau_2)}{\left[ 1 - ( 2 \tau_2 - \zeta) \right] \cdot (1 +
 \tau_1 )^3} \cdot s \left( \mathbf{\Gamma} \boxtimes
 \S^{\fc{n}{2}} Q^* \right) \right]_{n + 3} \cap \left[ G \times D
 \right] \\
 & = \sum_{ \begin{array}{ll}
            0 \leq m \hspace{0.1cm}; \hspace{0.1cm} l \\
            n-2 \leq l+m \leq n + 2 \\
            \end{array} } (-1)^m \hspace{0.1cm} \binom{2 n + 5}{m}
            \hspace{0.1cm} \fc{ \left[ 1 + \fc{n}{2} \left( \tau_1
            + \tau_2 \right) \right]^{2 n + 5 - m} \cdot \left( 1 + \tau_1
            - \tau_2 \right) \cdot \zeta }{ \left[ 1 - ( 2 \tau_2 - \zeta )
            \right] \cdot ( 1 + \tau_1 )^3 } \cdot \gamma_{1}^m \cdot
            s_{l} \left( \mathbf{\Gamma} \boxtimes \S^{\fc{n}{2}} Q^* \right)
            \\
            \end{array} $$
 On pose $\mu = \binom{\fc{n}{2} + 1}{2}$. On prend les images directes des classes $\gamma_{1}^m \cdot
            s_{l} \left( \mathbf{\Gamma} \boxtimes \S^{\fc{n}{2}} Q^*
            \right)$ par le morphisme de projection $G \times
            \P_{2}^* \ri \P_{2}^* \times \P_{2}^*$ qui d'apr\`es le lemme
            \ref{secondeimagedirecte} sont \'egales
            \`a
$$ \fc{1}{\mu +1} \hspace{0.1cm} \sum_{ \begin{array}{ll}
                                          0 \leq l_0 \leq l_1 \leq l_2 \leq l_3 \leq l \\
                                          l + m - n \leq l_0 + (l_3 - l_2) \leq \mathrm{min}(l + m + 2 - n,2) \\
                                          \end{array} }
                                          f_{0,n}(\underline{\mathbf{l}},m) \hspace{0.1cm}
 \tau_{1}^{l + m + 2 - n - l_0 - (l_3 - l_2)} \cdot \tau_{2}^{l_0 +
 (l_3 - l_2)} $$
 Dans le lemme \ref{secondeimagedirecte} la notation $\zeta$ d\'esignait
 le nombre $\binom{\fc{n}{2} - \alpha + 1}{2}$; pour $\alpha = 0$ on remplace
 ici cette \'ecriture par $\mu$. Apr\`es avoir pos\'e la r\'eindexation $m:=l+m - (n-2)$, $l:=l$, on
 obtient l'expression suivante de $\Delta$
$$ \begin{array}{llcc}
\fc{(-1)^n}{\mu + 1} \hspace{0.1cm} \sum_{
\begin{array}{lll}
0 \leq m \leq 4 \\
0 \leq l_0 \leq l_1 \leq l_2 \leq l_3 \leq l \leq m + (n-2) \\
m-2 \leq l_0 + (l_3 - l_2) \leq \mathrm{min}(m,2) \\
 \end{array} } (-1)^{m+l} \hspace{0.1cm} \binom{2 n + 5}{n + 7 + l-m}
 \hspace{0.1cm} f_{0,n} (\underline{\mathbf{l}}, m + (n-2) - l) \\
 \hspace{8cm} \fc{ \left[ 1 + \fc{n}{2} \left( \tau_1
            + \tau_2 \right) \right]^{2 n + 5 - m} \cdot \left( 1 + \tau_1
            - \tau_2 \right) \cdot \zeta }{ \left[ 1 - ( 2 \tau_2 - \zeta )
            \right] \cdot ( 1 + \tau_1 )^3 } \cdot \tau_{1}^{m - l_0
            - (l_3 - l_2)} \cdot \tau_{2}^{l_0 + (l_3 - l_2)} \\
            \end{array} $$
 On \'ecrit alors
 $$ \fc{\zeta}{1 - 2 \tau_2 + \zeta} = \fc{1}{1 - 2 \tau_2} \cdot
 \left[ \fc{\zeta}{1 + \fc{\zeta}{1 - 2 \tau_2}} \right] =
 \fc{\zeta}{1 - 2 \tau_2} \cdot \left( 1 - \fc{\zeta}{1 - 2 \tau_2}
 \right) $$
 Tenant compte enfin du fait que seuls les nombres d'intersection
 sur $D$
 $$ \tau_{2}^2 \cdot \zeta \cap \left[ D \right] = \tau_2 \cdot
 \zeta^2 \cap \left[ D \right] = 1 $$
 sont non nuls, on obtient en posant les fonctions rationnelles
 $H_i (\tau_1 , \tau_2 )$ suivantes pour $i = 1 , 2$
 $$ H_i (\tau_1 , \tau_2 ) = \fc{ \left[ 1 + \fc{n}{2} \left( \tau_1
            + \tau_2 \right) \right]^{n + 7 + l - m} \cdot \left( 1 + \tau_1
            - \tau_2 \right) }{ ( 1 + \tau_1 )^3 \cdot
            (1 - 2 \tau_2)^i } $$
 l'expression finale de $\Delta$
 $$ \fbox{$\displaystyle \begin{array}{llcc}
 \Delta & = \fc{(-1)^n}{\mu + 1} \hspace{0.1cm} \sum_{
\begin{array}{llll}
0 \leq m \leq 4 \\
0 \leq l_0 \leq l_1 \leq l_2 \leq l_3 \leq l \leq m + (n-2) \\
m-2 \leq l_0 + (l_3 - l_2) \leq \mathrm{min}(m,2) \\
1 \leq i \leq 2
 \end{array} } (-1)^{m+l} \hspace{0.1cm} \binom{2 n + 5}{n + 7 + l-m}
 \hspace{0.1cm} f_{0,n} (\underline{\mathbf{l}}, m + (n-2) - l) \\
 & \hspace{8cm} \left[ H_i (\tau_1 , \tau_2 ) \right]^{2 + l_0 + (l_3
 - l_2 ) - m}_{i - l_0 - (l_3 - l_2)} \\
 \end{array} $} $$
 la notation $\left[ H_i (\tau_1 , \tau_2) \right]^{x}_{y}$
 d\'esignant le coefficient devant le mon\^ome $\tau_{1}^x \cdot
 \tau_{2}^y$ dans le d\'eveloppement en s\'erie formelle de la
 fraction rationnelle $H_i$. Naturellement si $i=1$ et $l_0 + (l_3 - l_2) = 2$
 ce coefficient est nul. \\

 Il reste \`a calculer $\I_{\mathcal{G}} = c_{1}
(\mathcal{D}_{\mathcal{G}})^{2 n + 5} \cdot \zeta \cap
[\mathcal{G}]$. Posons maintenant $q: \mathcal{G} \ri \mathcal{C}$
la projection canonique, et $\gamma = c_1 (\Gamma_{\mathcal{G}})$,
    $\tau = c_1 (\O_{\mathcal{C}}(1,0))$, ainsi que $\zeta = c_1 (\O_{\mathcal{C}}(0,1))$.

 D'apr\`es le lemme \ref{expressDG} on a $ c_1 (\mathcal{D}_{\mathcal{G}}) =
 - \gamma + \fc{n}{2} \tau - \zeta$, on a donc
 $$ \begin{array}{lllcc}
 \I_{\mathcal{G}} & = - \left( \gamma - \fc{n}{2} \tau + \zeta
 \right)^{2 n + 5} \cdot \zeta \cap \left[ \mathcal{G} \right] \\
 & = - \left[ \left( \gamma - \fc{n}{2} \tau \right)^{2 n + 5} \cdot
 \zeta \hspace{0.1cm} + \hspace{0.1cm} (2 n + 5)
 \hspace{0.1cm} \left( \gamma - \fc{n}{2} \tau \right)^{2 n + 4}
 \cdot \zeta^2 \right] \cap \left[ \mathcal{G} \right] \\
 \end{array} $$
 Calculons d'abord $\I_{\mathcal{G}}^{'} = \left( \gamma - \fc{n}{2} \tau \right)^{2 n + 5} \cdot
 \zeta \cap \left[ \mathcal{G} \right]$. Dans ce qui suit, on note $\mathcal{F}$
 le fibr\'e vectoriel $p_{*} \mathcal{J} \left( \fc{n}{2} + 1 \right) $
 o\`u $p: \mathcal{C} \times_{\mathbf{P}_5} \mathcal{C} \ri \mathcal{C}$
 est la premi\`ere projection et $\mathcal{J}$ est le faisceau d'id\'eaux
 du plongement diagonal
 $\mathcal{C} \simeq \Delta \subseteq \mathcal{C} \times_{\mathbf{P}_5} \mathcal{C}$.
 Ci-dessous on prend
 \`a la deuxi\`eme ligne l'image directe des classes $\gamma^m$ par
 la morphisme $q$
 $$ \begin{array}{llllcc}
 \I_{\mathcal{G}}^{'} & = \sum_{2 n \leq m \leq 2 n + 5}
 \binom{2 n + 5}{m} \hspace{0.1cm} \left( - \fc{n}{2} \right)^{2 n +
 5 - m} \hspace{0.1cm} \gamma^m \cdot \tau^{2 n + 5 - m} \cdot \zeta
 \cap \left[ \mathcal{G} \right] \\
  & = \sum_{ \begin{array}{ll}
             2 n \leq m \leq 2 n + 5 \\
             n \leq p \leq m + 1 - n \\
             \end{array} }
             \binom{2 n + 5}{m} \hspace{0.1cm} \left( - \fc{n}{2}
             \right)^{2 n + 5 - m} \hspace{0.1cm} \binom{m}{p} \hspace{0.1cm}
             \left| \begin{array}{llcc}
                    s_{m-p-n}(\mathcal{F}) & s_{p-n-1}(\mathcal{F})
                    \\
                    s_{m-p-n+1}(\mathcal{F}) & s_{p-n}(\mathcal{F})
                    \\
                    \end{array} \right| \cdot \tau^{2 n + 5 - m}
                    \cdot \zeta \cap \left[ \mathcal{C} \right] \\
                    \end{array} $$
 Apr\`es avoir pos\'e la r\'eindexation $m := m - 2 n $, $p := p -
 n$ dans la somme pr\'ec\'edente, on obtient
 $$ \I_{\mathcal{G}}^{'} = \sum_{ \begin{array}{ll}
  0 \leq m \leq 5 \\
  0 \leq p \leq m + 1 \\
  \end{array} } \hspace{0.1cm} \binom{2 n + 5}{5-m} \hspace{0.1cm}
  \left( - \fc{n}{2} \right)^{5 - m} \hspace{0.1cm} \binom{m + 2
  n}{p+n} \hspace{0.1cm} \left| \begin{array}{llcc}
                    s_{m-p}(\mathcal{F}) & s_{p-1}(\mathcal{F})
                    \\
                    s_{m-p +1}(\mathcal{F}) & s_{p}(\mathcal{F})
                    \\
                    \end{array} \right| \cdot \tau^{5 - m}
                    \cdot \zeta \cap \left[ \mathcal{C} \right] $$
 Maintenant en prenant l'image directe par $p$ de la suite exacte de
 faisceaux coh\'erents sur $\mathcal{C} \times_{\mathbf{P}_5}
 \mathcal{C}$ donn\'ee par
 $$ 0 \ri \mathcal{J} \left( \fc{n}{2} + 1 \right) \ri \O_{\mathcal{C} \times_{\mathbf{P}_5}
 \mathcal{C}} \left( \fc{n}{2} + 1 \right) \ri \O_{\Delta} \left(
 \fc{n}{2} + 1 \right) \ri 0 $$
 on obtient la suite exacte suivante de fibr\'es vectoriels sur $\mathcal{C}$
$$ 0 \ri \mathcal{F} \ri p_* \O_{\mathcal{C} \times_{\mathbf{P}_5}
\mathcal{C}} \left( \fc{n}{2} + 1 \right) \ri \O_{\mathcal{C}}
\left( \fc{n}{2} + 1 \right) \ri 0 $$
 En effet $p$ induit un isomorphisme $p |_{\Delta} : \Delta \simeq
 \mathcal{C}$. On a de plus la suite exacte de faisceaux
 coh\'erents sur $\mathcal{C} \times \P_2$
 $$ 0 \ri \O_{\mathcal{C}}(-1) \boxtimes \O_{\mathbf{P}_2} \left(
 \fc{n}{2} -1 \right) \ri \O_{\mathcal{C}} \boxtimes
 \O_{\mathbf{P}_2} \left( \fc{n}{2} + 1 \right) \ri \O_{\mathcal{C} \times_{\mathbf{P}_5}
\mathcal{C}} \left( \fc{n}{2} + 1 \right) \ri 0 $$
 En prenant l'image directe $p_*$, on obtient encore la suite exacte
 $$ 0 \ri \O_{\mathcal{C}} \left( -1 \right)^{\binom{\fc{n}{2}
 +1}{2}} \ri \O_{\mathcal{C}}^{\binom{\fc{n}{2} + 3}{2}} \ri p_* \O_{\mathcal{C} \times_{\mathbf{P}_5}
\mathcal{C}} \left( \fc{n}{2} + 1 \right) \ri 0 $$
 On tire de ce qui pr\'ec\`ede que $s(\mathcal{F}) = \left( 1 + \tau
 \right)^{\mu} \cdot \left( 1 - \left( \fc{n}{2}
 + 1 \right) \zeta \right) $. On a alors
 $$ \begin{array}{llllcc}
    s_{m-p} (\mathcal{F}) \cdot s_{p}(\mathcal{F}) - s_{m-p+1}
    (\mathcal{F}) \cdot s_{p-1} (\mathcal{F})  & = \left[
    \binom{\mu}{m-p} \binom{\mu}{p} - \binom{\mu}{m-p+1}
    \binom{\mu}{p-1} \right] \tau^m \\
    & + \left( \fc{n}{2} + 1 \right) \left[ \binom{\mu}{m-p+1}
    \binom{\mu}{p-2} - \binom{\mu}{m-p-1} \binom{\mu}{p} \right]
    \tau^{m-1} \cdot \zeta \\
    & = \fc{m - 2 p + 1}{\mu + 1} \hspace{0.1cm} \binom{\mu +1}{m-p+1}
    \binom{\mu +1}{p} \hspace{0.1cm} \tau^m \\
    & + \left( \fc{n}{2} + 1 \right) \left[ \binom{\mu}{m-p+1}
    \binom{\mu}{p-2} - \binom{\mu}{m-p-1} \binom{\mu}{p} \right]
    \tau^{m-1} \cdot \zeta \\
    \end{array} $$
    Ce qui donne finalement, compte tenu du fait que $\tau^5 \cdot
    \zeta \cap \left[ \mathcal{C} \right] = 2$ et $\tau^4 \cdot
    \zeta^2 \cap \left[ \mathcal{C} \right]$. Cela donne le
    r\'esultat final pour $\I_{\mathcal{G}}^{'}$
    $$ \fbox{$\displaystyle \begin{array}{lllcc}
    \I_{\mathcal{G}}^{'} = & \sum_{ \begin{array}{ll}
                                    0 \leq m \leq 5 \\
                                    0 \leq p \leq m+1 \\
                                    \end{array} }
                                    \hspace{0.1cm} \binom{2 n
                                    +5}{5-m} \hspace{0.1cm} \left( -
                                    \fc{n}{2} \right)^{5-m}  \hspace{0.1cm}
                                    \binom{m + 2 n}{p+n} \\
   & \hspace{3cm} \left[ 2 \hspace{0.1cm} \fc{m-2p+1}{\mu + 1} \hspace{0.1cm}
   \binom{\mu +1}{m-p+1} \hspace{0.1cm}\binom{\mu + 1}{p}
   + \left( \fc{n}{2} + 1 \right) \hspace{0.1cm}
   \left[ \binom{\mu}{m-p+1} \hspace{0.1cm} \binom{\mu}{p-2}
   \hspace{0.1cm} - \hspace{0.1cm} \binom{\mu}{m-p-1}
   \hspace{0.1cm}\binom{\mu}{p} \right] \right] \\
   \end{array} $} $$

 Posons maintenant et calculons $\I_{\mathcal{G}}^{''} = \left(
 \gamma - \fc{n}{2} \tau \right)^{2 n + 4} \cdot \zeta^2 \cap \left[
 \mathcal{G} \right]$. Des calculs identiques aux pr\'ec\'edents
 donnent
 $$ \I_{\mathcal{G}}^{''} =  \sum_{ \begin{array}{ll}
                                   0 \leq m \leq 4 \\
                                   0 \leq p \leq m+1 \\
                                   \end{array} } \hspace{0.1cm}
                                   \binom{2 n + 4}{4-m} \hspace{0.1cm}
                                   \left( - \fc{n}{2} \right)^{4-m}
                                   \hspace{0.1cm} \binom{m+2n}{p+n}
                                   \hspace{0.1cm}\left| \begin{array}{llcc}
                    s_{m-p}(\mathcal{F}) & s_{p-1}(\mathcal{F})
                    \\
                    s_{m-p +1}(\mathcal{F}) & s_{p}(\mathcal{F})
                    \\
                    \end{array} \right| \cdot \tau^{4 - m}
                    \cdot \zeta^2 \cap \left[ \mathcal{C} \right] $$
 ce qui donne finalement
 $$ \fbox{$\displaystyle \I_{\mathcal{G}}^{''} = \sum_{ \begin{array}{ll}
                                   0 \leq m \leq 4 \\
                                   0 \leq p \leq m+1 \\
                                   \end{array} } \hspace{0.1cm}
                                   \binom{2 n + 4}{4-m} \hspace{0.1cm}
                                   \left( - \fc{n}{2} \right)^{4-m}
                                   \hspace{0.1cm} \binom{m+2n}{p+n}
                                   \hspace{0.1cm}\fc{m - 2 p +
                                   1}{\mu+1} \hspace{0.1cm}
                                   \binom{\mu +1}{m-p+1}
                                   \hspace{0.1cm} \binom{\mu
                                   +1}{p} $} $$
 On a donc obtenu la proposition suivante.

 \begin{prop}
 On a avec les expressions de $\I_{\mathcal{G}}^{'}$ et $\I_{\mathcal{G}}^{''}$
 ci-dessus
 $$\I_{\mathcal{G}} = - \I_{\mathcal{G}}^{'} -
 (2 n + 5) \I_{\mathcal{G}}^{''} $$
 \end{prop}

 \subsection{Le cas $n$ impair}

  Ce cas est plus simple \`a traiter que le pr\'ec\'edent, gr\^ace
  \`a la proposition qui suit. On note $\mathcal{S}$ l'espace de
  modules $\mathcal{S}_{\epsilon}$ pour un param\`etre $\epsilon >
  0$ strictement inf\'erieur \`a la plus petite valeur critique, c'est-\`a-dire
  $1/2$.
  Soit $p$ l'entier tel que $2 p + 1 = n + 2$, et $p: \mathcal{C} \ri
 \P_5$ la conique universelle introduite au paragraphe
 pr\'ec\'edent. Soit $\O_{\mathcal{C}}(1)$ le fibr\'e inversible
 sur $\mathcal{C}$ image r\'eciproque du fibr\'e $\O(1)$ par la projection
 $q: \mathcal{C} \ri \P_2$.

  \begin{prop}
 Si $\mathcal{E}$ d\'esigne le fibr\'e
 vectoriel $p_* \left( \O_{\mathcal{C}}(p) \right)$ sur $\mathcal{C}$, on
 a $\mathcal{S} \simeq \mathrm{Grass} \left( 2 , \mathcal{E} \right)$.
 \end{prop}

 \begin{pv} La preuve s'appuie sur le lemme suivant.

 \begin{lem}
  Si $\Theta$ est un faisceau stable de support sch\'ematique
  une conique, de caract\'eristique d'Euler-Poincar\'e $n + 2$, et
  si $\Gamma \subset \H^0 \left( \Theta \right)$ est un pinceau de sections le
  syst\`eme coh\'erent $(\Gamma, \Theta)$ est $\epsilon-$stable.
 \end{lem}

 \begin{pv} Pour tout sous-faisceau coh\'erent $\Theta^{'} \subset
 \Theta$ de multiplicit\'e 1 on a par hypoth\`ese de stabilit\'e
 $\chi \left( \Theta^{'} \right) \leq \fc{n + 1}{2}$. Alors pour $\epsilon
 < 1/2$ et $\Gamma^{'} = \Gamma \cap \H^0 \left( \Theta^{'} \right)$
 l'in\'egalit\'e $\epsilon \hspace{0.1cm} \mathrm{dim} \left(\Gamma^{'} \right)
 + \chi \left( \Theta^{'} \right) \leq \fc{n + 2}{2} + \epsilon$ est clairement
 v\'erifi\'ee. \end{pv} \\

  Maintenant pour chaque conique $C$ du plan, le faisceau
  $\O_{C}(p)$ v\'erifie les hypoth\`eses du lemme. On a donc un
  morphisme modulaire $f: \mathrm{Grass} \left( 2 , \mathcal{E} \right)
  \ri \mathcal{S}$ qui est clairement injectif. Comme les deux
  vari\'et\'es ont m\^eme dimension, et que $\mathcal{S}$
  est une vari\'et\'e normale (cf th\'eor\`eme
  \ref{normaliteSalpha}), le morphisme $f$ est un isomorphisme.
  \end{pv} \\

  On note $\gamma$ la premi\`ere classe de Chern du sous-fibr\'e
  tautologique de rang 2 sur $\mathcal{S}$, et $\tau$ l'image
  r\'eciproque de la classe du fibr\'e $\O_{\mathbf{P}_5}(1)$ par le morphisme de
  projection $\pi: \mathcal{S} \ri \P_5$. On note enfin $\mathcal{D}$ le
  fibr\'e inversible "d\'eterminant" sur l'espace de modules
  $\mathcal{S}$ et $d = c_1 (\mathcal{D})$.

  \begin{lem}
 On a $d = (p-1) \tau - \gamma$.
  \end{lem}

 \begin{pv} Soit $\left( \mathbf{\Gamma} , \mathbf{\Theta} \right)$
 une famille universelle de syst\`emes coh\'erents param\'etr\'ee par
 $\mathcal{S}$. On peut supposer que le faisceau $\mathbf{\Theta}$ est
 $\O_{\mathcal{C}}(p)$. Alors si $l$ est une droite quelconque du plan, on a
 $$ \mathcal{D} \simeq \mathrm{det} \left( p_{!} (
 \mathbf{\Theta} \cdot \left[ \O_l \right] ) \right) \otimes
 \mathrm{det} \left( \mathbf{\Gamma}^* \right) \simeq \tau^{\otimes p-1}
 \otimes \mathrm{det} \left( \mathbf{\Gamma} \right)^* $$
 Ceci se voit en effet en utilisant la r\'esolution suivante du faisceau
 $\O_{\mathcal{C}}(p)$ comme $\O_{\mathbf{P}_5 \times \P_2}-$module
 $$ 0 \ri \O(-1, p-2) \ri \O(0,p) \ri \O_{\mathcal{C}}(p) \ri 0
 \hspace{2cm}(*)$$
 \end{pv} \\

  On doit donc pour finir calculer $\I = d^{2 n + 5} \cap \left[
  \mathcal{S} \right]$. On a
  $$ \begin{array}{lllcc}
     \I & = \sum_{0 \leq k \leq 2 n + 5} (-1)^k \hspace{0.1cm}
     \binom{2 n + 5}{k} \hspace{0.1cm} (p-1)^{2 n + 5 - k}
     \hspace{0.1cm} \tau^{2 n + 5 - k} \cdot \pi_* \left( \gamma^k \right)
     \cap \left[
     \P_5 \right] \\
     & = \sum_{0 \leq l \leq k \leq 2 n + 5} (-1)^k \hspace{0.1cm}
     \binom{2 n + 5}{k} \hspace{0.1cm} \binom{k}{l} \hspace{0.1cm}
     (p-1)^{2 n + 5 - k} \left| \begin{array}{llcc}
                                s_{k-l-n} (\mathcal{E}) & s_{l-n-1}
                                (\mathcal{E}) \\
                                s_{k-l-n+1} (\mathcal{E}) & s_{l-n}
                                (\mathcal{E}) \\
                                \end{array} \right| \cdot \tau^{2 n + 5 - k} \\
                                \end{array} $$
 De la suite exacte (*) on d\'eduit la r\'esolution suivante du
 fibr\'e $\mathcal{E}$, en posant $a = p (p-1) / 2$
 $$ 0 \ri \O_{\P_5}(-1)^{a} \ri \O_{\P_5}^{\binom{p+2}{2}} \ri \mathcal{E}
 \ri 0 $$
 ce qui donne $s (\mathcal{E}) = c(\mathcal{E}^{\vee})^{-1} =
 (1 + \tau)^a$. On obtient finalement
 $$ \fbox{$\displaystyle \I = \fc{1}{a + 1} \hspace{0.1cm} \sum_{0 \leq l \leq k \leq 2 n + 5} \hspace{0.1cm}
 (-1)^k \hspace{0.1cm} \binom{2 n + 5}{k} \hspace{0.1cm}
 \binom{k}{l} \hspace{0.1cm} (p-1)^{2 n + 5 - k} \hspace{0.1cm}
 (k - 2 l + 1) \hspace{0.1cm} \binom{a + 1}{k-l-n+1} \hspace{0.1cm}
 \binom{a + 1}{l-n} $} $$

 \subsection{Applications num\'eriques pour
 $n \leq 6$}

 Suivant les notations des deux derniers paragraphes, nous
 r\'ecapitulons les valeurs de $\I_{<}$ et
 $\Delta$ pour $n \leq 6$. Le symbole $\I_{<}$ d\'esigne
 $\I_{\mathcal{G}}$ si $n$ est pair, et $\I$ si $n$ est impair,
 c'est-\`a-dire le nombre d'intersection sur la grassmanienne
 relative; le symbole $\Delta$ ne concerne que $n=4,6$ et est
 \'egal au saut entre $\mathcal{S}_{1-}$ et
 $\mathcal{G}$. \\

 On obtient finalement le degr\'e cherch\'e pour $n$ pair par
 la formule deg $\overline{P_n } = \frac{1}{2} \left(
 \I_{<} + \Delta \right) + \sum_{\alpha > 0}
 \Delta_{\alpha}$, la seconde somme portant sur les valeurs
 critiques $\alpha > 0$; pour $n$ impair, en l'occurrence $n=5$
 on a deg $\overline{P_n } = \I_{<} + \sum_{\alpha > 0} \Delta_{\alpha}$.
 Les valeurs $\Delta_{\alpha}$ ont \'et\'e report\'ees dans
 la table du paragraphe \ref{appnumDelta}. On obtient ainsi
 le th\'eor\`eme \ref{theoprincipal}. Celui-ci sugg\`ere que le
 comportement de deg $\overline{P_n }$ d\'epend de la parit\'e de
 $n$. Nous ne pouvons cependant confirmer cette observation.

 \subsection{Remerciements}
 Cet article est une version remani\'ee et augment\'ee d'un
 chapitre de la th\`ese de l'auteur. La question du calcul du
 degr\'e de la vari\'et\'e des courbes de Poncelet m'avait \'et\'e
 propos\'ee par Joseph Le Potier, h\'elas disparu en D\'ecembre
 2005. Je lui dois de nombreux conseils et suggestions.

 \begin{center}
  \begin{table} \caption{Valeurs de $\Delta$ et $\I_{<}$
  pour $n \leq 6$}
   \begin{tabular}{|r|l|l|l|} \hline
   & $n=4$ & $n=5$ & $n=6$ \\
   \hline\hline
   $\Delta$ &  6558 & \hspace{0.35cm}$\times$ &  10278 \\
   \hline
   $\I_{<}$ & -6006 & - 5148 &  - 8580 \\
   \hline
   \end{tabular}
  \end{table}
  \end{center}

\end{document}